\documentclass[11pt]{ut-thesis}
\usepackage{my-thesis}
\degree{Doctor of Philosophy}
\department{Mathematics}
\gradyear{2020}
\author{Anne Dranowski}
\title{Comparing two perfect bases}
\begin{document}
%
\begin{preliminary}
\maketitle
%
\cleardoublepage
%
\begin{abstract}
We study a class of varieties which generalize the classical orbital varieties of Joseph. We show that our \textit{generalized} orbital varieties are the irreducible components of a Mirkovi\'c--Vybornov slice to a nilpotent orbit, and can be labeled by semistandard Young tableaux.

Furthermore, we prove that Mirkovi\'c--Vilonen cycles are obtained by applying the Mirkovi\'c--Vybornov isomorphism to generalized orbital varieties and taking a projective closure, refining Mirkovi\'c and Vybornov's result. 

As a consequence, we are able to use the Lusztig datum of a Mirkovi\'c--Vilonen cycle to determine the tableau labeling the generalized orbital variety which maps to it, and, hence, the ideal of the generalized orbital variety itself. By homogenizing we obtain equations for the cycle we started with, which is useful for computing various equivariant invariants such as equivariant multiplicity and multidegree. 

As an application, we show that in type \(A_5\) the Mirkovi\'c--Vilonen basis and Lusztig's dual semicanonical basis are not equal as perfect bases of the coordinate ring of the unipotent subgroup. 
This is significant because it shows that perfect bases are not unique.
%
Our comparison relies heavily on the theory of measures developed in \cite{baumann2019mirkovic}, so we include what we need. 

Finally, we state a conjectural combinatorial ``formula'' for the ideal of a generalized orbital variety in terms of its tableau. 
%
\end{abstract}
%
%
\begin{dedication}
To my best friend, Chester

\includegraphics[width=4.5cm]{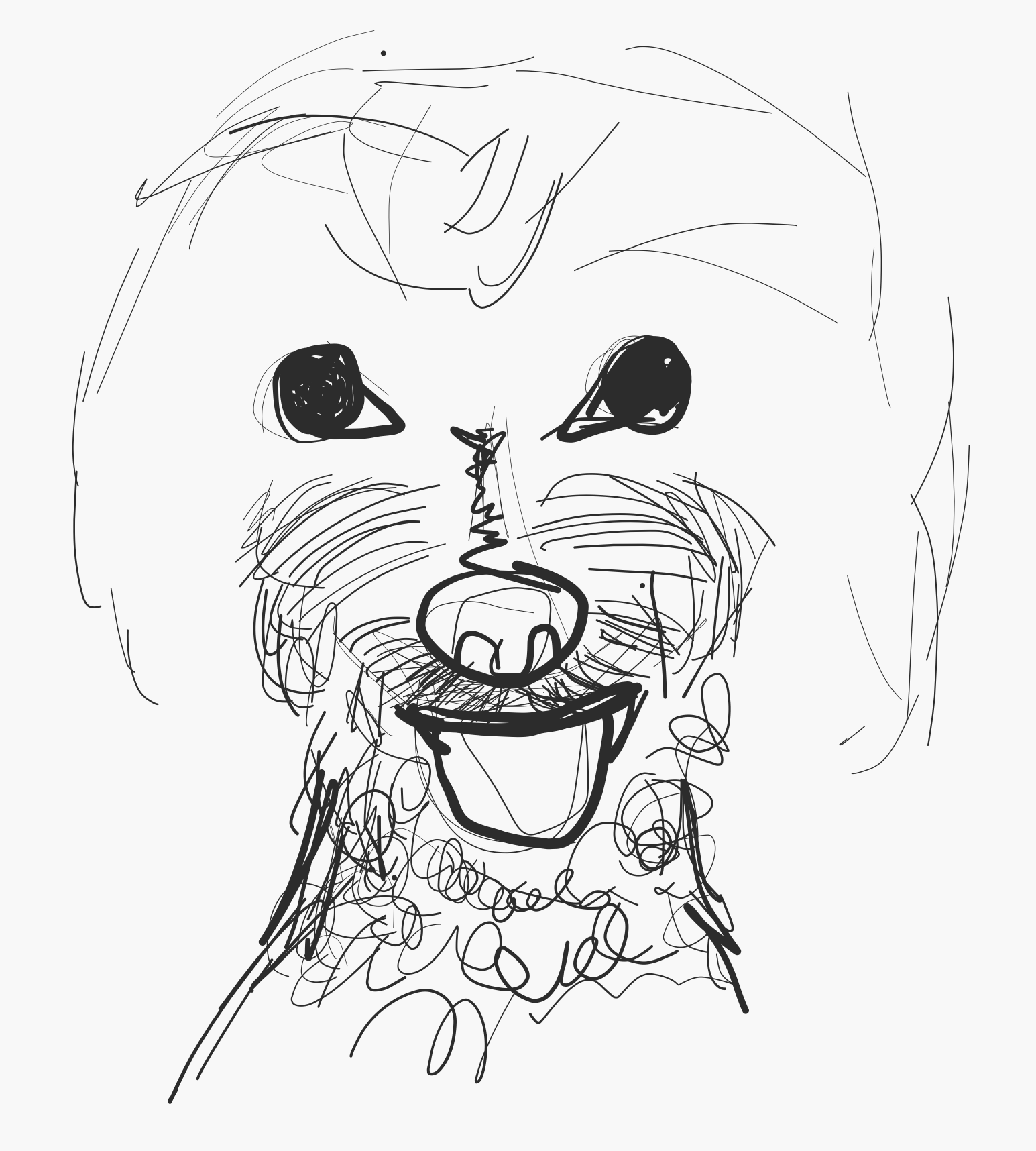}
\end{dedication}
%
\newpage  
%
%
\begin{acknowledgements}








Thank you to my advisor. For your guidance, for your patience, for your encouragement, and for generously sharing your (immense and open) knowledgebase. Needless to say none of this would have come to fruition without your immeasurable support. 

Thank you to Jemima. For all of your help, for your good humoured conversation, and for caring for our well-being. 
%
%

Thank you to my parents. Wealth of wisdom, well of love.

Thank you to my peers. For our discussions, and for your feedback. You make math a little bit easier. 

Thank you to all the awesome teachers at Dublin Heights EMS and William Lyon Mackenzie CI. For opening several doors, of various sizes.
Thank you to all the classmates whose good company and healthy competition has helped feed my reckless and relentless curiousity. 

Thanks are due to all the senior mathematicians who deigned to talk to me down the road:
Jim Arthur, Sami Assaf, Pierre Baumann, Alex Braverman, Almut Burchard, Greg Chambers, Ivan Cherednik, Ben Davison, Claus Fieker, Nora Ganter, Eugene Gorsky, Megumi Harada, David Hernandez, Nick Hoell, Abe Igelfeld, Yael Karshon, Allen Knutson, Anna Lachowska, Michael Lau, Samuel Leli\`evre, Michael McBreen, Fiona Murnaghan, Dinakar Muthiah, Ivan Mirkovi\'c, Adrian Nachmann, Vidit Nanda, Andrei Okounkov, Gerhard Pfister, Viviane Pons, Nick Proudfoot, Joe Repka, Ben Salisbury, Travis Scrimshaw, Mark Shimozono, Jason Siefken, Eric Sommers, Isabel Stenger, Bernd Sturmfels, Nicolas Thiery, Hugh Thomas, Peter Tingley, Monica Vazirani, Ben Webster, Paul Zinn-Justin. 
Even if just over email---thank you for taking the time to answer my questions and thanks for supporting me, one way or another.
%

Finally, thanks if you're reading this. Any mistakes are my own.
%
%
%
\end{acknowledgements}
%
\tableofcontents
%
\listoftables
%
%
%

\chapter*{Preface}
\label{ch:preface}

\epigraph{Curiosity killed the cat,\\
But satisfaction brought it back.}{James Newell Osterberg Jr.}

In their foundational paper \cite{mirkovic2007geometric}, Mirkovi\'c and Vilonen provided a geometric interpretation of the representation theory of reductive algebraic groups over arbitrary rings. Along the way, they constructed a basis of algebraic cycles 
for the cohomology of ``standard sheaves''
on the affine Grassmannian. 
While Ginzburg 
had already established an equivalence of categories of representations and sheaves over the complex numbers, the basis of the now eponymous MV cycles crucial to Mirkovi\'c and Vilonen's generalization was new. 
As such, it readily joined the ranks of a small but influential group of 
``good geometric bases'' of representations, including:
\begin{enumerate}
    \item Lusztig's canonical basis given by simple perverse sheaves on spaces of quiver representations \cite{lusztig1990canonical}
    \item Lusztig's semicanonical basis of constructible functions on Lusztig's nilpotent variety \cite{lusztig1992affine,lusztig2000semicanonical} and its dual basis of generic modules for the preprojective algebra studied by \cite{geiss2005semicanonical}
    \item Nakajima's basis of cycles in Borel--Moore homology of (Lagrangian subvarieties of) quiver varieties \cite{nakajima2001quiver} \label{it:nak}
\end{enumerate}
Much of the last decade in geometric representation theory has been devoted to the study of how these bases and their underlying spaces are related. 
The basic combinatorial structure that they are known to share is that of a crystal \cite{kashiwara1995crystal}. 
This includes restricting to bases for weight spaces and tensor product multiplicity spaces inside weight spaces.
%
The basic geometric structure that all but Nakajima's construction are known to share 
is that of a perfect basis: each of them behaves well with respect to actions of Chevalley generators, in addition to restriction and tensor product. 
Perfect bases were developed by Berenstein and Kazhdan in \cite{berenstein433geometric} in order to relate combinatorial crystals with actual bases.
In this thesis we will focus on comparing the basis of MV cycles and the dual semicanonical basis. Both give a perfect basis of the coordinate ring of the unipotent subgroup, 
%
whose crystal, \(B(\infty)\), is equal to that of the Verma module of highest weight zero.

The key to this comparison is the polytope construction \({\Pol}\) of \(B(\infty)\) conjectured by Anderson in \cite{anderson2003polytope} and supplied by Kamnitzer in \cite{kamnitzer2007crystal}. 
Any perfect basis comes with a god-given map \(\Pol\) to a certain set of polytopes (the MV polytopes). 
The two perfect bases that we are interested in are constructed from geometric spaces that come with intrinsically defined \(\Pol\) maps such that the triangle
\[
\begin{tikzcd}
    & B(\infty) \ar[dr,"{\Pol}"]\\
    \text{geometric spaces} \ar[ur,"b"] \ar[rr,"\Pol"] & & \text{polytopes}
\end{tikzcd}
\] 
commutes. Here, by a slight abuse of the notation, \(b\) denotes the construction at hand---we consider two. Similarly, because they have the same image, the remaining two maps are both denoted \(\Pol\). Of course, the two \(\Pol\) maps are not the same---for one thing they have different domains. 
The basis of MV cycles and the dual semicanonical basis were shown to be two such bases in \cite{baumann2014affine}. In particular, the moment polytope of an MV cycle \(Z\) agrees (up to sign) with the Harder--Narasimhan polytope of a generic preprojective algebra module \(M\) whenever \(Z\) and \(M\) define the same element in \(B(\infty)\).

To probe just how rich the ``polytope structure'' is, Baumann, Kamnitzer and Knutson, in \cite{baumann2019mirkovic}, associate measures to perfect (actual) basis elements, and ask whether equal polytopes have equal equivariant volumes.
\begin{question*}
    \label{q1}
    If \(M\) and \(Z\) define the same polytope, \(P\), do they define the same measure, \(D\), on \(P\)?
\end{question*}
%
We answer this negatively in Theorem~\ref{mthm:t3v1} below. 
For an element of the dual semicanonical basis, 
the measure \(D\) manifests 
as a piecewise polynomial measure whose coefficients are Euler characteristics of composition series in the corresponding module. Computing \(D\), given \(M\), is an exercise in basic hyperplane arrangements.

For an element of the basis of MV cycles, 
\(D\) can be reinterpreted as the Duistermaat--Heckman measure of the corresponding MV cycle.
This turns out to be related to a representation theoretic quantity called equivariant multiplicity (c.f.\ \cite[\S~4.2]{brion1997equivariant}). 
Given a local coordinate system, equivariant multiplicities are straightforward to compute. 
To find a local coordinate system in type \(A\), we rely on \cite{mirkovic2007geometric}.

Mirkovi\'c and Vybornov, in \cite{mirkovic2007geometric},
showed that certain subvarieties of the affine Grassmannian
of \(GL_m\) are isomorphic to slices
of \(N\times N\) nilpotent matrix orbits which we will temporarily denote \(\overline{\OO}\cap\TT\). 
Here \(N\) depends on the connected component of the affine Grassmannian that we start in. 
This isomorphism is what we will refer to as \textit{the} Mirkovi\'c--Vybornov isomorphism.
Mirkovi\'c and Vybornov also showed that the slices \(\overline{\OO}\cap\TT\) are isomorphic to affine Nakajima quiver varieties.
%
Their combined results provide a geometric framework for skew Howe duality, identifying
the natural basis of weight spaces in Nakajima's construction (``good basis'' no.~\ref{it:nak} above) with the natural basis of multiplicity spaces in skew Howe dual tensor products arising from
affine Grassmannians (c.f. Section~\ref{ss:sympl}).

In Chapter~\ref{ch:mvs}, we utilize the Mirkovi\'c--Vybornov's isomorphism to show that 
MV cycles can be obtained as projective closures
of irreducible components of a half-dimensional subvariety of \(\overline{\OO}\cap\TT\). We refer to these irreducible components as 
generalized orbital varieties. 

First, inspired by Spaltenstein's decomposition of the fixed point set of a flag variety under a regular unipotent element \cite{spaltenstein1976fixed},
we prove that generalized orbital varieties are in bijection with semistandard Young tableaux.
\begin{maintheorem}
\label{mthm:t1v1}
    (Theorem~\ref{thm:ad1})
    To each semistandard Young tableau \(\tau\) 
    we can associate a matrix variety whose unique top dimensional irreducible component is dense in a generalized orbital variety. Conversely, every generalized orbital variety is the Zariski closure of such a matrix variety.
    %
\end{maintheorem}
Interesting in its own right, this result lays the groundwork for a systematic description of ideals of generalized orbital varieties. Paired with the Mirkovi\'c--Vybornov isomorphism, and the fact that semistandard Young tableaux can be ``embedded'' in \(B(\infty)\), this result enables us to compute equivariant multiplicities of MV cycles.

\begin{maintheorem}
\label{mthm:t2v1}
    (Theorem~\ref{thm:ad2})
Let \(X_\tau\) be the generalized orbital variety corresponding to the semistandard Young tableau \(\tau\). 
Let \(\mathring Z_\tau\) be the image of \(X_\tau\) under the Mirkovi\'c--Vybornov isomorphism. 
Then the projective closure of \(\mathring Z_\tau \) is
the MV cycle whose polytope defines the same element of \(B(\infty)\) as \(\tau\). 
In other words, we have the following commutative diagram. 
\[
    \begin{tikzcd}
        \text{tableaux} \ar[r] \ar[d] & B(\infty) \ar[d] \\
        \text{generalized orbital varieties} \ar[r]  &\text{stable MV cycles}
    \end{tikzcd}
\] 
\end{maintheorem}

In Chapter~\ref{ch:calculs}, we present examples of equal basis vectors in types \(A_4\) and \(A_5\), and we conclude with an example of inequality in type \(A_5\).
This chapter is largely based on \cite[Appendix A]{baumann2019mirkovic} which was coauthored with J.\ Kamnitzer and C.\ Morton-Ferguson.
The calculations involved are carried out 
using the computer algebra systems SageMath, Macaulay2 and Singular.

\begin{maintheorem}
\label{mthm:t3v1}
    (Theorem~\ref{thm:adcmf})
There exists an MV cycle \(Z\) and a generic preprojective algebra module \(M\) such that \(Z\) and \(M\) have equal polytopes, but define different measures, and therefore distinct vectors in associated representations.
\end{maintheorem}

Theorem~\ref{mthm:t3v1} is related to the following general conjecture.
\begin{conjecture*}
    \cite[Conjecture 1.10]{baumann2019mirkovic}
    For any preprojective algebra module \(M\), 
    there exists a coherent sheaf \(\cF_M\) on the affine Grassmannian that is supported on a union of MV cycles such that 
    \[
        H^0(\cF_M\otimes\cO(n)) \cong H^\bullet(\mathbb G(M[t]/t^n))
    \] 
    is an isomorphism of torus representations.
    Here \(\mathbb G(M[t]/t^n)\) is a variety of submodules for the preprojective algebra tensored with \(\CC[t]\). 
\end{conjecture*}
An earlier version of this conjecture (due to Kamnitzer and Knutson) motivated a number of recent works by Mirkovi\'c and his coauthors on the subject of local spaces. In Chapter~\ref{ch:calculs} we explore a stronger version of the correspondence dictated by this conjecture called ``extra-compatibility''. It matches the cohomology of a flag variety of a preprojective algebra module and global sections of a corresponding MV cycle. 

\end{preliminary}
%
%
%

\chapter{Crystals, polytopes and bases}
\label{ch:bases}

\epigraph{It admits all deformation is justified by the search for the invariant.}{Ozenfant and Le Corbusier}

Let \(G\) be a complex algebraic simply-laced semisimple adjoint group of rank \(r\) and denote by \(\g\) its Lie algebra. 
A \(\g\)-crystal is a combinatorial object associated to a representation of \(\g\). 
The \(B(\infty)\) crystal is associated to the highest weight zero Verma module of \(\g\). 
The first construction of this crystal used Kashiwara's crystal basis for the positive part of the universal enveloping algebra of \(\g\) as the underlying set \cite{kashiwara1995crystal}.

In this chapter we review the dictionary of combinatorial models for \(B(\infty)\) that instructs and warrants subsequent geometric comparisons. 
Every concrete computation we are able to carry out starts with a ``Lusztig datum''.  
Traditionally, one uses these data to establish or compare crystal structure---something which we take for granted. We instead only use these data as crystal passports, to pinpoint ``corresponding'' elements in different models.

Towards defining abstract crystals, we begin with an elementary construction of Lusztig's canonical basis for the positive part of the universal enveloping algebra of \(\g\) (which is the same as Kashiwara's crystal basis).
By dualizing, we meet a ``perfect'' basis whose underlying crystal structure is again \(B(\infty)\).  
We then explain how to determine Lusztig data 
in each of the following models: MV polytopes and generic modules for preprojective algebras;
and semistandard Young tableaux when \(G\) is type \(A\).
%

\section{Notation}
\label{s:notation}
%
%
%
%
Let \(B\) be a Borel subgroup with unipotent radical \(U\) and let \(T\) be a maximal torus of \(B\). Denote by \(\fb\), \(\fu\) and \(\t\) their Lie algebras. 
Let \(P=\Hom(T,\CC^\times)\subset\t^\ast \) be the weight lattice of \(T\). 
%
%

Denote by \(\Delta \subset P\) the subset of roots (nonzero weights of \(T\) on \(\g\) with respect to the adjoint action) 
and by \(\Delta_+\subset\Delta\) the subset of positive roots (weights of \(T\) on \(\fu\) with respect to the adjoint action).
Denote by \(\{\alpha_i\}_I\subset\Delta_+\) the subset of simple roots (roots that cannot be expressed as the sum of any two positive roots). 
Since \(\Delta\) is a root system, the simple roots are a basis for \(\t^\ast\).
%
We denote by \(Q_+\) the positive root cone (nonnegative integer linear combinations of simple roots), and by \(Q\) the root lattice (integer linear combinations of simple roots). Note \(Q\subseteq P \). 
The positive root cone gives rise to the \new{dominance order} on \(P\). This is the partial order \(\lambda\ge\mu \) defined by the relation \(\{(\lambda,\mu)\in P\times P\big|\lambda - \mu \in Q_+\}\).

The simple reflection through the root hyperplane \(\Ker\alpha_i\) is denoted \(s_i\) and 
the set \(\{s_i\}_I\) generates the Weyl group \(W\) of \(G\) which acts on \(\t^\ast\) as well as \(\t\).
The latter action is defined by way of a \(W\)-invariant bilinear form \(\bil{\phantom{0}}:P\times P\to\QQ\). Given \(\alpha\in\Delta\) its coroot \(\alpha^\vee\in\t\) is defined by the requirement that \(\lambda(\alpha^\vee) = 2\frac{\bil[\lambda]{\alpha}}{\bil[\alpha]{\alpha}}\) for all \(\lambda\in P\). 
In particular, the simple coroots \(\{\alpha^\vee_i\}_I\) form a basis for \(\t\). 
Furthermore, the \((i,j)\)th entry of the Cartan matrix of \(\g\) is given by \(\alpha_j(\alpha_i^\vee) \). 
Indiscriminately, we sometimes also write \(\langle\alpha^\vee,\lambda\rangle\) for the pairing \(\lambda(\alpha^\vee)\). 
We denote by \(\Delta^\vee\) the set of coroots \(\{\alpha^\vee \big|\alpha\in \Delta\}\). 
The pairing 
\(\langle\phantom{0},\phantom{0}\rangle : \Delta^\vee\times P \to \ZZ\) 
extends linearly to a pairing on \(\t\times\t^\ast\) which is perfect. 
In the simply-laced case we will take \(\bil{\phantom{0}}\) to be the dot product on the ambient vector space \(\t^\ast\cong\CC^r\) so that \((\alpha,\alpha) = 2\) and \(\lambda(\alpha^\vee) = (\alpha,\lambda)\) for all \(\alpha\in\Delta\). 

The \new{fundamental weights} \(\{\omega_i\}_I\) are defined by the system \(\{\omega_i(\alpha_j^\vee) = \delta_{i,j}\}_{I\times I}\) forming a basis for \(\t^\ast\) dual to the basis of simple coroots for \(\t\).
%
%
%
If \(\lambda(\alpha_i^\vee)\ge 0\) for all \(i\in I\) then \(\lambda \) is called \new{dominant}, and 
the subset of all such weights is denoted \(P_+\subset P \). It parametrizes finite-dimensional highest weight irreducible representations of \(G\). 

The coweight lattice is denoted \(P^\vee\) and defined to be the set of coweights \(\xi\) such that \(\lambda(\xi)\in\ZZ\) for all \(\lambda\in P\). 
Equivalently, \(P^\vee = \Hom(\CC^\times, T)\). 
The fundamental coweights, \(\{\omega_i^\vee\}_I\),  are likewise defined as the basis for \(\t\) that is dual to the basis of simple coroots for \(\t^\ast\). 
%

Given \(\lambda\in P_+\) we denote the irreducible representation of \(G\) having highest weight \(\lambda\) by \(L(\lambda)\).
Rather than work with individual representations, we will work in the coordinate ring of the unipotent subgroup, \(\CC[U]\). 
Fix a highest weight vector \(v_\lambda\in L(\lambda)\) and denote by \(v_\lambda^* : L(\lambda) \to \CC\) the linear form such that \(v^\ast_\lambda(v_\lambda) = 1\) and \(v^\ast_\lambda(v) = 0\) for any weight vector \(v\in L(\lambda)\) other than \(v_\lambda\). 
Berenstein and Kazhdan \cite{berenstein433geometric} define a family of inclusions \(\Psi_\lambda : L(\lambda) \to \CC[U]\) by requiring that \(\Psi_\lambda(v)(u) = v_\lambda^\ast(u\cdot v)\) for any \(u\in U\). Thus, \(\Psi(L(\lambda)_\mu)\subset \CC[U]_{-\lambda + \mu}\).
In \cite[Corollary 5.43]{berenstein433geometric} they show that \(\Psi_\lambda\) is a map of \(\fu\)-representations and that
\begin{equation}
  \label{eq:CUallreps}
  \CC[U] = \bigcup_{P_+}\Psi_\lambda(L(\lambda))\,.
\end{equation}

Let \(e_i,f_i\in\g\) be root vectors of weights \(\alpha_i,-\alpha_i\), respectively, such that \([e_i,f_i] = \alpha_i^\vee\). 
The action of \(\fu \) on \(\CC[U]\) by left invariant vector fields 
\(e_i\cdot f(u) = \frac{d}{dt}\big|_{t=0}f(u\exp(te_i))\)
gives rise to the perfect pairing 
\begin{equation}
  \label{eq:perfpair}
  \begin{split}
  \cU(\fu)\times \CC[U] &\to \CC \\
  (a,f) &\mapsto \langle a,f\rangle := (a\cdot f)(1_U)\,.
  \end{split}
\end{equation}
Here \(\cU(\fu)\) denotes the universal enveloping algebra of \(\fu\). Generated by the root vectors \(e_i\), it is graded by \(Q_+\), with \(\deg e_i = \alpha_i\) for any \(i\in I\). 
The conjugation action of \(T\) on \(U\) makes \(\CC[U]\) into a \(Q_+\) graded algebra, too, and the pairing defined by Equation~\ref{eq:perfpair} induces isomorphisms 
\begin{equation}
  \label{eq:udual}
  \CC[U]_{-\nu}\cong \cU(\fu)_\nu^\ast 
\end{equation}
for all \(\nu\in Q_+\). 

We denote by \(\bar s_i \) the lift \(\exp(-e_i)\exp(f_i)\exp(-e_i)\) of \(s_i\in W\) to \(G\) for any \(i\in I\). 
By \cite[Proposition 3]{tits1966constantes}, the set \(\{\bar s_i\}_I\) satisfies the braid relations, so the lift \(\bar w\) makes sense for any \(w\in W\). 
Fix a longest element \(w_0\in W\) and write \(\ell\) for its length. Note \(\ell = \# \Delta_+ = \dim U\). 
We call \(\uvi = (i_1,\dots,i_\ell)\) a \new{reduced word} if \(w_0 = s_{i_1}\cdots s_{i_\ell}\) is a reduced expression. 
Reduced words induce convex orders on $\Delta_+$ where \new{convex order} means that the sum of any two positive roots \(\alpha + \beta \) falls between \(\alpha \) and \(\beta\). 
To a reduced word \(\uvi\) and a number \(0\le k\le\ell\), we can associate:
\begin{itemize}
  \item the expression \(s_{i_1}\cdots s_{i_{k}}\), denoted \(w_k^{\uvi} \), and
  \item the positive root \(w_{k-1}^{\uvi} \alpha_{i_k}\), denoted \(\beta_k^{\uvi} \).
\end{itemize}
Not to be confused with the longest element, the value of \(w^{\uvi}_k\) when \(k = 0\) is defined to be the identity element \(e\in W\). 
With
this notation the convex order on \(\Delta_+ \) coming from \(\uvi \) is the enumeration
\begin{equation}
  \label{eq:convex}
  \beta_1^{\uvi} \le \beta_2^{\uvi}\le\cdots \le\beta_\ell^{\uvi}\,.
\end{equation}
%
%
\begin{example}
  \label{eg:betanum}
  Let \(\g = \sl_3\) and take \(\uvi = (1,2,1)\) then 
  \[
  \beta_1 = \alpha_1 \le \beta_2 = s_1\alpha_2 = \alpha_1 + \alpha_2 \le \beta_3 = s_1s_2\alpha_1 = s_1(\alpha_1 + \alpha_2) = \alpha_2 \,. 
  \]
\end{example}
As a consequence of Equation~\ref{eq:udual}, we have that \(\CC[U]\cong\cU(\fu)^\ast\) as \(\fu\)-modules, making \(\CC[U]\) a representation which has the \(B(\infty)\) ``crystal structure''.
%
\begin{definition}
  \label{def:crystal}
\cite{kashiwara1995crystal}
  The data of a combinatorial \(\g\)-\new{crystal} is a set \(B\) endowed with a weight map \(\wt : B\to P \) and a family of maps
\begin{equation}
  \label{eq:crystalops}
  \tilde e_i,\, \tilde f_i : B \to B\sqcup \{0\} 
  \qquad\qquad 
  \varepsilon_i ,\, \varphi_i : B\to\ZZ
\end{equation} 
for all \(i\in I\), such that \(\langle\alpha_i^\vee,\wt(b)\rangle = \varphi_i(b) - \varepsilon_i(b) \) for all \(b \in B\). 
The element 0 is added so that \(\tilde e_i\) and \(\tilde f_i\) are everywhere defined.
Otherwise, the operators \(\tilde e_i,\, \tilde f_i\) are mutually inverse partial bijections such that 
\begin{equation}
  \label{eq:efid}
  b'' = \tilde e_i b' \Longleftrightarrow \tilde f_i b'' = b'
\end{equation}
and then 
\begin{equation}
  \label{eq:wtsofcrystalops}
  \wt (b'') = \wt(b') + \alpha_i 
  \qquad 
  \varepsilon_i(b'') = \varepsilon_i(b') - 1 
  \qquad
  \varphi_i(b'') = \varphi_i(b') + 1 
\end{equation}
\(B\) is called \new{upper semi-normal} if, for all \(b\in B\) and for all \(i\in I\), there exists \(n\in\NN\) such that \(\tilde e_i^n b = 0\) and \(\varepsilon_i(b) = \max\{k\big|\tilde e_i^k b \ne 0\}\). In other words, \(p = \varepsilon_i(b)\) is the largest integer such that  \(\tilde e_i^p b\) is defined, and, for all \(n\ge p + 1\), \(\tilde e_i^n b = 0\). 
\end{definition}

\section{Lusztig data as an artifact of quantization} 
\label{s:ldviaUq}
As the section title suggests, quantum groups play an auxiliary role for us: they are used to define a Lusztig datum of an element of \(B(\infty)\). 
According to Lusztig \cite{lusztig2010introduction}, ``the theory of quantum groups is what led to extremely rigid structure, in which the objects of the theory are provided with canonical bases with rather remarkable properties'' specializing for \(q = 1\) to canonical bases for objects in the classical theory.
Our exposition is based on~\cite{tingley2017elementary}. 

Let \(\cU_q(\g)\) be the quantized universal enveloping algebra of \(\g\). 
Its standard generators over \(\CC(q)\) are denoted \(\{E_i,F_i,K_i^{\pm 1}\}_I\) in analogy with the root vectors \(\{e_i,f_i,\alpha_i^\vee\}_I\) generating the classical enveloping algebra, which can be obtained from \(\cU_q(g)\) by taking \(q = 1\).

Two elements of the index set, \(i,j\in I\), are said to be \new{connected} if they are joined by an edge in the Dynkin diagram of \(\g\).
The ``upper triangular'' \(\CC(q)\)-subalgebra of \(\cU_q(\g)\) is generated by \(\{E_i\}_I\) subject to the following relations. 
\begin{equation}
  \label{eq:Uqrel}
  E_i^2 E_j + E_j E_i^2 = (q + q^{-1}) E_i E_j E_i  \text{ if } i, j\text{ are connected, and } 
  E_i E_j = E_j E_i \text{ otherwise.} 
\end{equation}
We denote this algebra \(\cU^+_q(\g) \). The abstract braid group is employed in the very elementary construction of a family of PBW bases for \(\cU^+_q(\g) \), one for each choice of reduced word. 

Let \(\langle T_i\rangle_I\) be a presentation of the abstract braid group on \(\# I\) strands. 
and define 
\begin{equation}
  \label{eq:braidaction}
  T_i\cdot E_j = \begin{cases}
    E_j & i\ne j\text{ not connected } \\
    E_iE_j - q^{-1}E_jE_i & i,j\text{ connected}\\
    -F_j K_j& i = j \,.
  \end{cases}
\end{equation}
Fix a reduced word \(\uvi = (i_1,\dots,i_\ell)\) and consider 
\begin{equation}
  \label{eq:braidbasisgens}
  \begin{split}
    E^{\uvi}_{1} &:= E_{i_1} \\ 
    E^{\uvi}_{2} &:= T_{i_1} \cdot E_{i_2} \\
    E^{\uvi}_{3} &:= T_{i_1} T_{i_2} \cdot E_{i_3} \\
    &\,\,\,\,\,\vdots\,.
  \end{split}
\end{equation}
Each \(E^{\uvi}_{k}\) is a weight vector in \(\cU_q(\g)\) of weight \(\beta_k^{\uvi}\).

For \(n \ge 1\) we denote by \([n]\) the quantum integer \(q^{n-1} + q^{n-3} + \cdots + q^{-n+3} + q^{-n+1} \) and by \([n]!\) the quantum factorial \([n][n-1]\cdots[2][1]\). Note, when \(q = 1 \), \([n] = n \) and \([n]! = n!\). 
Consider the set 
\begin{equation}
  \label{eq:ipbw}
  B_{\uvi} = \left\{{E^{\uvi}_{1}}^{(n_1)} {E^{\uvi}_{2}}^{(n_2)}\cdots {E^{\uvi}_{\ell}}^{(n_\ell)} \big| n_1,\dots,n_\ell\in\NN\right\}\,,
\end{equation}
with \(X^{(n)}\) denoting the divided power \(X^n/[n]!\). Given \(n\in\NN^{\ell}\) we will sometimes denote by \(E^{\uvi}(n)\) the element 
\( {E^{\uvi}_{1}}^{(n_1)} {E^{\uvi}_{2}}^{(n_2)}\cdots {E^{\uvi}_{\ell}}^{(n_\ell)} \). 
%
\begin{example}
  \label{eg:sl3basis}
Let
\(\g = \sl_3\) and take \(\uvi = (1,2,1) \) and \(\uvi' = (2,1,2)\). Then 
\begin{align*}
  (E^{\uvi}_{1}, E^{\uvi}_{2}, E^{\uvi}_{3}) &= (E_1, E_1 E_2 - q^{-1} E_2 E_1, E_2 )\\
  (E^{\uvi'}_{1}, E^{\uvi'}_{2}, E^{\uvi'}_{3}) &= (E_2, E_2 E_1 - q^{-1} E_1 E_2, E_1 )
\end{align*}
and
\((\beta^{\uvi}_1, \beta^{\uvi}_2, \beta^{\uvi}_3) = (\alpha_1, \alpha_1 + \alpha_2, \alpha_2 )\) and 
\((\beta^{\uvi'}_1, \beta^{\uvi'}_2, \beta^{\uvi'}_3) = (\alpha_2, \alpha_2 + \alpha_1, \alpha_1)\).
\end{example}
\begin{theorem}\label{thm:isbasis}
\cite[Propositions 1.8 and 1.13]{lusztig1990finite}
\(B_{\uvi}\) is a \(\CC(q)\)-linear PBW basis of \(\cU^+_q(\g)\) for any choice of reduced word \(\uvi\). 
\end{theorem}
%
%
The exponent \(n = (n_1,\dots,n_\ell)\) on the element \(E^{\uvi}(n)\) of \(B_{\uvi}\) is called its \(\uvi\)-\new{Lusztig datum} and we say that \(n\) has weight \(\nu\) if \(\nu = \sum_k n_k \beta^{\uvi}_k \). 
%

The PBW bases in general do not have crystal structure, but can be used to define a basis which does, whose elements can be represented uniquely by elements of \(B_{\uvi}\) for every \(\uvi\) as follows.
\begin{theorem}
  \label{thm:welldefofcb}
\cite[Proposition 2.3(b)--(c)]{lusztig1990canonical}
  The \(\ZZ[q^{-1}]\)-submodule of \(\cU_q^+(\g)\) generated by the basis \(B_{\uvi}\) does not depend on the choice of \(\uvi\). It is denoted \(\cL\). 
  Similarly, the set \(B_{\uvi} + q^{-1}\cL \) is a \(\ZZ\)-basis of the quotient \(\cL/q^{-1}\cL\) independent of \(\uvi\). 
\end{theorem}
The extremely rigid structure underlying the family \(\{B_{\uvi} + q^{-1} \cL \}\) is uniquely specified if we add the requirement that it be invariant with respect to an algebra antihomomorphism. We fix the bar antihomomorphism \(\bar{\phantom{0}} : \cU_q^+\to\cU_q^+\) defined by \(\bar q = q^{-1}\) and \(\bar E_i = E_i \). 

For each \(n\in\NN^{\ell}\), there is a unique bar-invariant element 
\begin{equation}
  \label{eq:barinvar} 
  b^{\uvi}(n) = \sum_{m\in\NN^\ell} c(n,m) E^{\uvi}(m)
\end{equation}
with \(c(n,n) = 1\) and \(c(n,m) \in q^{-1} \ZZ[q^{-1}]\) for any \(m\ne n\). Together the elements \(\{b^{\uvi}(n)\}_{n\in\NN^{\ell}}\) form \new{Lusztig's canonical basis} \(\cB \) of \(\cU_q^+(\g) \). 
Remarkably, it does not depend on the choice of \(\uvi\). 
\begin{theorem}
  \label{thm:cbexists}
  \cite[cf.\ Theorem 3.2]{lusztig1990canonical}
\(\cB \) is the unique basis of \(\cU^+_q(\g)\) such that
  \begin{enumerate}[label=(\roman*)]
    \item \(\cB\subset\cL\) and \(\cB + q^{-1}\cL\) is a \(\ZZ\)-basis of \(\cL/q^{-1}\cL\) agreeing with \(B_{\uvi} + q^{-1}\cL\) for any \(\uvi\). \label{thm:cbexistsi}
    \item \(\cB\) is pointwise bar invariant. 
  \end{enumerate}
\end{theorem}
%
This theorem provides a family of inverse bijections to the defining parametrization of \(\cB\) in Equation~\ref{eq:barinvar}: 
the \(\uvi\)-Lusztig datum of a canonical basis element, \(b\mapsto n^{\uvi}(b) \), where \(n=n^{\uvi}(b)\) is such that \(E^{\uvi}(n) \equiv b\bmod q^{-1}\cL\).

Specializing \(\cB\) at \(q = 1\) gives a basis of \(\cU(\fu)\) whose dual basis is a basis of \(\CC[U]\) having enough nice properties to merit a definition (c.f.\ \cite[Theorems 1.6, 7.5]{lusztig1990canonical}). 
But first, by a slight abuse of notation, let us denote this dual basis \(B(\infty)\) because it is the prototypical construction and our landmark model 
of the \(B(\infty)\) crystal structure.
%
%
\begin{definition}
  \label{def:perfb}
  \cite[Definition 5.30]{berenstein433geometric}
  $B \subset \CC[U]$ is called a \new{perfect basis} if it is a \(\CC\)-linear basis of \(\CC[U]\) endowed with an upper semi-normal crystal structure and having the following properties.
\begin{itemize}
  \item The constant function \(1_U\) belongs to \(B\).
  \item Each \(b\in B\) is homogeneous of degree \(\wt(b)\) with respect to the \(Q_+\) grading of \(\CC[U]\).
  \item For each \(i\in I\) and for each each \(b\in B\), the expansion of \(e_i\cdot b\) in the basis \(B\) has the form 
\begin{equation}
  \label{eq:perfeiop}
  e_i\cdot b = \varepsilon_i(b) \tilde e_i(b) + \Sp_{\CC}(\{b'\in B\big|\varepsilon_i(b') < \varepsilon_i(b) - 1\})\,.
\end{equation}
\end{itemize}
\end{definition}
%
%
Indeed, \(B(\infty)\) is a pefect basis of \(\CC[U]\), and the following theorem allows us to define (a family of) bijections \(B\to B(\infty) \to\NN^\ell\) for any perfect basis \(B\) of \(\CC[U]\). 
\begin{theorem}
  \label{thm:perfbunique}
  \cite[Theorem~5.37]{berenstein433geometric}
Let \(B\) and $B'$ be two perfect bases of $\CC[U]$. Then there is a unique
bijection $B\cong B'$ that respects the crystal structure.
\end{theorem}
%
%
One consequence of this thesis is that perfect bases are not unique, and Theorem~\ref{thm:perfbunique} justifies our comparison of two particular perfect bases on combinatorial grounds.

\section{Three models for \texorpdfstring{\(B(\infty)\)}{B-infinity}}
\label{s:bmodels}
Here we describe three models for the \(B(\infty)\) crystal structure (underlying the \(B(\infty)\) dual canonical basis); two in the simply-laced case, and one in type \(A\). Each comes with its own notion of Lusztig data (the southeasterly arrow below)
\begin{equation}
  \label{eq:bmodel}
\begin{tikzcd}
  \text{model}\ar[dr,"n^{\uvi}"']\ar[rr,dashed,"b"] & & B(\infty)\ar[dl,"n^{\uvi}"] \\
  & \NN^\ell & 
\end{tikzcd}  
\end{equation}
and the requirement that the above triangle commute for every reduced word \(\uvi\)
uniquely determines the horizontal arrow. 
%
\subsection{MV polytopes}
\label{ss:mvps}
We begin with the MV polytopes. This part of the exposition is based on \cite{kamnitzer2007crystal,kamnitzer2010mirkovic}. 

Let \(\Gamma\) denote the set of Weyl orbits of fundamental coweights, aka \new{chamber coweights},
and let \(M_\bullet = \{M_\gamma\}_\Gamma\) be a collection of integers. 
%
%
Recall that \(a_{ij} = \alpha_j(\alpha_i^\vee)\) are the entries of the Cartan matrix of \(\g\).

\(M_\bullet \) is said to satisfy the \new{edge inequalities} if 
\begin{equation}
  \label{eq:eiwi}
  M_{w\omega_i^\vee} + M_{ws_i\omega_i^\vee} + \sum_{j\ne i} a_{ji} M_{w\omega_j^\vee} \le 0
\end{equation}
for all \(i\in I\) and for all \(w\in W\). 

For completeness we recall here a convenient characterization of the Bruhat order on \(W\). Let \(u,w\in W\) be arbitrary, and fix a reduced expression \(w = s_{i_1}\cdots s_{i_b}\). Then \(u\le w\) if and only if there exists a reduced expression \(u = s_{j_1}\cdots s_{j_a}\) such that \((j_1,\dots,j_a)\) is a subword of \((i_1,\dots,i_b)\). 

Let \((w,i,j)\in W\times I\times I\) be such that \(w \le ws_i \), \(w\le ws_j\) in the Bruhat order, and suppose \(i\ne j\). 
\(M_\bullet\) is said to satisfy the \new{tropical Pl\"ucker relation} at \((w,i,j)\) if \(a_{ij} = 0 \), or \(a_{ij} = a_{ji} = -1\) and 
\begin{equation}
  M_{ws_i\omega_i^\vee} + M_{ws_j\omega_j^\vee} = \min\{M_{w\omega_i^\vee} + M_{ws_is_j\omega_j^\vee}, M_{ws_js_i\omega_i^\vee} + M_{w\omega_j^\vee}\}\,.
\end{equation}
The edge inequalities and tropical Pl\"ucker relations are due to Kamnitzer and can be found in \cite{kamnitzer2010mirkovic}, which covers groups that aren't simply-laced as well.

\begin{definition}
  \label{def:bzdata}
  \(M_\bullet\) is called a BZ (Berenstein--Zelevinsky) datum of weight \((\lambda,\mu)\) if
  \begin{enumerate}[label=(\roman*)]
    \item \(M_\bullet\) satisfies the tropical Pl\"ucker relations,
    \item \(M_\bullet\) satisfies the edge inequalities,
    \item \(M_{\omega_i^\vee} = \langle \mu,\omega_i^\vee\rangle \) and \(M_{w_0\omega_i^\vee} = \langle \lambda,w_0\omega_i^\vee\rangle\) for all \(i\).
  \end{enumerate}
\end{definition}

BZ data define MV polytopes in \(\t^\ast_{\RR} = P\otimes\RR\)
as follows. 
Given a BZ datum \(M_\bullet = (M_\gamma)_{\Gamma}\) of weight \((\lambda,\mu)\) we define for each \(w\in W\) the weight \(\mu_w \in \t_{\RR}^\ast\) by requiring that \(\langle \mu_w, w\omega_i^\vee\rangle = M_{w\omega_i^\vee}\) for each \(i\in I\). 
Then 
\begin{equation}
  \label{eq:polbz}
      P(M_\bullet) = \{\nu\in\t_{\RR}^\ast \big| \langle\nu,\gamma\rangle\ge M_\gamma \text{ for all }\gamma\in\Gamma\} 
      = \Conv(\{\mu_w \big| w\in W\})
\end{equation}
is an \new{MV polytope} of weight \((\lambda,\mu)\). Here \(\Conv\) denotes convex hull.
If \(\mu\le\lambda\) in the dominance order
then we say that \(P(M_\bullet) \) has lowest vertex \(\mu_e = \mu\) and highest vertex \(\mu_{w_0} = \lambda\). Otherwise \(P(M_\bullet)\) is empty.
\begin{example}
  \label{eg:sl2polbz}
Let \(G = PGL_2\) and \(W = \langle s\rangle = \{e,s\}\), 
Then we can represent a coweight \((a,b)\bmod (1,1)\in P^\vee\) by the number \(a-b\). In particular, the chamber coweights are \(\omega^\vee = 1\) and \(s\omega^\vee = -1\).
%
%
The tuple \(M_\bullet = (M_{\omega^\vee}, M_{s\omega^\vee}) \) 
satisfies the edge inequalities if 
\[
\begin{cases}
  M_{\omega^\vee} + M_{s\omega^\vee} = M_{\omega^\vee} + M_{s\omega^\vee} \le 0 & w = e \\
  M_{s\omega^\vee} + M_{\omega^\vee} = M_{s\omega^\vee} + M_{\omega^\vee} \le 0 & w = s\,. 
\end{cases}  
\]
The tropical Pl\"ucker relations hold automatically, so
\(M_\bullet\) makes sense as a BZ datum of weight \((\lambda,\mu)\) if and only if \(M_{\omega^\vee} = \mu\), \(M_{s\omega^\vee} = -\lambda\) and \(\mu\le \lambda\).
The corresponding MV polytope 
has vertices \(\mu_e = \mu\) and \(\mu_s = \lambda\).
\[
  P(M_\bullet) = \{\nu\big|\langle\nu,\omega^\vee\rangle \ge M_{\omega^\vee}, \langle\nu,s\omega^\vee\rangle \ge M_{s\omega^\vee}\} = \{\nu\big| \mu\le\nu\le\lambda\}\,. 
\]
\end{example}
%
The \(\uvi\)-Lusztig datum of an MV polytope \(P(M_\bullet)\)
is the tuple \(n^{\uvi}(M_\bullet)\in\NN^\ell\) defined by 
\begin{equation}
  \label{eq:bzld}
  n_k^{\uvi}(M_\bullet) = -M_{w^{\uvi}_{k-1}\omega^\vee_{i_k}} - M_{w^{\uvi}_k\omega^\vee_{i_k}} - \sum_{j\ne i} a_{ji} M_{w^{\uvi}_k \omega^\vee_j}\,.
\end{equation}
%
The weights \(\{\mu_k^{\uvi} := \mu_{w_{k}^{\uvi}}\}_{0\le k\le\ell}\) determine the path 
\(\mu_0^{\uvi} = \mu \to\mu_2^{\uvi}\to\cdots\to\mu_{\ell}=\lambda\) 
in the 1-skeleton of \(P(M_\bullet)\) 
and one can show that 
\(n_k\) is the length of the \(k\)th leg of the path, 
\begin{equation}
  \label{eq:mvleg}
  \mu^{\uvi}_k - \mu^{\uvi}_{k-1} = n_k \beta_k^{\uvi} \,. 
\end{equation}
Note that under our assumptions on \(G\), \(Q = P\) and so the MV polytopes lie in \(Q\), and can be translated to lie in \(Q_+\). 
%
\begin{example}
  \label{eg:sl2mvleg}
  Continuing with \(G = PGL_2\), we have no choice but to take \(\uvi = (1)\), and then \(n_1 = - M_{\omega^\vee} - M_{s\omega^\vee} = -\mu - (-\lambda) = \lambda - \mu\). 
  Indeed, this is also the length of the only edge of \(P(M_\bullet)\)---itself!
\end{example}
%
  An MV polytope of weight \((\nu,0)\), i.e.\ one whose lowest vertex \(\mu_e = 0\), is called \new{stable} and the collection of all such polytopes is denoted \(\cP\). 
%
\begin{example}
  \label{eg:sl3stabmv}
Let \(G = PGL_3\). Consider the stable MV polytope \(S\) of weight \((\alpha_1, 0)\). It has 
\[
\begin{aligned}
  M_{\omega^\vee_1} &= 0 & M_{w_0\omega^\vee_1} &= 0 \\
  M_{\omega^\vee_2} &= 0 & M_{w_0\omega^\vee_2} &= -1 \,.
\end{aligned}  
\]
Since \(S\) is just a line segment, we can infer from the conditions 
\[
  \langle\nu,s_1\omega^\vee_1\rangle\ge M_{s_1\omega^\vee_1} 
  \text{ and }
  \langle\nu,s_2\omega^\vee_2\rangle\ge M_{s_2\omega^\vee_2} 
\] 
on \(\nu\in P\) that \(M_{s_1\omega^\vee_1} = -1 \text{ and } M_{s_2\omega^\vee_2} = 0\).
Note, 
the (this time nontrivial) tropical Pl\"ucker relation 
\[
 -1 + 0 = M_{s_1\omega^\vee_1} + M_{s_2\omega^\vee_2} = \min (M_{\omega^\vee_1} + M_{s_1s_2\omega^\vee_2}, M_{\omega^\vee_2} + M_{s_2s_1\omega^\vee_1}) = \min(0-1,0+0)
\]
is satisfied.
\end{example}
%
The following theorem of Kamnitzer characterizes the set of all stable MV polytopes according to their Lusztig data.
\begin{theorem}
  \label{thm:ldpol}
  \cite[Theorem 7.1]{kamnitzer2010mirkovic}
  Taking \(\uvi\)-Lusztig datum is a bijection \(\cP\to \NN^\ell\). 
\end{theorem}
\begin{example}
  \label{eg:pol11}
  Once again take \(G = PGL_3\).
  The stable MV polytopes of weight \((\alpha_1 + \alpha_2,0)\) are depicted below.
  $$
  \begin{tikzpicture}[scale = 0.5]
    \draw (0,0) -- (1,1) node[right] {$\alpha_1$} -- (0,2) node[above] {$\alpha_1 + \alpha_2$} -- (0,0); 
    \draw (7,0) -- (6,1) node[left] {$\alpha_2$} -- (7,2) node[above] {$\alpha_1 + \alpha_2$} -- (7,0); 
  \end{tikzpicture}
  $$ 
  With respect to the reduced word \(\uvi = (1,2,1)\) their Lusztig data are \((1,0,1)\) and \((0,1,0)\) respectively. 
\end{example}
\subsection{The dual semicanonical basis}
\label{ss:dsc}
Let \((I,H)\) denote the double of the simply laced Dynkin quiver of \(G\). Given \(h= (i,j)\in H\), let \(\bar h\) denote the oppositely oriented edge \((j,i)\).
To each edge \(h \) we attach a ``sign'' \(\epsilon(h)\in\{1,-1\}\) such that \(\epsilon(h) + \epsilon(\bar h) = 0 \). 
\begin{definition}
  \label{def:preprojalg}
  The preprojective algebra \(\cA\) is the quotient of the path algebra of \((I,H)\) by the relation 
  \begin{equation}
    \label{eq:preprojrel}
    \sum_{h\in H} \epsilon(h) h\bar h = 0 \,.
  \end{equation}
\end{definition}
Thus, an \(\cA\)-module is an \(I\)-graded vector space \(M = \bigoplus_{i\in I} M_i\) together with linear maps
  \[
    \{M_{(i,j)}: M_i\to M_j\}_{(i,j)\in H}
  \]
  such that 
\begin{equation}
  \label{eq:preprojmod}
\sum_{ j: (i,j)\in H } \epsilon(i,j) M_{(j,i)} M_{(i,j)} = 0 
\end{equation}
for all \(i\in I\).

In type \(A\) we depict such modules by 
stacking \(\dim M_i\) digits \(i\) over each \(i\in I\) and 
drawing 
arrows connecting the digits in adjacent columns. 
The digits denote basis vectors, and the arrows may be adorned with small matrices indicating the actions \(M_h\) of \(h\in H\). 
If an arrow is unadorned, then the generic action is understood. In particular, the arrows are always injections.
\begin{example}
  \label{eg:preprojmod}
  Consider the module \(M = 1\to 2\) for the preprojective algebra of the Dynkin quiver of type \(A_2\). Once we fix \(e_i\in M_i\setminus\{0\}\) for \(i = 1,2\), i.e.\ a basis of \(M\), the action of \(\cA\) on \(M \) is well-defined by \(M = 1\to 2\):
  \[
    (1,2)e_1 = e_2 \qquad (1,2)e_2 = 0 \qquad (2,1) M = 0\,. 
  \]
\end{example}
Identifying the vertex set \(I\) with the set of simple roots \(\alpha_i\) in $\Delta_+$, 
we denote by \(\dimvec M\)
the graded dimension \(\sum_I \dim M_i \alpha_i \) of a module. 
In particular, the module that has dimension vector \(\alpha_i\) is called the simple module supported at \(i\) and denoted \(S_i\). 
The map \(\dimvec\) thus gives an isomorphism of the Grothendieck group of the category of \(\cA\)-modules
and the root lattice \(Q)\). 

Given \(\nu = \sum_I\nu_i \alpha_i \in Q_+\), let 
\(E_\nu = \bigoplus_{(i,j)\in H} \Hom(\CC^{\nu_i},\CC^{\nu_j})\). 
The adjoint action of \(\prod_I GL_{\nu_i}\) on \(E_\nu\) has moment map \(\psi = (\psi_i)_I : E_\nu \to \prod_I\gl_{\nu_i}\) defined by 
\begin{equation}
  \label{eq:luszmm}
  \psi_i(x) = \sum_{j:(i,j)\in H} \epsilon(i,j) x_{(j,i)} x_{(i,j)}\,.
\end{equation}
The closed subvariety of points \(x\in E_\nu \) such that 
\(\psi_i(x)= 0\) for all \(i\in I\) is called Lusztig's \new{nilpotent variety} of type \(\nu\) and denoted \(\Lambda(\nu)\). Its irreducible components are denoted \(\irr\Lambda(\nu)\).

We view points \(x\in\Lambda(\nu)\) as modules by sending \(x_{(i,j)} \in \Hom (\CC^{\nu_i},\CC^{\nu_j})\) to \(M_{(i,j)}\) such that in some basis of \(M\) the matrix of \(M_{(i,j)}\) is equal to the matrix of \(x_{(i,j)}\) in the standard basis of \(\oplus_I\CC^{\nu_i}\).
%

Fix \(\nu = \sum \nu_i\alpha_i \in Q_+\) having height \(\sum_i\nu_i = p \). 
The set of all sequences \((i_1,\dots,i_p)\) such that \(\sum_{k=1}^p\alpha_{i_k} = \nu\) is denoted \(\Seq(\nu)\). 
%
Given \(M\in\Lambda(\nu)\) and \(\vi\in\Seq(\nu)\), we can consider the variety \(F_\vi(M)\) of composition series of \(M\) of type \(\vi\),
\begin{equation}
  \label{eq:FiM}
  F_\vi(M) = 
  \{M = M^1 \supseteq M^2\supseteq \cdots \supseteq M^{p+1} = 0 \big| M^k/M^{k+1} \cong S_{i_k}\}\,.
\end{equation}
The varieties defined by Equation~\ref{eq:FiM} complement the paths of (co)weights in the 1-skeleton of an MV polytope witnessed in section~\ref{ss:mvps} and will be used to define the dual semicanonical basis.

Let \(\chi\) denote the topological Euler characteristic. Given \(Y\in\irr\Lambda(\nu)\), we say that \(M\) is a \new{general point} in \(Y\) if it is generic for the constructible function \(M\mapsto \chi (F_\vi (M)) \) on \(Y\), i.e.\ if it belongs to the dense subset of points of \(Y\) where \(\chi (F_\vi(M))\) is constant. 
%
%

Consider the function \(\xi_M\in\CC[U]\) that is defined by the system 
\begin{equation}
  \label{eq:cY}
    \langle e_\vi, \xi_M \rangle = \chi (F_\vi (M)) \qquad \vi \in \Seq(\nu)
\end{equation}
where \(e_\vi\) is shorthand for the product \(e_{i_1}e_{i_2}\cdots e_{i_p}\in\cU(\fu)\). 
The following result is an immediate consequence of the definition and will come in handy in Chapter~\ref{ch:calculs}. 
\begin{lemma}
  \label{lem:glshom}
  Given \(M,M'\in\Lambda\), we have \(\xi_M\xi_{M'} = \xi_{M \oplus M'}\).
\end{lemma}
Since
\begin{equation}
  \label{eq:xiM}
  \Lambda(\nu) \to \CC[U]_{-\nu} \qquad M\mapsto \xi_M 
\end{equation}
is constructible, 
for any component \(Y\in\irr\Lambda(\nu)\), we can define \(c_Y \in\CC[U]_{-\nu}\) by setting \(c_Y = \xi_M\) for any any choice of general point \(M\) in \(Y\). 
Note that \(c_Y\in\CC[U]_{-\nu}\) since \(\deg e_\vi = \nu\) for all \(\vi \in \Seq(\nu)\). 

We write \(\Lambda\) for the disjoint union \(\bigcup_{Q_+}\Lambda(\nu)\) and \(\irr\Lambda\) for the set of all irreducible components \(\bigcup_{Q_+}\irr\Lambda(\nu)\). 
\(\irr\Lambda\) indexes the dual semicanonical basis \(B(\Lambda) = \{c_{Y}\}_{\irr\Lambda}\). 

Saito, in \cite{saito2002crystal}, defined a crystal structure on \(\irr\Lambda\) and showed that it was the \(B(\infty)\) crystal structure. The following result of Lusztig \cite{lusztig2000semicanonical} is an upgrade of Saito's result. See also \cite[Theorem 11.2]{baumann2019mirkovic}.
\begin{theorem}
  \label{thm:cYperf}
  The set 
  \(B(\Lambda)\) is a perfect basis of \(\CC[U]\). 
\end{theorem} 
%

By Theorem~\ref{thm:perfbunique}, Theorem~\ref{thm:cYperf} implies that \(B(\Lambda)\) is in unique bijection with \(B(\infty)\).
%
\begin{example}
  \label{eg:Lambda12}
  Let's examine some of the ingredients of this section in the case \(G = PGL_3\). 
  The preprojective algebra \(\cA\) attached to the \(A_2\) quiver
  %
  $$
  \begin{tikzcd}
    \bullet\ar[r,bend right,"h" description] & \bullet\ar[l,bend right, "\bar h" description]
  \end{tikzcd}
  $$
  is the algebra generated by the arrows \(h\) and \(\bar h\) subject to the relations \(h\bar h = 0\) and \(\bar h h = 0\). 
The resulting nilpotent variety has connected components 
\[
  \Lambda(\nu_1,\nu_2) = \{(x_h,x_{\bar h}) \in \Hom(\CC^{\nu_1},\CC^{\nu_2})\oplus \Hom(\CC^{\nu_2},\CC^{\nu_1}) \big| x_hx_{\bar h} = 0 \text{ and } x_{\bar h} x_h = 0\}\,,
\] 
indexed by roots \(\nu = \nu_1 \alpha_1 + \nu_2 \alpha_2 \in Q_+\). 
For instance, the connected component of \(\Lambda\) made up of those modules which have \(\dimvec M = \alpha_1 + \alpha_2\) decomposes into two irreducible components as follows. 
  \[
  \begin{aligned}
      \{(x_h,x_{\bar h} )\in \CC^2 \big| x_h x_{\bar h} = 0\} = \{x_h = 0\}\bigsqcup\{x_{\bar h} = 0\} 
      = \{1 \xleftarrow{\lambda} 2\}\bigsqcup_{\{1\oplus 2\}}\{1\xrightarrow{\lambda} 2\}
      = \CC^\times\bigsqcup_{\{0\}}\CC^\times\,. 
  \end{aligned} 
  \]  
\end{example}
\begin{definition}
  \label{def:hnfiltr}
The \new{HN (Harder--Narasimhan) filtration of a module} \(M\in\Lambda \) is defined to be the unique decreasing filtration 
\begin{equation}
  \label{eq:hnfiltr}
  M = M^1\supset M^2\supset\cdots\supset M^\ell\supset M^{\ell+1} = 0
\end{equation}
such that 
\(M^{k}/M^{k+1}\cong (B_{k}^{\uvi})^{\oplus n_k} \) for all \(1\le k\le \ell \), and 
the \(\uvi\)-\new{Lusztig datum of \(M\)}, \(n^{\uvi}(M)\), is the tuple of multiplicities \(n^{\uvi}(M)_{k} = n_k\).
Here \(B_k^{\uvi}\)
is a certain indecomposable module of dimension
\(\beta_k^{\uvi}\).
Just like \(\beta_k^{\uvi}\) is obtained by applying a sequence of reflections to \(\alpha_k\), \(B_k^{\uvi}\) is built up from the simple module \(S_k\) by applying a sequence of reflection functors.
For further details, see \cite[Theorem 5.10(i), Example 5.14]{baumann2014affine}.
For the slope-theoretic definition of an HN filtration, see \cite[\S 1.4]{baumann2014affine}. 
\end{definition}
%

It turns out that \(n^{\uvi}\) is also constructible, and we can 
define \(n^{\uvi }(Y) \) by \(n^{\uvi}(M)\) for any \(M\) that is a general point in \(Y\).
%
\begin{theorem}
~\cite[Remark 5.25(ii)]{baumann2014affine}
\label{thm:modld}
Taking \(\uvi\)-Lusztig datum is a bijection \(\irr\Lambda\to \NN^\ell\). It agrees with the composite map \(\irr(\Lambda) \to B(\infty) \to \NN^\ell\).
\end{theorem}
\begin{example}
  \label{eg:Lambda12lds}
  Continuing with the previous example, fix the reduced word \(\uvi = (1,2,1)\). The indecomposable modules for this choice of reduced word are \(B_1 = 1\), \(B_2 = 1 \leftarrow 2\), and \(B_3 = 2\).
We would like to use the definition to find the Lusztig data of \(\irr\Lambda(\alpha_1+\alpha_2)\).
\(M = 1 \leftarrow 2\) has HN filtration \(M \supset M \supset 0 \supset 0 \). The resulting ``HN table'' of \(B_k\)-subquotients is
  \[
    \begin{BMAT}{c}{c|c|c}
      M^1/M^2 \\ 
      M^2/M^3 \\ 
      M^3 
    \end{BMAT} = 
    \begin{BMAT}{c}{c|c|c}
      0 \\
      B_2 \\
      0
    \end{BMAT} \,.
  \]
  Therefore the irreducible component of \(\Lambda(\alpha_1 + \alpha_2)\) for which \(M\) is generic has Lusztig datum \((0,1,0)\). 
  Similarly, the irreducible component for which \(M = 1\to 2\) is generic has Lusztig datum \((1,0,1)\), since \(M\) has HN filtration \(M \supset B_2 \supset B_2 \supset 0 \), and, hence, HN table
  \[
    \begin{BMAT}{c}{c|c|c}
      M^1/M^2 \\ 
      M^2/M^3 \\ 
      M^3 
    \end{BMAT} = 
    \begin{BMAT}{c}{c|c|c}
      B_1 \\
      0 \\
      B_2
    \end{BMAT} \,. 
  \]
\end{example}
Henceforth, we will represent modules and their HN filtrations by HN tables. 
\begin{definition}
  \label{def:hnpol}
Let \(Y\in\irr\Lambda\) and take \(M\in Y\) generic. The set 
\begin{equation}
  \label{eq:hnpol}
\nu - \Conv(\{\dimvec N \big| N \subseteq M \text{ as }\cA\text{-submodule}\})
\end{equation}
is called the HN polytope of \(Y\) and denoted \(\Pol(Y)\). \(\Pol\) is yet another constructible map.%
\footnote{In \cite{baumann2014affine} the HN polytope of a module \(M\) is defined as \(\Conv(\{\dimvec N \big| N \subseteq M \text{ as }\cA\text{-submodule}\})\). 
The discrepancy is due to the way we filter our modules by quotients instead of submodules. Indeed, we could equally have defined our polytope as \(\Conv(\{\dimvec N\big| N\text{ is a quotient of } M\})\).}
\end{definition}
\begin{example}
  \label{eg:hnpol11}
  Take \(\nu = \alpha_1 + \alpha_2\). We can check that the HN polytopes of \(\irr\Lambda(\nu)\) match up with the MV polytopes computed in Example~\ref{eg:pol11} according to their Lusztig data. 
  \[
  \begin{aligned}
    & \nu - \Conv(\{\dimvec N\big| N\subset 1\to 2\}) & & \nu - \Conv(\{\dimvec N\big| N\subset 1\leftarrow 2\}) \\
    &= \nu - \Conv(\{\dimvec 0, \dimvec 2, \dimvec (1\to 2)\}) & &= \nu - \Conv(\{\dimvec 0, \dimvec 1, \dimvec (1\leftarrow 2)\}) \\
    &= \nu - \Conv(\{0,\alpha_2,\alpha_1 + \alpha_2 \} ) & &= \nu - \Conv(\{0 , \alpha_1 , \alpha_1 + \alpha_2\}) \\
    &= \Conv(\{\alpha_1 + \alpha_2, \alpha_1 , 0 \}) & &= \Conv(\{\alpha_1 + \alpha_2, \alpha_2 , 0 \})
  \end{aligned}  
  \]
\end{example}

\subsection{Semistandard Young tableaux}
\label{ss:ssytld}
The final model we consider is only available in type \(A\). We therefore fix throughout this section \(G = PGL_m\). 
%
Consider the set of Young diagrams having at most \(m-1\) rows. 
\begin{equation}
  \label{eq:yng}
  \cY(m) = 
  \{\lambda = (\lambda_1,\dots,\lambda_m)\in\ZZ^m\big|\lambda_1\ge\cdots\ge\lambda_m = 0\}\,.
\end{equation}
Later we will identify this set with the set of dominant coweights of \(G\). 
Given \(\lambda\in\cY(m)\), we can consider the set of semistandard Young tableaux of shape \(\lambda\), which we denote by \(\cT(\lambda)\). 
We say that \(\tau\in\cT(\lambda)\) has ``coweight'' \(\mu= (\mu_1,\dots,\mu_{m})\in\NN^m\) if it has \(\mu_1\) 1s, \(\mu_2\) 2s,\dots, and \(\mu_{m}\) \(m\)s and we write \(\cT(\lambda)_\mu\) for the set of semistandard Young tableaux of shape \(\lambda\) and coweight \(\mu\).

Given such a tableau \(\tau\in\cT(\lambda)_\mu\), we write \(\tau^{(i)}\) for the tableau obtained from \(\tau\) by deleting all boxes numbered \(j\) for any \(j > i \), and we write \(\lambda^{(i)}_\tau\) for its shape. The list of tuples \((\lambda_\tau^{(i)})\) (when it is arranged in a triangular configuration that looks locally like 
\[
  \begin{smallmatrix}
& \lambda^{(i)}_k \\
\lambda^{(i+1)}_k & & \lambda^{(i+1)}_{k+1}
\end{smallmatrix}
\]
satisfying \(\lambda^{(i+1)}_{k+1}\le \lambda^{(i)}_k\le \lambda^{(i+1)}_k\) for all \(1\le k\le \mu_1 + \cdots + \mu_i\), for all \(1\le i\le m-1\)) is referred to by Berenstein and Zelevinsky in \cite[\S 4]{BERENSTEIN1988453} as the GT (Gelfand--Tsetlin) pattern of \(\tau\). It is needed to define the Lusztig datum of \(\tau\). 

Throughout this thesis we will only work with the reduced word 
\(\uvi = (1,2,\dots,m-1,\dots,1,2,1)\). For this reason we often omit it from the notation, writing 
\(n(\phantom{0})\) in place of \(n^{\uvi}(\phantom{0})\) and \(\beta_k\) in place of \(\beta^{\uvi}_k\).
This choice yields \(n(\tau) = (n_1,\dots,n_{\ell})\) where
\[
  n_k
  = \lambda^{(j)}_i - \lambda^{(j-1)}_i \text{ whenever }k\text{ is such that } \beta_k = \alpha_i + \cdots + \alpha_j 
  \,.
\] 
Other choices of reduced word will result in a more complicated (possibly piecewise linear) formula for \(n(\tau)\). 
We will sometimes write \(n_{(i,j)}\) for \(n_k\) and \(\alpha_{i,j}\) for \(\beta_k\). 
%
\begin{example}
\label{eg:ldgt}
The tableau 
\[
  \young(112,233)
\] 
has GT pattern 
\[
\begin{tikzpicture}
  [scale=0.5]
  \node (22) at (2,2){$\lambda_1^{(1)}$}; \node (11) at (1,1){$\lambda_1^{(2)}$}; \node (00) at (0,0){$\lambda_1^{(3)}$};
  \node (31) at (3,1){$\lambda_2^{(2)}$}; \node (20) at (2,0){$\lambda_2^{(3)}$};
  \node (40) at (4,0){$\lambda_3^{(3)}$};
  \node (equals) at (5.3,1.3) {$=$}; 
\node (22) at (8,2){$2$}; \node (11) at (7,1){$3$}; \node (00) at (6,0){$3$};
\node (31) at (9,1){$1$}; \node (20) at (8,0){$3$};
\node (40) at (10,0){$0$};
\end{tikzpicture}
\]
and Lusztig datum \(n(\tau) = (1,0,2)\). 
Note, \(n(\tau)\) can be easily read off of the GT pattern of \(\tau\): it is the string of Southwesterly differences starting with (the top of) the left side of the triangle, and moving inward toward the bottom right corner, as pictured below. 
\[
  \begin{tikzpicture}[scale=0.5]
    \draw[->,gray] (2.4,1.6) -- node[near start] {1} node[near end] {0} (0.4,-0.4);
    \draw[->,gray] (3.4,0.6) -- node[midway] {2} (2.4,-0.4);
    \node (22) at (2,2){$2$}; \node (11) at (1,1){$3$}; \node (00) at (0,0){$3$};
    \node (31) at (3,1){$1$}; \node (20) at (2,0){$3$};
    \node (40) at (4,0){$0$};
  \end{tikzpicture}  
\]
\end{example}
\begin{theorem}
\label{thm:tabmseg}
(cf.\ \cite[Theorem 3.11]{claxton2015young}) 
There exists an embedding \(\cT(\lambda)\to B(\infty)\) such that the composite map \(\cT(\lambda)\to B(\infty)\to \NN^\ell\) is given by \(\tau\mapsto n(\tau)\).
%
Conversely, given \(n\in\NN^\ell\) of weight \(\nu\in Q\), there exists \(\mu\in P\) and \(\tau\in\cT(\nu + \mu)\) such that \(n(\tau) = n\) and the shape of \(\tau\) is the smallest possible. In particular, \(\nu + \mu\) is effective dominant. 
\end{theorem}
Let us describe the embedding \(\cT(\lambda) \to B(\infty)\) in terms of the so-called multisegment 
model for \(B(\infty)\). 

A multisegment is a collection of segments \([i,j]\) for various \(1\le i\le j\le m-1\) allowing multiplicity. The basic segments \([i,j]\) are in bijection with the simple roots \(\alpha_{i,j}\). 
We depict them as columns of height \(j-i+1\) having base \(i\). The \(\uvi\)-Lusztig datum of a multisegment is the tuple of multiplicities with which the basic segments \([i,j]\) occur in it. The tuple is ordered so that the multiplicity of a segment \([i,j] \) precedes the multiplicity of a segment \([i',j']\) whenever \(\alpha_{i,j}\le_{\uvi}\alpha_{i',j'}\) in the convex order on \(\Delta_+\) induced by \(\uvi\).
%
The tableau \(\tau\in\cT(\lambda)\) determines the multisegment whose \(\uvi\)-Lusztig datum is equal to that of \(\tau\). 
%

The true value of the multisegment model for us is it is convenient for describing the partial inverse \(B(\infty)\dashrightarrow \cT(\lambda)\)
gauranteed by the second half of Theorem~\ref{thm:tabmseg}. Consider the following example. 
\begin{example}
\label{eg:tabmseg}
In type \(A_5\) we fix the reduced expression \(\uvi = (1,2,3,4,5,1,2,3,4,1,2,3,1,2,1)\) and take \(n = (0,1,0,0,0,0,0,1,0,1,0,0,0,1,0)\). 
This determines the multisegment \(m\) depicted below.
\[
  %
  \Yvcentermath0
\young(3)\,\young(5,4)\,\young(2,1)\,\young(4,3,2)  
\]
Note that the columns (basic segments) are strategically ordered such that they are increasing in height and decreasing in base, with height taking priority.
From this picture we can compute \(\varepsilon_2^\ast(m) = \varepsilon_4^\ast(m) = 1\) and all other \(\varepsilon_i^\ast( m) = 0\). 

These quantities\footnote{Related to \(\varepsilon_i\) appearing in the definition of a crystal; see \cite{claxton2015young} for details.} are calculated by placing ``('' under each segment with base $i$, ``)'' under each segment with base $i+1$, and counting the unpaired ``(''s. Here, a ``('' and a ``)'' are called paired if the ``('' is directly to the left of the ``)'', or if the only parentheses between them are paired. 

For example, to see that \(\varepsilon_2^\ast(m) = 1\), we note that there is only one column with base \(2\). 
We count \(\varepsilon_2^\ast(m) = 1\) unpaired open parenthesis:
\[
  \Yvcentermath0
  \underset{)}{\young(3)}\,\underset{\phantom{0}}{\young(5,4)}\,\underset{\phantom{0}}{\young(2,1)}\,\underset{(}{\young(4,3,2)} 
\]

The tuple \((\varepsilon_i^\ast(m))\) contains the information of the smallest possible \(\lambda\) such that some \(\tau\in\cT(\lambda)\) has Lusztig datum equal to that of \(m\). Namely, \(\lambda = \sum\varepsilon_i^\ast(m)\omega_i \).
In this case, the smallest possible shape of a tableau having prescribed Lusztig datum \(n\) is therefore \(\omega_2 + \omega_4 = (2,2,1,1)\) and the unique such tableau is 
\[
  \young(13,25,4,6)
\]
The reader can check that the GT pattern of \(\tau \) does in fact yield Lusztig datum \(n\).
\end{example}

\section{More type \texorpdfstring{\(A\)}{A} show and tell}
\label{s:gue1}
As 
our creative take on a leitfaden, let us summarize the entirety of the data we will have covered by the end of this thesis in the context of 
our running \(G = PGL_3\) example of weight \(\nu = \alpha_1 + \alpha_2\), including Example~\ref{eg:preprojmod}.
Take the reduced word \(\uvi = (1,2,1)\) and zoom in on the Lusztig datum \(n = (1,0,1)\). 
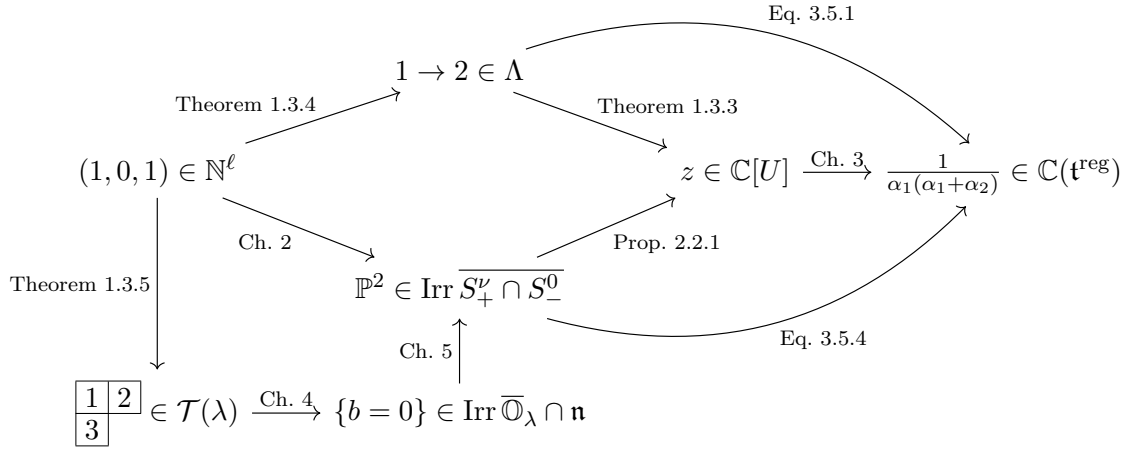
\begin{figure}[ht!]
  \centering
  \begin{tikzcd}
    & 1\to 2 \in \Lambda \ar[dr,"\text{Theorem~\ref{thm:cYperf}}"] \ar[drr,bend left,"\text{Eq.~\ref{eq:barDc}}"]\\
    (1,0,1) \ar[dd,"\text{Theorem~\ref{thm:tabmseg}}"']\in \NN^\ell \ar[ur,"\text{Theorem~\ref{thm:modld}}"] \ar[dr,"\text{Ch.~\ref{ch:gr}}"'] &  & z\in\CC[U]\ar[r,"\text{Ch.~\ref{ch:comparison}}"] & \frac 1 {\alpha_1(\alpha_1 + \alpha_2)} \in \CC(\treg)\\
    & \PP^2 \in \irr\overline{S^\nu_+\cap S^0_-} \ar[ur,"\text{Prop.~\ref{prop:uniquebZ}}"'] \ar[urr,bend right,"\text{Eq.~\ref{eq:barDb}}"'] \\
    \young(12,3) \in \cT(\lambda) \ar[r,"\text{Ch.~\ref{ch:govs}}"] & \{b = 0\} \in \irr \overline{\OO}_\lambda\cap \n \ar[u,"\text{Ch.~\ref{ch:mvs}}"]
  \end{tikzcd}
  \caption{Leitfaden}
  \label{fig:leitfaden}
\end{figure}
%

Here \(\OO_\lambda \) denotes the nilpotent orbit of \(3\times 3\) matrices having Jordan type \(\lambda = (2,1)\), \(\n\) is the set of upper triangular matrices with coordinates
\(\left[\begin{smallmatrix}
  0 & a & c \\
  0 & 0 & b \\
  0 & 0 & 0 
\end{smallmatrix}\right]\), and \(U \) is the set of unitriangular matrices with coordinates 
\(\left[\begin{smallmatrix}
  1 & x & z \\
  0 & 1 & y \\
  0 & 0 & 1
\end{smallmatrix}\right]\). 

By the end of Chapter~\ref{ch:mvs}, the reader should be able to fill in the complementary data associated to the other \(\uvi\)-Lusztig datum of weight \(\nu\), \(n = (0,1,0)\). 

The MV polytope model is absent from our spaghetti map because, while necessary to make sense of the measures defined in Chapter~\ref{ch:comparison}, the information it contains is equivalent to the information of the Lusztig data.

%
%
%
We conclude this chapter with a playful attempt to make up for this gap by presenting one 
of the five polytopes of weight \(\nu = \alpha_1 + 2\alpha_2 + \alpha_3\) in type \(A_3\). 
Let \(M\in\Lambda(\nu)\) be a module of the following form.  
\[
  \begin{tikzpicture}
    \node (00) at (0,0){1}; 
    \node (11) at (1,1){22};
    \node (20) at (2,0){3};
    \draw[->] (11) -- (00);
    \draw[->] (11) -- (20); 
  \end{tikzpicture}
\]
Then, with respect to \(\uvi = (1,2,3,1,2,1)\), \(M\) has HN table 
\[
  \begin{BMAT}{c}{c|c|c|c|c|c}
    0 \\
    1\leftarrow 2\\
    0 \\
    2 \\
    0 \\
    3
  \end{BMAT}
\]
and, hence, \(\uvi\)-Lusztig datum \((0,1,0,1,0,1)\).
%
%
The HN polytope of \(M\)---also the MV polytope with this Lusztig datum---is depicted below.
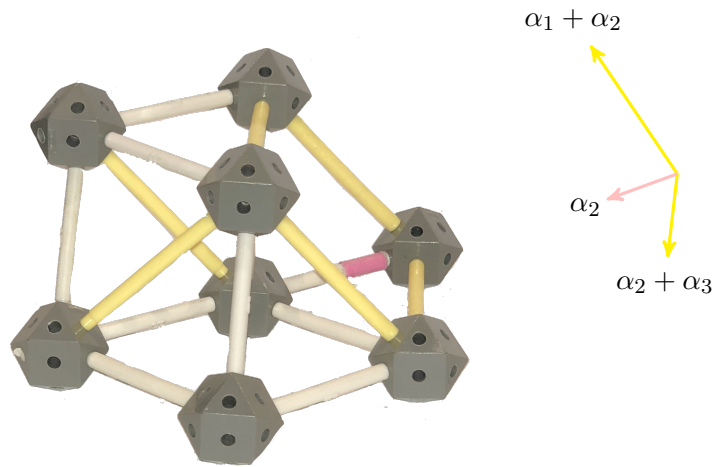
\begin{figure}[ht!]
  \centering
\raisebox{-.4\height}{\includegraphics[width=7cm]{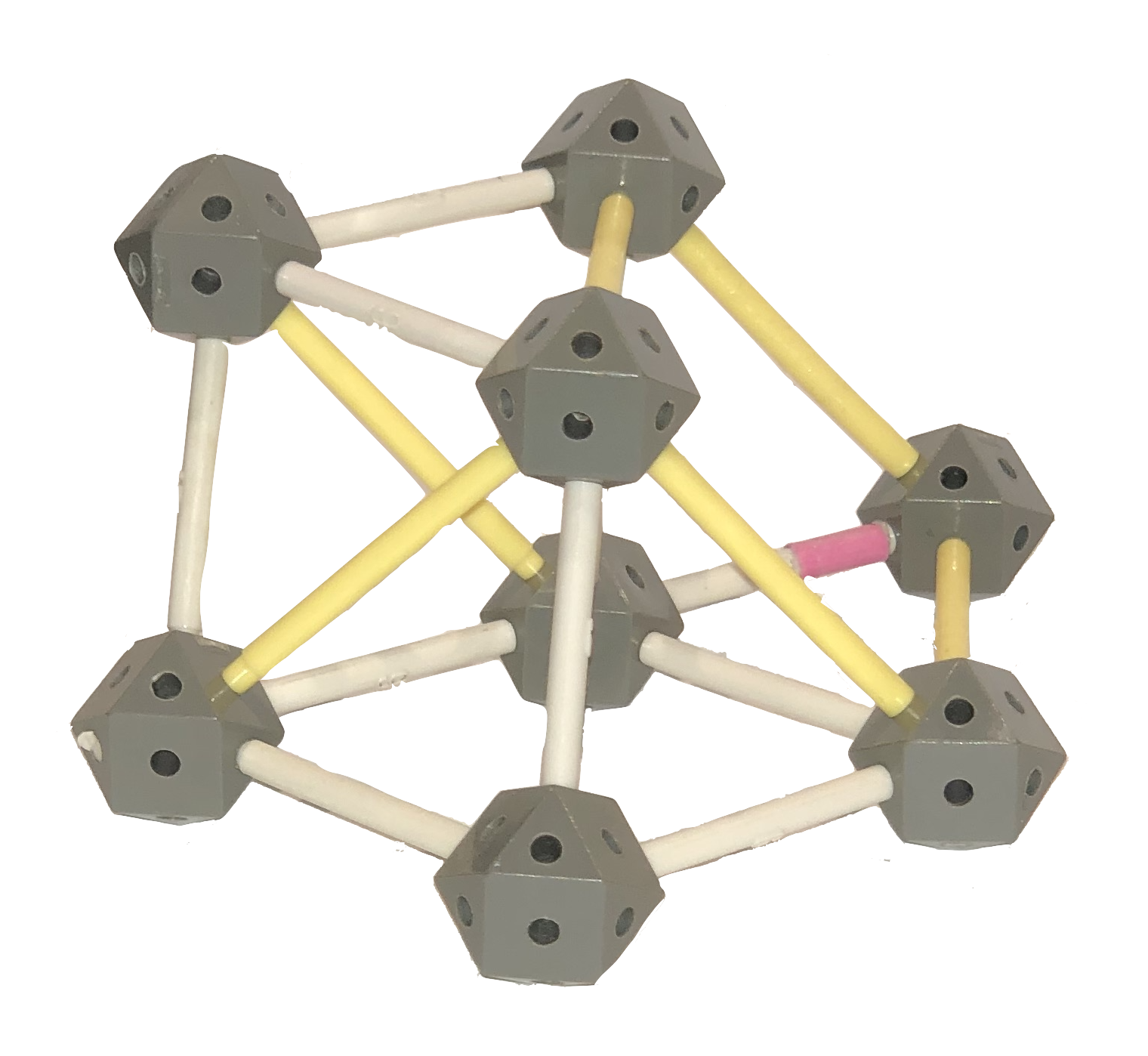}}
\begin{tikzpicture}[rotate=20,transform shape]
  \draw[very thick,->,>=stealth',yellow] (5,0) -- node[rotate=-20,black,pos=1.2] {$\alpha_1+\alpha_2$} (4.5,2);
  \draw[very thick,->,>=stealth',yellow] (5,0) -- node[rotate=-20,black,pos=1.3] {$\alpha_2 + \alpha_3$} (4.5,-1);
  \draw[very thick,->,>=stealth',pink] (5,0) -- node[rotate=-20,black,pos=1.3] {$\alpha_2$} (4,0); 
\end{tikzpicture}
\caption{The polytope of \(M\), \(\Conv(\{0,\alpha_2 + \alpha_3, 2\alpha_2 , 2\alpha_2 + \alpha_3, \alpha_1 + \alpha_2, \alpha_1 + 2\alpha_2, \alpha_1 + 2\alpha_2 + \alpha_3\})\)}
\label{fig:3212}
\end{figure}

\chapter{Background on the affine Grassmannian}
\label{ch:gr}
%
%
In this chapter we introduce the second perfect basis with which this thesis is concerned. Its definition requires a bit more setup. 

\section{MV cycles as a basis of \texorpdfstring{\(L(\lambda)\)}{highest weight representations}}
\label{s:mvb}
Let \(G^\vee\) denote the canonical smooth complex semisimple simply connected group whose root datum is dual to that of \(G\). (See \cite{demazure1965schemas} for the construction.) 
\(G^\vee\) comes with a maximal torus whose coweight lattice is \(P\), the weight lattice of \(G\). 
Our choice of positive and negative roots provide a pair of dual opposite Borel subgroups \(B^\vee_{\pm}\) with unipotent radicals \(U^\vee_{\pm}\).

Fix the indeterminate \(t\). Let \(\cO \) denote the power series ring \(\CC\xt\), and \(\cK\) its fraction field, the ring of Laurent series \(\CC\xT\). 
The affine Grassmannian of \(G^\vee\) is defined to be the quotient \(G^\vee(\cK)/G^\vee(\cO) \), and is denoted by \(\Gr \). 
We identify \(G^\vee \) with its functor of points, so that \(G^\vee(R) = \Hom(\spec R, G^\vee)\) for any ring \(R\).

%
The torus \(T^\vee\) acts on \(\Gr\) as a subgroup of \(G^\vee(\cK)\) which acts on \(\Gr \) by left multiplication. 
Given \(\lambda\in P\equiv\Hom(\CC^\times,T^\vee)\),
we denote by \(t^\lambda\) its image in \(T^\vee (\cK)\) and by \(L_\lambda\) 
its image \(t^\lambda G^\vee (\cO)\) in \(\Gr\).
Consider the orbits \(\Gr^\lambda = G^\vee(\cO)L_\lambda\) and \(S_{\pm}^\lambda = U_{\pm}^\vee(\cK) L_\lambda\) for any \(\lambda\in P\). The former are finite-dimensional quasi-projective varieties, while the latter are infinite-dimensional ``semi-infinite cells''. 
Like Schubert cells in a flag variety, the semi-infinite cells are attracting for a \(\CC^\times\) action on \(\Gr\).

The geometric correspondence of Mirkovi\'c and Vilonen describes an equivalence of categories between representations of \(G\) and certain sheaves on \(\Gr\), relating the highest weight irreducible representations \(L(\lambda)\) of \(G\) to the \(G^\vee(\cO)\) orbits \(\Gr^\lambda\). More precisely, the correspondence proceeds as follows.
\begin{theorem}
  \label{thm:mvs}
\cite[Proposition 3.10]{mirkovic2007geometric}
For each \(\lambda\in P_+\) there is a decomposition 
\begin{equation}
  \label{eq:ihgr}
  IH_\bullet(\overline{\Gr^\lambda}) = \bigoplus_{\mu\le\lambda}H_{\topp}(\overline{\Gr^\lambda\cap S^\mu_-})
\end{equation}
with \(IH\) denoting intersection homology. Moreover, 
\begin{equation}
  \label{eq:geosat}
  L(\lambda)_\mu \cong H_{\topp}(\overline{\Gr^\lambda\cap S^\mu_-}) 
\end{equation}
for all \(\mu\in P\) such that \(\mu\le\lambda\).
\end{theorem}
Let \(\rho^\vee\) denote the half-sum of positive coroots of \(G\). By \cite[Theorem 3.2]{mirkovic2007geometric}, the intersection \(\Gr^\lambda\cap S^\mu_-\) has pure dimension 
\(\langle\rho^\vee,\lambda-\mu\rangle \). 
Therefore, the set of irreducible components of
\(\overline{\Gr^\lambda\cap S^\mu_-}\), which we denote by 
\(\cZ(\lambda)_\mu\), gives a basis for the vector space \(H_{\topp}(\overline{\Gr^\lambda\cap S^\mu_-})\). 
These are called the MV (Mirkovi\'c--Vilonen) cycles of type \(\lambda\) and weight \(\mu\).
%
We denote by \(\cZ(\lambda)\) the disjoint union of all MV cycles of type \(\lambda\), \(\cZ(\lambda) = \bigcup_{\mu\le\lambda} \cZ(\lambda)_\mu \). 
Given \(Z\in \cZ(\lambda)\) we write \([Z]\) for the corresponding basis vector in \(L(\lambda)\) with the isomorphisms in Equation~\ref{eq:geosat} normalized so that 
\(v_\lambda = [\{L_\lambda\}]\). 
Note that \(\{L_\lambda\} = \Gr^\lambda\cap S^\lambda_-\).
\section{Stable MV cycles as a basis of \texorpdfstring{\(\CC[U]\)}{the ring of functions on the unipotent subgroup}}
\label{s:stabmvbasis}
Consider the basis \(B(\lambda) = \{[Z]\}_{Z\in \cZ(\lambda)}\) of \(L(\lambda)\). 
We know from Chapter~\ref{ch:bases} that the family of bases \(\bigcup_{\lambda\in P_+}B(\lambda)\)
can be studied simultaneously,
and in this section we describe how this is done.

Given \(\nu\in Q_+\), we define the \new{stable MV cycles} of weight \(\nu\) to be the irreducible components of \(\overline{S^\nu_+\cap S^{0}_-}\). We denote these irreducible components by \(\cZ(\infty)_{\nu}\) and we denote their disjoint union \(\bigcup_{Q_+} \cZ(\infty)_{\nu}\) by \(\cZ(\infty)\). 

The weight lattice \(P\) 
acts on \(\Gr\) by left multiplication, with \(\mu \cdot L = t^\mu L \)
for any \(L\in \Gr\) and for any \(\mu\in P\). 
Since \(t^\mu\in T^\vee \subset N_G (U^\vee_w)\), this action translates semi-infinite cells, with \(\nu\cdot S^\mu_w = S^{\mu + \nu}_w\).
In particular, we get an isomorphism of \(\overline{S^\nu_+\cap S^{0}_-}\) and \(\overline{S^{\nu + \mu}_+\cap S^{\mu}_-}\) and, hence, a bijection of irreducible components. 

By \cite[Proposition 3]{anderson2003polytope}, an MV cycle of type \(\lambda\) and weight \(\mu\) can be described as an irreducible component of \(\overline{S^\lambda_+\cap S^\mu_-}\) that is contained in \(\overline{\Gr^\lambda}\). 
Thus, the set of stable MV cycles of weight \(\nu = \lambda-\mu\) whose \(\mu\)-translates are contained in \(\overline{\Gr^\lambda}\) is in bijection with MV cycles of type \(\lambda\):
\begin{equation}
  \label{eq:inftylambda}
  \{Z\in\cZ(\infty)\big|\mu\cdot Z \subset \overline{\Gr^\lambda}\}\rightarrow\cZ(\lambda)\,.
\end{equation}
This observation is key to transporting the family of bases \(\bigcup_{\lambda\in P_+} B(\lambda)\) to a basis of \(\CC[U]\). 
%
\begin{proposition}
  \label{prop:uniquebZ}
\cite[Proposition 6.1]{baumann2019mirkovic}
For each \(\nu\in Q_+\) and for each \(Z\in\cZ(\infty)_\nu\), there exists a unique element \(b_Z\in\CC[U]_{-\nu}\) such that, for any \(\mu\) such that \(\nu + \mu \in P_+\),
\begin{equation}
  \label{eq:whatisbZ}
  \mu\cdot  Z\subset\overline{\Gr^{\nu + \mu}}\,
  \Longrightarrow\, 
  b_Z = \Psi_{\nu + \mu}([\mu\cdot Z]) \in \CC[U]_{-\nu} \,.
\end{equation}
\end{proposition}
By \cite[Proposition 6.2]{baumann2019mirkovic}, the elements \(b_Z\) guaranteed by Proposition~\ref{prop:uniquebZ} form a perfect basis of \(\CC[U]\). This basis is called the \new{MV basis} of \(\CC[U]\) for it is obtained by gluing the MV bases of the representations \(L(\lambda)\). We denote it by \(B(\Gr)\). 
%

\section{MV polytopes from MV cycles} 
\label{s:xva}
Let \(A\) be a torus and let \(V\) be a possibly infinite-dimensional representation of \(A\). 
Consider a finite-dimensional \(A\)-invariant closed subvariety \(X\) of \(\PP(V)\). 
We denote the set of \(A\)-fixed points of \(X\) by \(X^A\) and assume that every \(x\in X^A\) is isolated.

The moment polytope of \(X\) is defined with the help of the dual determinant line bundle \(\cO(n)\to\PP(V)\). 
%
\(\cO(n)\) has a natural \(A\)-equivariant structure, with the weight of \(A\) on the fibre of \(\cO(n)\) over $L = [v]\in\PP(V)$ equal to \(-n\) times the weight of \(A\) on the representative \(v\).
This follows by definition: \(\cO(n)\) has fibres \(\cO(n)_L = (L^{\otimes n})^* = (L^\ast)^{\otimes n}\) and the weight of \(A\) on \(L^\ast\) is equal to minus the weight of \(A\) on \(L\). 

Equivariance of \(\cO(n)\) makes the space of sections \(H^0(X,\cO(n))\) into a representation, with the action of \(a\in A\) given pointwise by $(a\cdot f)(x) = a\cdot f(a^{-1}\cdot x)$ for any $f\in H^0(X,\cO(n))$ and for any $x\in X$.

Given a fixed point \(x\in X^A\), we denote by \(\Phi_A(x)\) the weight of \(A\) on the fibre of \(\cO(1)\) at \(x\). 
Since \(V^\ast \cong H^0 (\PP(V),\cO(1))\), this is the same as minus the weight of \(A\) on \(v\in V\) for any \(v\) such that \([v] = x \).

\begin{definition}
  \label{def:mompol}
The \new{moment polytope} of the triple \((X,A,V)\) is given by
\begin{equation}
  \label{eq:mompol}
  \Conv(\{\Phi_A(x)\big| x\in X^A\}) 
\end{equation}
and denoted \(\Pol(X)\). 
\end{definition}
We would like to apply this general setup to MV cycles. 
To do so we require a projective action of \(A = T^\vee\) on \(\PP(V)\) that comes from a linear action of \(T^\vee\) on \(V\). 
The following exposition is based on \cite{mirkovic2007geometric}.

By definition, the roots of \(G^\vee \) are equal to the coroots \(\Delta^\vee\) of \(G\).
%
Let \(\CC^\times\) be the subgroup of automorphisms of \(\cK\) that rotate \(t\). 
Consider the semidirect product \(G^\vee (\cK)\rtimes\CC^\times\) with respect to this rotation action. It has maximal torus \(T^\vee \times\CC^\times\).

Let \(\widehat{G^\vee (\cK)}\) be a central extension of 
\(G^\vee (\cK)\rtimes\CC^\times\)
by \(\CC^\times\). 
It turns out that the root system of \(\widehat{G^\vee(\cK)}\) is equal to the root system of \(G^\vee (\cK)\rtimes\CC^\times\). Indeed, the central $\CC^\times$ acts by a scalar 
called the ``level'' on any given irreducible representation.
Moreover, the eigenspaces of 
\(T^\vee \times\CC^\times\) in the loop Lie algebra \(\g^\vee(\cK)\) are 
\[
    \g^\vee (\cK)_{k\delta^\vee + \alpha^\vee} = t^k \g^\vee_{\alpha^\vee} \qquad (k\in\ZZ,\,\alpha^\vee \in\Delta^\vee\cup\{0\})
\]
with \(\delta^\vee\) denoting the character of \(T^\vee\times\CC^\times\) which is trivial on \(T^\vee\) and the identity on \(\CC^\times \).
It follows that the root system of \(G^\vee(\cK)\) is given by nonzero characters of \(T^\vee\times\CC^\times\) of the form \(\alpha^\vee + k\delta^\vee \) for some \(\alpha^\vee\in\Delta^\vee\cup\{0\}\) and for some \(k\in\ZZ\). 

Consider the level 1 highest weight zero irreducible representation of $\widehat{G^\vee(\cK)}$, denoted \(L(\Lambda_0)\).\footnote{\(L(\Lambda_0)\) is sometimes also called the basic representation} 
Let \(v_{\Lambda_0}\) be a highest weight vector for \(L(\Lambda_0)\). 
Since \(U_\pm^\vee,\, T^\vee\) and the rotating \(\CC^\times\) all fix \(v_{\Lambda_0}\) the map \(g\mapsto gv_{\Lambda_0} \) induces an embedding of \(\Gr\) into \(\PP(L(\Lambda_0))\). See Proposition~\ref{prop:grembed}. 
%

We can therefore take \((X,A,V) = (Z,T^\vee,L(\Lambda_0))\) for any MV cycle \(Z\).
Remarkably, if \(Z\) is an MV cycle of 
of type \(\lambda\) and weight \(\mu\), then applying the moment polytope construction of Equation~\ref{eq:mompol} to \((Z,T^\vee,L(\Lambda_0))\) results in an MV polytope of weight
\((\lambda,\mu)\). 
\begin{theorem}(\cite[Theorem 3.1]{kamnitzer2010mirkovic})
  \label{thm:kamvgeo}
  Let \(P\) be an MV polytope with vertices \((\mu_w)_W\). Then
  \begin{equation}
    \label{mvcgeo}
    Z = \overline{\bigcap S^{\mu_w}_w}
  \end{equation}
  is an MV cycle of weight of type \(\mu_{e}\) and weight \(\mu_{w_0}\),
  and its moment polytope,
\begin{equation}
  \label{eq:mompolZ}
  \Pol(Z) = \Conv(\{\mu\big| L_\mu\in Z^{T^\vee}\})\,,
\end{equation} 
is such that \(P = \Pol(Z)\). 
\end{theorem}
It follows that we can define the Lusztig datum of an MV cycle to be the Lusztig datum of its moment polytope.
On the other hand, the proof of this theorem relies on certain constructible functions which can be used to define Lusztig data of MV cycles directly.
In Chapter~\ref{ch:mvs} we will see explicit formulae for these functions in type \(A\). 
%
  
  As a special case we see that if \(Z\in\cZ(\infty)_\nu\) then 
  \(\Pol(Z)\) is a stable MV polytope of weight \((\nu,0)\)
  having lowest vertex \(\mu_e = 0\) and highest vertex \(\mu_{w_0} = \nu\).

\chapter{Further comparison}
\label{ch:comparison}

In this chapter we set up the geometric comparison of the perfect bases \(B(\Gr)\) and \(B(\Lambda)\) of \(\CC[U]\). 
We will use the nondegenerate \(W\)-invariant bilinear form \((\phantom{0},\phantom{0})\) on \(P\) fixed in Chapter~\ref{ch:bases} to identify \(\t^\ast\) and \(\t\). 

\section{Measures from \texorpdfstring{\(\CC[U]\)}{the ring of functions on the unipotent subgroup}}
\label{s:mfromCU}
In this section we use the pairing defined by Equation~\ref{eq:perfpair} to associate to elements of \(\CC[U]\) measures on the real weight lattice \(\t^\ast_{\RR} = P\otimes\RR\). The information of these measures is an ``upgrade'' from that of the MV polytopes.

Denote by \(\cPP\) the subspace of distributions on \(\t_{\RR}^\ast\) spanned by 
linear combinations of piecewise-polynomial functions times Lebesgue measures on (not necessarily full-dimensional) polytopes.
%

Given \(\nu\in Q_+\) with \(p=\hgt\nu\), let \(\vi\in \Seq(\nu)\) and consider the map \(\pi_\vi : \RR^{p+1}\to \t_\RR^\ast\) that takes the \(k\)th standard basis vector to the \(k\)th partial sum \(\alpha_{i_1} + \cdots + \alpha_{i_k}\) for each \(1\le k\le p\). 
Let \(\delta_{\Delta^p}\) denote Lebesgue measure on the standard \(p\)-simplex \(\Delta^p = \{x_1 + \cdots + x_{p + 1} = 1\}\) in \(\RR^{p+1}\) and define \(D_\vi \) to be its pushforward along \(\pi_\vi\).

The measures \(D_\vi\) satisfy a ``shuffle'' product as follows. Given \(\vj\in \Seq(\nu')\) and \(\vk\in\Seq(\nu'')\), we have that
\begin{equation}
    \label{eq:shufd}
    D_\vj\ast D_\vk = \sum_{\vi\in \vj\shuffle \vk} D_\vi
\end{equation}
where \(\vj\shuffle \vk\) denotes the set of permutations \(\vi\) of \(\vj\sqcup \vk\) that maintain the same relative order among elements of \(\vj\) and \(\vk\).
Note, if \(\vj\) has length \(a\) and \(\vk \) has length \(b\), then \(\vj\shuffle\vk\) has \(\binom{a + b}{a}\) elements, each of which belongs to \(\Seq(\nu' + \nu'')\). 

Consider the free Lie algebra \(\f\) on the set \(\{ e_i\}_I\). Its universal enveloping algebra \(\cU(\f)\) is the free associative algebra on the set \(\{ e_i\}_I\). In particular, \(\cU(\f)\) is graded by \(Q_+\). Moreover, for each \(\nu\in Q_+\),
the set \(\{e_\vi\}_{\Seq(\nu)}\), with \(e_\vi\) denoting the element \(e_{i_1}\cdots e_{i_p}\), is a basis for \(\cU(\f)_{\nu}\). 

The algebra \(\cU(\f)\) is in fact a graded Hopf algebra with finite dimensional components, so its graded dual \((\cU(\f))^\ast\) is also a Hopf algebra. For each \(\nu\in Q_+\) we can consider the dual basis \(\{ e_\vi^\ast\}\) indexed by \(\vi \in \Seq(\nu)\) for \((\cU(\f))_{-\nu}^\ast\). From the definition of the coproduct on \(\cU(\f)\) we get the following ``shuffle identify'' in \((\cU(\f))^\ast\) 
\begin{equation}
    \label{eq:shufe}
     e_\vj^\ast e_\vk^\ast = \sum_{\vi\in\vj\shuffle \vk} e_\vi^\ast\,.
\end{equation}
There is a unique Hopf algebra map \(\cU(\f)\to\cU(\fu)\) that sends \(e_i\) to \(e_i\) for each \(i\in I\). The fact that the dual map is an inclusion of algebras \(\CC[U]\to (\cU(\f))^\ast\) implies that the map \(D: \CC[U]_{-\nu} \to \cPP\) defined by
\begin{equation}
    \label{eq:Df}
    D(f) = \sum_{\vi\in\Seq(\nu)} \langle e_\vi , f\rangle D_\vi  
\end{equation}
extends to 
an algebra map \(D : \CC[U] \to \cPP\).

We can think of \(D(b_Z)\) and \(D(c_Y)\) as measures on \(\Pol(Z)\) and \(\Pol(Y)\), respectively. Indeed, as shown in \cite{baumann2019mirkovic} and recalled below, (modulo our chosen isomorphism of \(\t_\RR^\ast \) and \(\t_\RR\),) the former is exactly the Duistermaat--Heckman measure of \(Z\), and both are supported on their associated polytopes.
%
\section{Measures from \texorpdfstring{\(\Gr\)}{the affine Grassmannian}}
\label{s:mfromGr}
To define the Duistermaat--Heckman measure of an MV cycle we need to identify characters of representations and measures.
Let \((X,A,V)\) be a triple as in the general setup of Section~\ref{s:xva}. Given a weight \(\mu\) of \(A\), 
let \(\delta_\mu\) denote the distribution on 
\(\fa_{\RR}^\ast\)
defined by \(\delta_\mu(f) = f(\mu)\) for any 
\(f\in C^\infty(\fa^\ast_{\RR})\otimes\CC\).
%
Let 
\(R(A)\)
denote the representation ring of 
\(A\) 
and consider the embedding \(R(A) \to \cPP \) which sends the class of an irreducible \(A\)-representation to the inverse Fourier transform of its formal character. (See Section~\ref{s:barD}).
\begin{equation}
    \label{eq:repmes}
    R(A) \to \cPP \qquad 
    [V]\mapsto\sum_{\mu} \dim V_\mu~ \delta_\mu \,.
\end{equation}
%
%
The final ingredient in this definition of the Duistermaat--Heckman measure is the scaling automorphism \(\tau_n: \mu\mapsto\frac\mu n\) of \(\fa^\ast_{\RR}\). 
%
\begin{definition}
    \label{def:dhxav} 
The Duistermaat--Heckman
measure of \((X,A,V)\) is defined by 
\begin{equation}
    \label{eq:dhxav}
    DH(X) = \lim_{n\to\infty} \frac{1}{n^{\dim X}}(\tau_n)_* [H^0(X,\cO(n))]\,.
\end{equation}
%
\end{definition}
One can check that each \(\tau_n^\ast[H^0(X,\cO(n)]\), and, hence, \(DH(X)\) is supported on \(\Pol(X)\) (c.f.\ \cite{brion1990action}).
For \((X,A,V) = (Z,T^\vee,L(\Lambda_0))\), this means that \(DH(Z)\) is supported on the corresponding MV polytope, the moment polytope of \(Z\). Together with the following theorem, this fact ensures that the measures coming from the MV basis \(B(\Gr)\) are \textit{at least} upgrades of MV polytopes, in the sense that taking the support of a measure \(D(b_Z)\) recovers the corresponding MV polytope \(\Pol(Z)\).
%
\begin{theorem}\label{thm:DDh}
\cite[Theorem 10.2]{baumann2019mirkovic} 
The measure associated to an element of the MV basis is equal to the DH measure of the corresponding MV cycle,
    \(D(b_Z) = DH(Z)\). 
\end{theorem}
%

\section{Extra-compatibility}
\label{s:exco}
Let \(Y\in\irr\Lambda \). In this section we will see that the measure \(D(c_Y)\) is also a limit supported on a polytope. In other words, the comparison of polytopes undertaken in \cite{baumann2014affine} can be enhanced to a comparison of sequences of measures and their limits.

Let \(F_n(M)\) be the space of \((n+1)\)-step flags of submodules of \(M\) and denote by \(F_{n,\mu}(M)\) the locus
\begin{equation}
    \label{eq:Fnmu}
F_{n,\mu}(M) = 
\{M = M^1 \supset M^2 \supset \cdots \supset M^n  \supset M^{n+1}= 0\big|\sum\dimvec M^k = \mu \}\,.
\end{equation}
Then, formally, 
\begin{equation}
    \label{eq:HFnM}
  [H^\bullet(F_n(M))] = \sum_\mu [H^\bullet(F_{n,\mu}(M))] \delta_{\mu} 
\end{equation}
and one can check that this measure is supported on \(n\Pol(M)\).

\begin{theorem}\label{thm:DHFn}
\cite[Theorem 11.4]{baumann2019mirkovic} 
Let \(M\) be a general point in \(Y\) of dimension \(\dimvec M = \nu\). Then
\begin{equation}
    \label{eq:DHFn}
D(c_Y) = \lim_{n\to\infty} \frac{1}{n^{\rho(\nu)}} (\tau_n)_\ast [H^\bullet(F_n(M))]\,.
\end{equation}
In particular, \(D(c_Y)\) is supported on \(\Pol(Y)\).
\end{theorem} 

Given \(Z\in\cZ(\infty)\) and \(Y\in\irr\Lambda\) having common polytope \(P\), we can thus ask whether \(D(c_Y) = D(b_Z)\), and, 
if so, we can ask whether the sequences converging to either measure are in some sense compatible. We study the strongest possible form of compatibility, dubbed ``extra-compatibility'' by the authors of \cite{baumann2019mirkovic}. 
The pair \((Z,Y)\) is called \new{extra-compatible} if \(\dim H^0(Z,\cO(n))_\mu = \dim H^\bullet (F_{n,\mu}(M))\) for all \(n\in\NN\) and for all \(\mu\in P\) whenever \(M\) is a general point in \(Y\). In Chapter~\ref{ch:calculs}, we present evidence of extra-compatibility in type \(A\).

\section{Fourier transform of a measure}
\label{s:barD}
The equivariant measures developed so far are as yet unwieldy, so in this section we recall the poor man's alternative and our equivariant invariant of choice: a coefficient on a Fourier transform of an equivariant measure.
%
\begin{definition}
    \label{def:ft}
The Fourier transform is the map that takes a measure \(\eta\in\cPP\) to the meromorphic function \(FT(\eta)\in\CC(\t)\) defined by 
\begin{equation}
    \label{eq:ft}
    FT(\eta)(x) = \int_{\t^\ast_\RR} e^{\langle\mu,x\rangle} d\eta(\mu)
\end{equation}
for any \(x\in\t\). 
\end{definition}
One can show that the Fourier transform is one to one, multiplicative on convolution of distributions, and satisfies \(FT(\delta_\lambda) = e^{\lambda}\) \cite[Lemma 8.6]{baumann2019mirkovic}.
Henceforth, we write \(e^\mu\) for the function \(x\mapsto e^{\langle\mu,x\rangle}\) and \(FT(\eta)\) for the function \(x\mapsto FT(\eta)(x)\).

Given \(\vi = (i_1,\dots,i_p) \) in \(\Seq(\nu)\), let \(\beta^\vi_j\) be the partial sum \(\alpha_{i_1} + \cdots + \alpha_{i_j}\), not to be confused with the positive root \(\beta^{\uvi}_k\) depending on a reduced word \(\uvi\) defined in Chapter~\ref{ch:bases}. 
\begin{lemma}\label{lem:ftD}
    \cite[Lemma 8.7]{baumann2019mirkovic}
\begin{equation}
    \label{eq:ftD}
    FT(D_\vi) = \sum_{j = 0}^p 
    \frac{e^{\beta^\vi_j}}{\prod_{k\ne j}(\beta^\vi_k - \beta^\vi_j)}\,.
\end{equation}
\end{lemma}
Let \(\barD : \CC[U] \to \CC(\t)\) be the map that takes \(f\in\CC[U]_{-\nu}\) to constant term (the coefficient of 
\(e^0\))
in the Fourier transform of \(D(f)\). 
By linearity,
\begin{equation}
    \label{eq:barDf}
    \barD(f) = \sum_{\vi\in\Seq(\nu)}\langle e_\vi,f\rangle \barD_\vi
\end{equation}
and, as we can check by inspecting the \(j=0\) summand of Equation~\ref{eq:ftD}, 
\begin{equation}
    \label{eq:barDi}
    \barD_\vi = \prod_{k=1}^p\frac{1}{\alpha_{i_1} + \cdots + \alpha_{i_k}}\,.
\end{equation}
To see that \(\barD(f)\) is a rational function for any \(f\in\CC[U]\), one can rewrite it as follows. 

Let \(x\in\t\). Then \(x\) is called \new{regular} if
\(\langle\alpha,x\rangle\ne 0\) for all \(\alpha\in\Delta\).
We denote by \(\treg\) the subset of regular elements of \(\t\) and distinguish the regular element \(e = \sum e_i\) that is the sum of root vectors of weight \(\alpha_i\).
%
\begin{theorem}\label{thm:xux}
    \cite[Proposition 8.4]{baumann2019mirkovic}  
For any \(x\in\treg\) there exists unique \(u_x\in U\) such that \(u_x x u_x^{-1} = x + e\). 
Moreover \(\barD\) can equally be defined as the pullback of the map 
\begin{equation}\label{eq:xux}
    \treg\to U\qquad x\mapsto u_x\,.
\end{equation}
Namely, \(\barD (f)(x) = f(u_x)\).
\end{theorem}
%
\section{Reinterpreting \texorpdfstring{\(\barD\)}{\(D\) data}}
\label{s:barDanew}
%
Let \(Y\in\irr\Lambda(\nu)\) and let \(M\) be a general point in \(Y\). Then, by definition of \(c_Y\), 
\begin{equation}
    \label{eq:barDc}
    \barD(c_Y) = \sum_{\vi\in\Seq(\nu)}\chi(F_\vi(M)) \barD_\vi\,
\end{equation}
and we refer to the right-hand side as the \new{flag function} of \(Y\). 
More interesting is what we get when we evaluate \(\barD\) on elements of the MV basis. 

Once again, let \((X,A,V)\) be a general triple as in Section~\ref{s:xva}. 
Let \(S\) be the multiplicative set in \(H^\bullet_A(V)\cong\CC[\fa]\) generated by nonzero weights of \(A\).
Then, the inclusion of fixed points \(X^A\to X\) induces the isomorphism \(S^{-1}H_\bullet^A(X^A) \to S^{-1}H_\bullet^A(X)\). In turn, the class of \(X\) acquires the unique fixed point expansion
\begin{equation}
    \label{eq:equivloc}
    [X] = \sum_{x\in X^A} \varepsilon^A_x(X) [\{x\}]
\end{equation}
over \(S^{-1}\CC[\fa]\). 
Following \cite{brion1997equivariant}, the coefficient \(\varepsilon^A_x(X)\) is called the equivariant multiplicity of \(X\) at \(x\). 

\begin{theorem}\label{thm:ftdhz}
    \cite[Theorem 9.6]{baumann2019mirkovic}
    Let \((X,A,V) = (Z,T^\vee,L(\Lambda_0))\) with \(Z\) an an MV cycle of weight \(\nu\).
    Then 
    \begin{equation}
    \label{eq:ftdhz}
    FT(D(b_Z)) = \sum_{L\in Z^{T^{\vee}}} \varepsilon_L^{T^{\vee}} (Z) e^{\Phi_{T^\vee}(L)}
    = \sum_{\beta\in \Pol(Z)} \varepsilon_{L_\beta}^{T^{\vee}}(Z) e^\beta
    \,.
    \end{equation}
In particular, 
    \begin{equation}
        \label{eq:barDb}
        \barD(b_Z) = \varepsilon_{L_{0}}^{T^\vee}(Z)\,.
    \end{equation}
\end{theorem}   
\section{Equivariant multiplicities via multidegrees}
\label{s:mdegrules}
In this section, we see that equivariant multiplicities of MV cycles at fixed points can be replaced by multidegrees of hyperplane sections. In turn, multidegrees of hyperplane sections can be replaced by multidegrees of generalized orbital varieties, the subject of the following chapter. 

Let \((X,A,V)\) be as in Section~\ref{s:xva}. 
Given \(p \in X^A\) we can choose a weight vector \(\alpha \in V^\ast \) such that \(\alpha(p) \ne 0\) and 
\(\{x\in X \big| \alpha(x)\ne 0\}\) is \(A \)-invariant.
We write \(\mathring X_p \) for \(X\cap \{\alpha\ne 0\}\) and \(W\) for \( \PP(V)\cap \{\alpha\ne 0\}\). 
Note that \(W\)
has dimension \(\dim V - 1\). 
Taking the projective closure of \(\mathring X_p\) 
we recover \(X\). 
%

The \new{multidegree} of the triple \((\mathring X_p, A, W)\)
is defined to be the unique element \(\mdeg_W(\mathring X_p)\) of \(H_A^\bullet(W)\cong\CC[\fa]\) such that \(\mdeg_W (\mathring X_p) [W] = [\mathring X_p]\) in \(H^\bullet_A(W)\).
%
In practice, it can be computed according to the following list of rules. (See \cite[Chapter 2, Lemma 3.2]{joseph1997orbital} or \cite{knutson2014brauer}.) 
Given a hyperplane \(H\subset W\)
\begin{itemize}
    \item \(\mdeg_W(\mathring X_p) = \mdeg_{\mathring X_p\cap H} (\mathring X_p\cap H) \) if \(\mathring X_p\not\subset H\) 
    \item \(\mdeg_W(\mathring X_p) = \mdeg_{W\cap H} (\mathring X_p)\left(\text{weight of }T\text{ on }W/H\right)\) if \(\mathring X_p\subset H\)
    \item \(\mdeg_W(\mathring X_p) = \sum_{C \in \irr \mathring X_p}\mdeg_W(C)\) 
\end{itemize}

\begin{proposition}\label{prop:epX}
    \cite[Proposition 9.5]{baumann2019mirkovic} If \((X,A,V)\) and \((\mathring X_p, A, W)\) are as above, then 
    \(\mdeg_W(\{p\}) \) is the product of weights of \(A\) on \(W\), and
    we have the following equation in \(H^\bullet_A(W)\cong \CC[\fa]\). 
\begin{equation}
    \label{eq:epX}
    \varepsilon^A_p (X) = \frac{\mdeg_W(\mathring X_p)}{\mdeg_W(\{p\})}\,.
\end{equation} 
\end{proposition}
Let \(v_{\Lambda_0}^\ast\) be such that \(v_{\Lambda_0}^\ast(v) = 1 \) if \(v\) is the highest weight vector \(v_{\Lambda_0}\), and \(v_{\Lambda_0}^\ast(v) = 0\) if \(v\) is any other weight vector. In other words, \(v_{\Lambda_0}^\ast\) is projection onto the highest weight space of \(L(\Lambda_0) \). The intersection \(W = \PP(L(\Lambda_0)) \cap \{v_{\Lambda_0}^\ast \ne 0\} \) is an affine space. Moreover, if \(Z\) is a stable MV cycle embedded in \(\PP(L(\Lambda_0))\), then \(\mathring Z = Z \cap \{v_{\Lambda_0}^\ast \ne 0\} \) is an open neighbourhood of \(L_0\in Z\). 

Proposition~\ref{prop:epX} will allow us to compute \(\varepsilon_{L_{0}}^{T^\vee}(Z)\) using the Mirkovi\'c--Vybornov isomorphism in terms of the multidegree of the triple \((\mathring Z, T^\vee, W)\). 
%
But first, we shall pin down \(\mathring Z\) by studying generalized orbital varieties.
%

\chapter[Generalized orbital varieties]{Generalized orbital varieties for Mirkovi\'c--Vybornov slices}
\label{ch:govs}

In this chapter we give a Spaltenstein type decomposition of a subvariety of the Mirkovi\'c--Vybornov slice that lends itself to coordinatizing MV cycles in type \(A\). 

\section{Setup and context}
\label{s:setup}
Fix \(m\ge 0\). 
We say that \(\mu = (\mu_1,\dots,\mu_m)\in\ZZ^m\) is \new{effective} if \(\mu_i\ge 0\) for all \(1\le i\le m\).
%
Recall the notation \(\cY(m)\) for partitions having at most \(m\) parts.
%
Fix a positive integer \(N\ge m\) and denote by \(\cY(m)_N\) the subset of partitions of (size) \(N\). 
\begin{equation}
    \label{eq:effdom}
    \cY(m)_N = 
    \{\lambda = (\lambda_1\ge \lambda_2\ge\cdots\ge\lambda_m\ge 0) \big|\sum_1^m\lambda_i = N\}\,.
\end{equation}
Given \(\lambda\in \cY(m)_N\) we denote by \(\OO_\lambda \) the conjugacy class of \(N\times N\) matrices having Jordan type \(\lambda\). 

In \cite{mirkovic2003quiver}, the authors define a family of slices \(\TT_\mu\) to $\overline\OO_\lambda$ 
through the Jordan normal forms \(J_\mu\) associated to effective dominant weights \(\mu\in \cY(m)_N\), with the additional property that \(\mu_m\ne 0\).
%
They show that the intersection \(\overline\OO_\lambda\cap\TT_\mu\) is an irreducible affine variety isomorphic to a subvariety of the affine Grassmannian. 
Let \(\n\) denote the algebra of uppertriangular \(N\times N\) matrices. 

In this chapter, towards refining this isomorphism of Mirkovi\'c and Vybornov, we investigate the half-dimensional subvariety \(\ot\cap\n\) of \(\overline\OO_\lambda\cap\TT_\mu\).
Our first result is that the set of irreducible components \(\irr (\ot\cap\n)\) is in bijection with semistandard Young tableaux of shape \(\lambda \) and weight \(\mu\).
%
%

Recall that given \(\tau\in\cT(\lambda)_\mu\) we denote by \(\lambda^{(i)}_\tau\) the shape of the tableau \(\tau^{(i)}\) obtained from \(\tau \in \cT(\lambda)_\mu\) by deleting all entries \(j > i\).
Given an \(N\times N\) matrix \(A\in\OO_\lambda\) we denote by \(A_{(i)}\) the upper left \(\lvert\lambda^{(i)}\rvert\times\lvert\lambda^{(i)}\rvert\) submatrix of \(A\), with \(\lvert\lambda^{(i)}\rvert = \sum_{j=1}^i \lambda^{(i)}_j\).
The GT pattern \((\lambda^{(i)})\) of \(\tau\) is used to define a matrix variety
\(\mathring X_\tau\) as follows.
\begin{equation}\label{eq:blockyXt}
\mathring X_\tau = \left\{
A \in \TT_\mu\cap\n \big|
A_{(i)} \in \OO_{\lambda^{(i)}}
\text{ for } 1\le i\le m
\right\} \,.
\end{equation}
%
Our first result (Theorem~\ref{thm:ad1}) is that by taking closures of top dimensional components of \(\mathring X_\tau\), as \(\tau\) varies over the set of semistandard Young tableaux of shape \(\lambda \) and weight \(\mu\), we get a complete list of 
irreducible components of the transverse intersection $\ot\cap\n$.
Part of this claim appears in \cite{zinn2015quiver} where it is stated without proof.

We call the closure of the top dimensional component of \(\mathring X_\tau\) a \new{generalized orbital variety} for the Mirkovi\'c--Vybornov slice \(\TT_\mu\). 
The reason for this terminology is that when \(\tau\) is a \textit{standard} Young tableau, i.e.\ weight \(\mu = (1,1,\dots,1)\). 
In this case \(\TT_\mu = \gl_N\) and the decomposition 
\begin{equation}
    \label{eq:spadec} 
    \overline \OO_\lambda\cap\n = \bigcup_{\sigma \in \cT(\lambda)_{\mu}} X_\sigma
\end{equation}
is a special case 
of a correspondence established in \cite[Corollary 9.14]{joseph1984variety} which relates (ordinary) orbital varieties and Spaltenstein's \cite{spaltenstein1976fixed} decomposition of the fixed point set of a unipotent transformation on a flag variety. 
In particular, the varieties \(X_\sigma\) recover the ordinary orbital varieties, a priori defined to be 
the irreducible components of \(\overline\OO_\lambda \cap \n\).
In fact, conventions differ about whether ordinary orbital varieties are subvarieties of \(\n\) that are dense in a component of \(\OO_\lambda\cap\n\) or are the closed irreducible components of \(\overline{\OO}_\lambda\cap\n\). 
These issues will not be relevant for us and we will take generalized orbital varieties to be closed.
%

\section{Gaining traction}
\label{s:gains}
To give a more tractable description of the sets defined by Equation~\ref{eq:blockyXt}, we will need additional notation.
Fix \(\tau\in\cT(\lambda)_\mu\) and 
\((i,k)\in\{1,2,\ldots,m\}\times\{1,2,\ldots,\mu_i\}\). 
We denote by \(\lambda^{(i,k)}_\tau\) (respectively, by \(\mu^{(i,k)}_\tau\)) the shape (respectively, the weight) of the tableau \(\tau^{(i,k)}\) obtained from \(\tau\) by deleting all \(j>i\) and all but the first \(k\) occurrences of \(i\). 
Here, repeated entries of a tableau are ordered from left to right, so that the first occurrence of a given entry is its leftmost.

We make the identifications \(\lambda^{(i)}_\tau\equiv \lambda^{(i,\mu_i)}_\tau\), \(\mu^{(i)}_\tau \equiv \mu^{(i,\mu_i)}_\tau\) and \(\tau^{(i)}\equiv\tau^{(i,\mu_i)}\) for all \(i\). 
We also do away with the subscript \(\tau\) as (we never work with two tableaux at a time, so) it is always clear from from the context. 
\begin{example}
    \label{eg:tableauxnotns}
If we take 
\[
    \tau = \young(112,23)
\] 
as before, then we find that  
\[
    \tau^{(2)} = \young(112,2)
\] 
has shape \(\lambda^{(2)} = (3,1)\) and weight \(\mu^{(2)} = (2,2)\), 
while 
\[
    \tau^{(2,1)} = \young(11,2)
\] 
has shape \(\lambda^{(2,1)} = (2,1)\) and weight \(\mu^{(2,1)} = (2,1)\).
\end{example}

Now let \((e^1_1,\ldots,e^{\mu_1}_1,\ldots,e^1_m,\ldots,e^{\mu_m}_m)\)
be an enumeration of the standard basis of \(\CC^N\) which we call the \(\mu\)-numeration so as to have something to refer to later. 
Following \cite[\S 1.2]{mirkovic2003quiver}, we define the Mirkovi\'c--Vybornov slice \(\TT_\mu\) through the Jordan normal form \(J_\mu\) to be the set of \(A\in \gl_N\) such that
\begin{equation}
    \label{eq:mvyslice}
        \begin{aligned}
            &\text{for all } 1 \le a,s\le m\,,
            \text{for all } 1\le b\le \mu_a\,, 1\le t\le \mu_s\,, \\
            &\text{ if } 1\le t < \mu_s \text{ or } t = \mu_s < b \le \mu_a\,,
            \text{ then } (e^t_s)' (A-J_\mu) e^b_a = 0 \,.
        \end{aligned}
\end{equation}
Here \((e^t_s)'\) denotes the transpose of the column vector \((e^t_s)\). 
%
\begin{example}
    \label{eg:mvyslice423}
The weight of the tableau from the last example, \(\mu = (2,2,1)\), leads us to the \(\mu\)-numeration \((e_1^1,e_1^2,e_2^1,e_2^2,e_3^1)\) of \(\CC^5\), and elements of \(\TT_\mu\) in this basis take the block form
\[
    \left[\begin{BMAT}(e){ccccc}{ccccc}
        0 & 1 & 0 & 0 & 0 \\
        * & * & * & * & * \\
        0 & 0 & 0 & 1 & 0 \\
        * & * & * & * & * \\
        * & 0 & * & 0 & *
        \addpath{(2,5,:)ddrrddr}
        \addpath{(4,5,:)ddr} 
        \addpath{(0,3,:)rrddrrd}
        \addpath{(0,1,.)rrd}
    \end{BMAT}\right]    
\]
with $\ast$s denoting unconstrained entries. Thus, for instance, 
\[
\tau = \young(112,23)\,\Longrightarrow\,
\mathring X_\tau = 
\left\{
\left[\begin{BMAT}(e){ccccc}{ccccc}
    0 & 1 & 0 & 0 & 0 \\
    0 & 0 & a & b & c \\
    0 & 0 & 0 & 1 & 0 \\
    0 & 0 & 0 & 0 & d \\
    0 & 0 & 0 & 0 & 0
    \addpath{(2,5,:)ddrrddr}
    \addpath{(4,5,:)ddr} 
    \addpath{(0,3,:)rrddrrd}
    \addpath{(0,1,.)rrd}
\end{BMAT}\right] \text{such that } a,d = 0\text{ and } b,c\ne 0
\right\}\,.
\]
%
\end{example}
While due to \cite{mirkovic2003quiver}, the definition of this slice is motivated (and, as we see in the proof of Theorem~\ref{thm:mvy}, elucidated) by the lattice description of $\Gr$ dating at least as far back as \cite{lusztig1981green}.

Given \((i,k)\in \{1,2,\ldots,m\}\times\{1,2,\ldots,\mu_i\}\), we denote by \(V^{(i,k)}\) the span of the first \(\lvert\mu^{(i,k)}\rvert = \mu_1 + \cdots + \mu_{i-1} + k\) basis vectors of \(\CC^N\) in the \(\mu\)-numeration.
Similarly, we denote by \(A_{(i,k)}\) the restriction \(A\big|_{V^{(i,k)}}\) viewed as the top left \(\lvert\lambda^{(i)}\rvert \times \lvert\lambda^{(i)}\rvert \) submatrix of \(A\). We are justified in viewing the restriction of \(A\) as a submatrix since we will be working with upper triangular \(A\).
Finally, as we did for tableaux, we make the identification \(V^{(i)}\equiv V^{(i,\mu_i)}\) and \(A_{(i,\mu_i)}\equiv A_{(i)}\) for all \(1\le i\le m\).

\section{A boxy description of \texorpdfstring{\(\mathring X_\tau\)}{generalized orbital varieties}}
\label{s:aboxy}
In this section we show that the block-by-block restrictions
\(\{A_{(i)}\text{ such that }1\le i\le m\}\) of Equation~\ref{eq:blockyXt}, coming from deleting boxes of \(\tau\) one ``weight-block'' at a time, can be replaced by the box-by-box restrictions 
\(\{A_{(i,k)} \text{ such that } 1\le k\le\mu_i,1\le i\le m\}\) 
coming from deleting boxes of \(\tau\) one by one.
\begin{example}
    \label{eg:boxy423}
    Continuing with the previous example (Example~\ref{eg:mvyslice423}), we see that subset of \(\TT_{(2,2,1)}\cap\n\) determined by the conditions \(A_{(1)}\subset\OO_{(2)}\), \(A_{(2)}\subset\OO_{(3,1)}\), and \(A\subset\OO_{(3,2)}\) is unchanged if we add that \(A_{(1,1)}\in\OO_{(1)}\) and \(A_{(2,1)}\in\OO_{(2,1)}\). 
\end{example}
\begin{lemma}
    \label{lem:impliesboxyalg}
Let \(B\) be an \((N-1)\times(N-1)\) matrix of the form
\[
    \left[\begin{BMAT}(e){cc}{cc}
            C & v \\
            0 & 0 
    \end{BMAT}\right]
\]
for some \((N-2)\times(N-2)\) matrix \(C\) and column vector \(v\). 
Let \(A \) be an \(N\times N\) matrix of the form 
\[
    \left[\begin{BMAT}(e){ccc}{ccc}
        C & v & w \\
        0 & 0 & 1 \\ 
        0 & 0 & 0 
    \end{BMAT}\right]
\]
for some column vector \(w\).
Let \(p\ge 2\). If \(\rk C^p < \rk B^p  \), then \(\rk B^p < \rk A^p\). 
\end{lemma}
\begin{proof}
Let 
\[
    B = 
\left[\begin{BMAT}(e){cc}{cc}
        C & v \\
        0 & 0 
\end{BMAT}\right]
\]
and let
\[
A = 
\left[\begin{BMAT}(e){ccc}{ccc}
    C & v & w \\
    0 & 0 & 1 \\
    0 & 0 & 0 
\end{BMAT}\right] 
= 
\left[\begin{BMAT}(@){c:c}{c:c}
    B & \begin{BMAT}{c}{cc}
        w \\ 1
    \end{BMAT} \\
    \begin{BMAT}{cc}{c}
        0 & 0
    \end{BMAT} & 0 
\end{BMAT}\right]\,.
\]
Suppose \(\rk B^p > \rk C^p \) for \(p\ge 0\). Clearly \(\rk A^p > \rk B^p\) for \(p = 0,1\) independent of the assumption. Suppose \(p\ge 2\). 
Since 
\[
B^p = 
\left[\begin{BMAT}(@){cc}{cc}
    C^p & C^{p-1} v \\
        0 & 0 
\end{BMAT}\right]
\]
this means \(C^{p-1} v \not\in\Im C^p \). So \(C^{p-2} v\not\in \Im C^{p-1}\) and \(C^{p-2} v + C^{p-1} w \not\in \Im C^{p-1}\).
Since 
\[
A^p = 
\left[\begin{BMAT}(@){ccc}{ccc}
    C^p & C^{p-1} v & C^{p-1} w + C^{p-2} v\\
    0   & 0         & 0 \\
    0   & 0         & 0 
\end{BMAT}\right]
\] 
it follows that \(\rk A^p > \rk B^p \) as desired. 
\end{proof}

Now fix \(\lambda,\mu\) and \(\tau \in \cT(\lambda)_\mu\). Let \(A\in \mathring X_\tau\). 
Recall that \(V^{(i,k)}\) denotes the span of the first \(\mu_1+\cdots + \mu_{i-1} + k\) vectors of the \(\mu\)-numeration \((e^1_1,\ldots,e^{\mu_1}_1,\ldots,e^1_m,\ldots,e^{\mu_m}_m)\).
\begin{lemma}
    \label{boxyalg} 
    \(A_{(m,\mu_m-1)}\in\OO_{\lambda^{(m,\mu_m - 1)}}\).
\end{lemma}
\begin{proof}
Let 
\(B = A_{(m,\mu_m-1)}\). 
Assume \(\mu_m > 1\) or else there is nothing to show. Let \(C = A_{(m,\mu_m-2)}\).

By definition of \(\mathring X_\tau\), \(A\in \OO_\lambda\). 
Let \(\lambda(B)\) denote the Jordan type of \(B\) and \(\lambda(C)\) the Jordan type of \(C\).
Since \(\dim V/V^{(m-1)} = \mu_m \) is exactly the number of boxes by which \(\lambda\) and \(\lambda^{(m-1)}\) differ, \(\lambda(B)\) must contain one less box than \(\lambda\), and \(\lambda(C)\) must contain one less box than \(\lambda(B)\).
Let \(c(A)\) denote the column coordinate of the box by which \(\lambda\) and \(\lambda(B)\) differ, and let \(c(B) \) denote the column coordinate of the box by which \(\lambda(B)\) and \(\lambda(C)\) differ.
Then
\[
\rk B^p - \rk C^p = \begin{cases}
    1 & p < c(B) \\
    0 & p \ge c(B)
\end{cases}\,,
\]
so we can apply Lemma~\ref{lem:impliesboxyalg} to our choice of \((A,B,C)\) to conclude that \(\rk A^p > \rk B^p \) for \(p < c(B) \). At the same time,
\[
\rk A^p - \rk B^p = \begin{cases}
    1 & p < c(A) \\
    0 & p \ge c(A)
\end{cases}
\]
implies that \(c(A) > c(B) \). We conclude that \(B \in \OO_{\lambda^{(m,\mu_m - 1)}}\) as desired. 
\end{proof}

Thus we see that the \(\mu\times\mu\) block rank conditions defining \(\mathring X_\tau\) in Equation~\ref{eq:blockyXt} refine a possibly redundant set of \(a\times a\) box rank conditions for each \(1\le a\le N\). 
\begin{proposition}
\label{prop:boxyXt}
\begin{equation}\label{eq:boxyXt} 
    \mathring X_\tau = \left\{
        A\in\TT_\mu\cap\n \text{ such that } 
        A_{(i,k)}\in\OO_{\lambda^{(i,k)}} \text{ for } 1\le k\le \mu_i \text{ and } 1\le i\le m
    \right\}\,.
\end{equation}
\end{proposition}
\begin{proof}
The non-obvious direction of containment is an immediate consequence of Lemma~\ref{boxyalg}. 
\end{proof}
%
%
\section{Irreducibility of \texorpdfstring{\(X_\tau\)}{generalized orbital varieties}}
\label{s:irred}
Set \(\tau-\raisebox{-0.5ex}{\young(m)} \equiv \tau^{(m,\mu_m-1)}\), and let \(r\) equal to the row coordinate of the last \(m\), aka the row coordinate of the box by which \(\tau\) and \(\tau-\raisebox{-0.5ex}{\young(m)}\) differ.
Let \(L=\Sp(e^{\mu_1}_1,e^{\mu_2}_2,\ldots,e^{\mu_{m-1}}_{m-1})\).
\begin{lemma}
    \label{lem:sall}
    Let \(B\in \mathring X_{\tau-\raisebox{-0.5ex}{\scriptsize\young(m)}}\) and let \(S = {(B^{\lambda_r - 1})}^{-1}\Im B^{\lambda_r}\). Then
\begin{equation}
    \label{eq:sall}
    V^{(m,\mu_m-1)} = S + L\,.
\end{equation}
\end{lemma} 
\begin{proof}
Let \(1\le a\le m\) and \(1\le b\le \mu_a\).
If \(b\le\mu_a - 1 \) then 
\(e_a^b = B (e_a^{b+1}) - v\) for some \(v\in L\). 
Note that
\[
    B(e_a^{b+1}) \subset (B^{c-1})^{-1}\Im B^{c}
\]
for any \(c\).
Therefore \(e_a^b \in (B^{c-1})^{-1}\Im B^{c} + L \) 
unless of course for \((a,b) = (m,\mu_m)\).
\end{proof}
\begin{lemma}
\label{lem:altfibdes}
Let \(B\in \mathring X_{\tau-\raisebox{-0.5ex}{\scriptsize\young(m)}}\) and let \(S = {(B^{\lambda_r - 1})}^{-1}\Im B^{\lambda_r}\). %
Then
\begin{subequations}
    \begin{align}
        \label{eq:dim1}
        \dim S &= N - r\,,\text{ and} \\
        \label{eq:dim2}
        \dim S\cap L &= m - r\,.
    \end{align}
\end{subequations}
\end{lemma}
\begin{proof}
By Lemma~\ref{boxyalg}, 
\(B\)
has Jordan type \(\lambda^{(m,\mu_m - 1)}\) 
which differs from \(\lambda\) by a single box in position 
\((r,\lambda_r)\).
Denote by 
\[
    (f_1^1,\ldots,f^{\lambda_1}_1,\ldots,f_r^1,\ldots,f_r^{\lambda_r - 1},\ldots,f_\ell^1,\ldots,f_\ell^{\lambda_\ell})
\] 
a Jordan basis for \(B\).
Then 
\[
\begin{aligned}
    \Im B^{\lambda_r - 1}  
    &= \Sp(\{f_c^1 ,\ldots, f_c^{\lambda_c - \lambda_r + 1}\}_{1\le c\le\ell}) \\
    &= \Sp(\{f_c^1 ,\ldots, f_c^{\lambda_c - \lambda_r} \}_{1\le c\le\ell})
    + \Sp(\{f_c^{\lambda_c - \lambda_r + 1}\}_{1\le c\le\ell}) \\
    &= \Im B^{\lambda_r} + \Sp(\{f_c^{\lambda_c - \lambda_r + 1} \}_{1\le c\le\ell})\,,
\end{aligned}
\]
with \(f_c^p \equiv 0\) if \(p\le 0\). In particular, \(f_c^{\lambda_c - \lambda_r + 1}\) is equal to \(B^{\lambda_r - 1}(f_c^{\lambda_c})\) and is nonzero if \(c > r\). 
Thus \(\dim {(B^{\lambda_r - 1})}^{-1}\Im B^{\lambda_r} = N - r\).

By Lemma~\ref{lem:sall}, \(\dim (S + L) = N-1\). Therefore 
\begin{equation}
    \begin{split}
        \dim S\cap L &= \dim S + \dim L - \dim (S + L) \\
                     &= (N-r) + (m-1) - (N-1) \\ 
                     &= m-r\,.
    \end{split}
\end{equation}
\end{proof}
\begin{lemma}
\label{lem:irrfib}
The map
\begin{equation}
    \label{eq:irrfib}
    \mathring X_\tau\to \mathring X_{\tau - \raisebox{-0.5ex}{\scriptsize\young(m)}}\quad A \mapsto 
    A_{(m,\mu_m - 1)}
\end{equation}
has irreducible fibres of dimension \(m - r\).
\end{lemma}
\begin{proof}
Let \(B \in \mathring X_{\tau-\raisebox{-0.5ex}{\scriptsize\young(m)}}\) and let \(F_B\) denote the fibre over \(B\).
Assume \(\mu_m > 1\). Then \(A\in F_B\) takes the form 
\[
\left[\begin{BMAT}(@){c:c}{c:c}
    B & v + e^{\mu_m - 1}_m \\
    0 & 0
\end{BMAT}\right] 
\]
for some \(v\in L \). 
In other words, \(A\) is determined by some \(v\in L\). 

Let \(u\in\Ker A^{\lambda_r} \setminus (\Ker A^{\lambda_r - 1} \cup \Ker B^{\lambda_r})\). Of course the first \(\lambda_r-1\) columns of \(\lambda\) have less boxes than the first \(\lambda_r\) columns. In addition, the \(\lambda_r\)th column, by definition of \(r\), contains the ``last'' \(m\). So \(\dim \Ker B^{\lambda_r}\) and \(\dim \Ker A^{\lambda_r}\) differ by 1 (the 1 box at the end of row \(r\)), while \(\dim \Ker A^{\lambda_r}\) and \(\dim \Ker A^{\lambda_r -1}\) differ by \(r\) (the length of the \(\lambda_r\)th column, which includes the box at the end of row \(r\)).
Suppose without loss of generality that \(u = e^{\mu_m}_m + w\) for some \(w\in V^{(m,\mu_m-1)}\) where a coefficient of zero on \(e^{\mu_m}\) would result in \(u\) actually contributing to \(\Ker B^{\lambda_r}\) and a nonzero coefficient can be normalized.
Then 
\begin{equation}
    \label{eq:preimB}
    0 = A^{\lambda_r} (u) 
      = A^{\lambda_r} ( e^{\mu_m}_m + w ) 
      = B^{\lambda_r-1} (v + e^{\mu_m - 1}_m) + B^{\lambda_r} (w) \,.
\end{equation}
Let \(S =  {(B^{\lambda_r-1})}^{-1}(\Im B^{\lambda_r})\) and \(S' = (B^{\lambda_r - 2})^{-1}(\Im B^{\lambda_r - 1})\). 
Equation~\ref{eq:preimB} implies that 
\(v + e^{\mu_m - 1}_m \in S \), while the assumption \(u\not\in\Ker A^{\lambda_r - 1}\) implies \(v + e^{\mu_m}_m \not\in S'\).
These are the only conditions on \(v\), making \((S\setminus S' + e^{\mu_m -1}_m)\cap L\) the total space of permissible \(v\).

Lemma~\ref{lem:sall} lets us rewrite \(e^{\mu_m - 1}_m = s + l\) for some \(s \in S\setminus S'\) and some \(l\in L\). 
Moreover, a choice of such a pair \((s,l)\) provides an isomorphism
\begin{equation}\label{eq:sliso}
    \begin{split}
        (S\setminus S' + e^{\mu_m -1}_m)\cap L = (S + l) \cap L  &\to S\setminus S' \cap L \\
        v &\mapsto v - l\,.
    \end{split}
\end{equation}

The map
\begin{equation}
    \label{eq:slfib}
    \begin{split}
        (S\setminus S') \cap L &\to F_B \\
        c &\mapsto 
        \left[
            \begin{BMAT}(@){c:c}{c:c}
                B & c + s \\
                0 & 0
            \end{BMAT}
        \right]
    \end{split}
\end{equation}
is then also an isomorphism. Since \((S\setminus S')\cap L\) is locally closed, \(F_B\) is irreducible. 
Lemma~\ref{lem:altfibdes}, Equation~\ref{eq:dim2} then applies to give the desired dimension count.
\end{proof}
\section{Irreducibility of generalized orbital varieties}
\label{s:irrs3} 
%
We would like to use Lemma~\ref{lem:irrfib} to describe \(\irr(\ot\cap\n)\).
First, we will need the following proposition.
\begin{proposition}
\label{prop:myirrcl} 
Let \(X\) and \(Y\) be varieties. Let \(f:X\to Y\) be surjective, with irreducible fibres of dimension \(d\). Assume \(Y\) has a single component of dimension \(m\) and all other components of \(Y\) have smaller dimension. Then \(X\) has unique component of dimension \(m + d\) and all other components of \(X\) have smaller dimension.
\end{proposition}
To prove Proposition~\ref{prop:myirrcl} 
we will use %
\cite[I, \S 8, Theorem 2]{mumford1988red} and \cite[Lemma 005K]{stacks-project} which are recalled below.
\begin{theorem}
\label{thm:mumf2}
\cite[I, \S 8, Theorem 2]{mumford1988red}
Let \(f:X\to Y\) be a dominating morphism of varieties and let \(r = \dim X - \dim Y\). Then there exists a nonempty open \(U\subset Y \) such that:
\begin{enumerate}[label={(\roman*)}]
    \item \label{mump1} \(U\subset f(X)\) 
    \item \label{mump2} for all irreducible closed subsets \(W\subset Y\) such that \(W\cap U \ne \varnothing\), and for all components \(Z\) of \(f^{-1}(W)\) such that \(Z\cap f^{-1}(U)\ne\varnothing \), \(\dim Z = \dim W + r\). 
\end{enumerate}
\end{theorem}
\begin{lemma}
\label{lem:generic-point-in-constructible}
\cite[Lemma 005K]{stacks-project}
Let $X$ be a topological space. Suppose that
$Z \subset X$ is irreducible. Let $E \subset X$
be a finite union of locally closed subsets (i.e.\ $E$ is constructible). %
The following are equivalent
\begin{enumerate}
\item \label{stacks1} The intersection $E \cap Z$ contains an open
dense subset of $Z$.
\item The intersection $E \cap Z$ is dense in $Z$.
\end{enumerate}
\end{lemma}
\begin{proof}[Proof of Proposition~\ref{prop:myirrcl}]
Let \(X = \cup_{\irr X} C\) be a (finite) decomposition of \(X\). Consider the restriction \(f\big|_C : C \to \overline{f(C)}\) of \(f\) to an arbitrary component. It is a dominating morphism of varieties, with irreducible fibres of dimension \(d\).

We apply Theorem~\ref{thm:mumf2}. Let \(U\subset \overline{f(C)}\) be an open subset satisfying properties~\ref{mump1} and \ref{mump2}. 
Then, taking \(W = \{y\} \subset U\) for some \(y\in U\subset \overline{f(C)}\), we get that \(\dim f^{-1}(y) = \dim C - \dim \overline{f(C)}\). Since all fibres have dimension \(d\), the difference \(\dim C - \dim \overline{f(C)}\) is constant and equal to \(d\), independent of the component we're in. 

Since \(f\) is surjective, it is in particular dominant, so we have that 
\[
    Y = f(X) = f(\cup_{\irr X} C) = \cup_{\irr X} f(C) = \overline{\cup_{\irr X} f(C)}= \cup_{\irr X} \overline{f(C)}\,.
\]

We now show that no two irreducible components can have \(m\)-dimensional image.
Suppose towards a contradiction that \(C_1\ne C_2 \in \irr X\)
are two components of dimension \(m + d\). 
Since \(\dim C_i = d + \dim \overline{f(C_i)} = d + m\) for both \(i = 1,2\), \(\dim \overline{f(C_1)} = \dim\overline{f(C_2)} = m\), so \(\overline{f(C_1)} = \overline{f(C_2)}\). 
For each \(i = 1,2\), denote by \(f_i\) the restriction \(f\big|_{C_i}\), and let \(U_i \subset \overline{f_i(C_i)} \) be the open sets supplied by Theorem~\ref{thm:mumf2} 
(equivalently, by Lemma~\ref{lem:generic-point-in-constructible})
for the constructible sets \(E_i = f(C_i)\). 
Fix \( y\in U = U_1\cap U_2\), a dense open subset of \(\overline{f(C_1)}= \overline{f(C_2)}\). 
Since \(V_i = f_i^{-1}(U)=f^{-1}(U) \cap C_i \) contains \(f_i^{-1}(y) = f^{-1}(y)\cap C_i = f^{-1}(y)\) for both \(i = 1 ,2\), the set \(V = V_1\cap V_2\) is nonempty and open. Moreover, \(V\) is contained in \(C_1\cap C_2\), proving that \(C_1 = C_2\).
\end{proof}
With Lemma~\ref{lem:irrfib} and Proposition~\ref{prop:myirrcl} in hand, we are ready to state and prove the first part of our first main theorem, which says that our matrix varieties are irreducible in top dimension. 

For each \(1\le i\le m\), let 
\(\rho^{(i)} = (i,i-1,\dots,1)\).
%
%
Set \(\rho \equiv \rho^{(m)}\), and, unambiguously, denote by \(\langle\phantom{0},\phantom{0}\rangle\) the standard inner product on \(\ZZ^i\) for any \(i\). 
\begin{proposition}
\label{prop:Xtirr}
    \(\mathring X_\tau\) has one irreducible component \(\mathring X^{\text{top}}_\tau\) of maximum dimension \(\langle \lambda - \mu, \rho\rangle\).
\end{proposition}
\begin{proof}
Consider the restriction map
\[
    \mathring X_\tau \to \mathring X_{\tau^{(m-1)}} \,.
\]
By induction on \(m\), we can assume that \(\mathring X_{\tau^{(m-1)}}\) has one irreducible component of dimension 
\(d = \langle \lambda^{(m-1)} - \mu^{(m-1)}, \rho^{(m-1)}\rangle\),
and apply Proposition~\ref{prop:myirrcl} in conjunction with Lemma~\ref{lem:irrfib} to conclude that \(\mathring X_\tau\) has one irreducible component of dimension \(d+\sum_1^{\mu_m}(m-r_{m,k})\) where \(r_{m,k}\) is equal to the row coordinate of the \(k\)th \(m\) in \(\tau\). Note \(\mathring X_{\tau^{(1)}} = \{J_{\mu_1}\}\).

The reader can check that 
\[
    \sum_{k=1}^{\mu_m} (m-r_{m,k}) = \langle \lambda-\mu,\rho\rangle 
    - d\,. 
\]
\end{proof}
It follows that the generalized orbital variety \(X_\tau = \overline{\mathring X_\tau^{\text{top}}}\) is irreducible. It is also implicit in the proof of Proposition~\ref{prop:Xtirr} that there is a rational map \(\pi_b : X_\tau \to X_{\tau^{(b)}}\) for all \(1\le b\le m\). 
%
\section{Conjectured recursion}
\label{s:recursion}
In this section, we state the conjecture that \(\mathring X_\tau\) is itself irreducible. We use the fact that \(\overline{A\times B} = \overline{A}\times\overline{B}\) in the Zariski topology.
\begin{conjecture}\label{conj:rec}
The map in Lemma~\ref{lem:irrfib} is a trivial fibration. Consequently \(\mathring X_\tau^{\text{top}} = \mathring X_\tau\) and \(X_\tau = \overline{\mathring X_\tau} \) is defined by a recurrence \(X_\tau \cong X_{\tau-\raisebox{-0.5ex}{\scriptsize\young(m)}} \times \CC^{m-r}\) for \(r\) equal to the row coordinate of the last \(m\) in \(\tau\). 
%
\end{conjecture}
\begin{example}
    \label{eg:trivfibev}
    Let \(\tau = \young(12,34)\) so $m = 4$ and $r = 2$. 
    Then \(A\in \mathring X_\tau\) takes the form 
    \[
    \left[
    \begin{BMAT}(@){cccc}{cccc}
        0 & a & b & c \\
        0 & 0 & 0 & e \\
        0 & 0 & 0 & f \\ 
        0 & 0 & 0 & 0
    \end{BMAT}
    \right]
    \]
    with \(a,f\ne 0\) and \(ae + bf = 0 \). Since \(e = -bf/a\), i.e.\ the defining equations completely determine the coordinate \(e\) in the fibre over a point of 
    \(\mathring X_{\tau-\raisebox{-0.5ex}{\scriptsize\young(4)}}\), 
    the isomorphism
    \[
        \mathring X_\tau \to \mathring X_{\tau-\raisebox{-0.5ex}{\scriptsize\young(4)}} \times \CC\times\CC^\times \qquad A \mapsto (A\big|_{V^{(3)}}, (c,f))
    \] supplies the desired trivialization.
\end{example}
To state our conjecture more precisely we need to define an order \(<_\tau\) on the set \([m-1]=\{1,2,\dots,m-1\}\). We do so by the following 
tableau sort 
function. This function takes a tableau as input and outputs a permutation of the alphabet \([m-1]\) which we interpret as an ordering of the alphabet from largest to smallest induced by the tableau.
\begin{lstlisting}[language=Python]
    def sort(t):            
    MT = max(t)[0]
    L = []
    for i in range(1,MT):
        s = t.restrict(i)   
        MS,(r,c) = max(s)
        L.insert(r,MS)   
    return L
\end{lstlisting}
It depends on the {\texttt{max}} function
\begin{lstlisting}[language=Python]
    def max(t):                
    ME = -1,(0,0)
    for c in t.corners():
        if t.entry(c) > ME[0]:
            ME = t.entry(c),c
    return ME
\end{lstlisting}
which once again takes a tableau as input and outputs the $($value, position$)$ pair of the rightmost largest corner entry. Here, the position of an entry is its $($row, column$)$ coordinate, with the top left entry of a tableau anchored in position \((0,0)\). The resulting order is best illustrated by example.
\begin{example}
    The tableau 
    \[\young(113,22,4)\]
    has \texttt{maxentry} \texttt{4,(2,0)}. 
    The resulting list \texttt{L} is \texttt{3,1,2}. 
    The tableau \[\young(12,34)\] has \texttt{maxentry} \texttt{4,(1,1)}. 
    The resulting list \texttt{L} is \texttt{2,3,1}. 
\end{example}
We claim that the smallest \(m-r\) elements of the ordered set \(\left([m-1],<_\tau\right)\) correspond precisely to the unconstrained entries in the last column of a generic point \(A\in \mathring X_\tau\). In the case of \[\young(12,34)\] this agrees with the computation of Example~\ref{eg:trivfibev}: the first and third entries of the last column are unconstrained. 
\begin{conjecture}
    \label{conj:rec2}
Let \((A_1,\dots,A_{m-1})\) be coordinates on \(X_\tau\) given by the a priori nonzero entries in the last column of \(A\in\TT_\mu\). Let \(i_1,\dots,i_{m-r}\) be the \(m-r\) smallest elements of \(\left([m-1],<_\tau\right)\). Then the isomorphism \(X_\tau\cong X_{\tau-\raisebox{-0.5ex}{\scriptsize\young(m)}}\times\CC^{m-r}\) is defined by extending the map 
\[
    \begin{aligned}
        \mathring X_\tau &\to \mathring X_{\tau-\raisebox{-0.5ex}{\scriptsize\young(m)}}\times\CC^{m-r} \\
        A &\mapsto (A_{(m-1)}, (A_{i_1},\dots, A_{i_{m-r}}))\,.
    \end{aligned}
\]
\end{conjecture}

We intend to pursue this in the future, especially insofar as it can help us to (efficiently) deduce the equations, multidegrees, or initial ideals of generalized orbital varieties. 

Let \(P = \tau^{(m-1)}\). We conjecture that the permutation achieved by our sort applied to a standard tableau \(\tau\) is precisely the permutation associated to the pair \((P,P)\) by Viennot's geometric construction (RSK). In particular, this permutation is always an involution. 
Let \(\tilde\tau\) be the tableau obtained from \(\tau\) by deleting all but the max entries and performing a Jeu de Taquin to standardize the result in case it is skew.

When \(\tau \) is semistandard, the permutation achieved by our sort is that associated to the pair \((\tilde P,\tilde P)\).
Let \(JT\) denote the Jeu de Taquin function. Below is an example of the commuting operations just described.
\begin{example}\label{eg:rsk}
    Let 
    \[
        \tau = \young(113,22,4)\,.
    \] 
    Then our sort function results in the permutation \(\pi = {123\atop 213} \). Viennot's geometric construction applied to \(\pi\) produces a standard Young tableau as follows (c.f.\ \cite[Section 3.6]{sagan2013symmetric}).
\begin{equation}\label{eq:rsk}
    {123\atop 213} \quad \overset{\text{Viennot}}{\longmapsto} \quad 
    \begin{tikzpicture}[scale=0.5,baseline={([yshift=-.5ex]current bounding box.center)}]
        \draw[OliveGreen,thick] (1,4) -- (1,2) -- (2,2) -- (2,1) -- (4,1);
        \draw[red,thick] (2,4) -- (2,2) -- (4,2);
        \draw[OliveGreen,thick] (3,4) -- (3,3) -- (4,3);
        \foreach \x in {1,2,3} 
        {
            \node at (\x,-0.3) {\x};
            \node at (-0.3,\x) {\x};
        }
        \draw[help lines] (0,0) grid (3,3);
    \end{tikzpicture} \quad 
    \overset{\text{RSK}}{\longmapsto}\quad 
    \ytableausetup{nosmalltableaux}
    \begin{ytableau}
        *(OliveGreen!20!white) 1& *(OliveGreen!20!white)3 \\
        *(red!20!white) 2
    \end{ytableau}\,.
\end{equation}
This agrees with
\begin{equation}
    \widetilde{\young(113,22,4)}
    =
    JT\left(\young(:13,:2,4)\right) 
    = \young(13,2,4)
\end{equation}
\end{example}
\noindent up to the maximum entry of our input tableau. 

Note that the tableau 
\[
    \young(13,24)
\] 
for instance would give the same permutation. More generally, a feature of \texttt{sort} is that it depends only on a certain standard skeleton of a given semistandard tableau. In addition, it does not depend on (the position occupied by the largest entry) \(m\). 

We cite private communication with D.\ Muthiah in conjecturing furthermore that a distinguished facet in the Stanley--Reisner description of the initial ideal of \(X_\tau\) is given by the union of the \(i-r_i\) unconstrained entries over all \(1\le i\le m\).

\section{Equations of generalized orbital varieties}
\label{s:eqgov}
Here we record our method for finding equations with the help of the computer algebra system Macaulay2.
Given a tableau \(\tau\in \cT(\lambda)_\mu\) we initiate a matrix \(A\in\TT_\mu\cap\n\) and consider its \(\lvert\mu^{(i)}\rvert\times\lvert\mu^{(i)}\rvert\) submatrices \(A_i\) for all \(1\le i\le m\). 

The first nontrivial restriction on the entries of \(A\) may come from the shape \(\lambda^{(2)}\) of \(\tau^{(2)}\) as follows. Recall that this shape is recording the Jordan type of the \((\mu_1 + \mu_2)\times(\mu_1 + \mu_2)\) submatrix \(A_2\). This means that the number of boxes \(h_c^{(2)}\) in the first \(c\) columns of \(\tau^{(2)}\) is equal to the nullity of \(A_2^c\). Equivalently, the number of boxes \(r^{(2)}_c \) to the right of the \(c\)th column is equal to the rank of \(A_2^c\) so all \((r^{(2)}_c + 1)\times(r^{(2)}_c + 1)\) minors of \(A_2^c\) must vanish.

The ideal of our matrix variety is thus built up by iteratively adding equations for \((r^{(i)}_c + 1)\times(r^{(i)}_c + 1)\) minors of \(A_i^c\). We have to take care to discard any prime factors that hold on \(\overline\OO_\lambda\) but not on \(\OO_\lambda\) at every step. This is accomplished either by inspecting the output of a primary decomposition, or by preemptively performing a rip operation on (known) irrelevant substrata. 

For example, let 
\[
    \tau = \young(12,3)
\] 
and take
\[
    A = 
    \left[
        \begin{BMAT}(@){ccc}{ccc}
        0 & a & c \\
        & 0 & b \\
        &   & 0
        \end{BMAT}
    \right]
\] 
in \(\mathring X_\tau\). Clearly the ideal of \(2\times 2\) minors is \((ab)\). Since \(\rk A_2 = 1\) forces \(a\ne 0\), we must colon out ideal \((a)\). The desired variety is therefore described by the ideal \((ab):(a) = (b)\).
The colon operation is defined on pairs of ideals \(I,J\) by 
\[
    I:J = \{f\big| fJ \subset I\} 
\]
and loosely corresponds to ripping the hypersurface where \(J\) vanishes.

In the appendix we detail this method as applied to our counterexample. Below we take it up on a smaller example.
\begin{example}\label{eg:eqgov}
Let 
\[
    \tau = \young(113,22,4)\,.
\] 
Then \(A\in \mathring X_\tau\) takes the form
\[
    \left[
        \begin{BMAT}(@){cc:cc:c:c}{cc:cc:c:c}
            0 & 1 & 0 & 0 & 0 & 0 \\
            0 & 0 & a_1 & a_2 & a_3 & a_4 \\
            0 & 0 & 0 & 1 & 0 & 0 \\
            0 & 0 & 0 & 0 & a_5 & a_6 \\
            0 & 0 & 0 & 0 & 0 & a_7 \\
            0 & 0 & 0 & 0 & 0 & 0 
        \end{BMAT}
    \right]
\]
and
\[
A^2 = 
    \left[
        \begin{BMAT}(@){cc:cc:c:c}{cc:cc:c:c}
            0 & 0 & a_1 & a_2 & a_3 & a_4 \\
            0 & 0 & 0 & a_1 & a_2a_5 & a_2a_6 + a_3 a_7 \\
            0 & 0 & 0 & 0 & a_5 & a_6 \\
            0 & 0 & 0 & 0 & 0 & a_5a_7 \\
            0 & 0 & 0 & 0 & 0 & 0 \\
            0 & 0 & 0 & 0 & 0 & 0 
        \end{BMAT}
    \right]
    \quad 
    A^3 =
    \left[
        \begin{BMAT}(@){cc:cc:c:c}{cc:cc:c:c}
            0 & 0 & 0 & a_1 & a_2 a_5 & a_2a_6 + a_3a_7 \\
            0 & 0 & 0 & 0 & a_1a_5 & a_2a_5a_7 + a_1a_6 \\
            0 & 0 & 0 & 0 & 0 & a_5a_7 \\
            0 & 0 & 0 & 0 & 0 & 0 \\
            0 & 0 & 0 & 0 & 0 & 0 \\
            0 & 0 & 0 & 0 & 0 & 0 
        \end{BMAT}
    \right]
\]
Since \(\rk A_2^2 = 0\), all \(1\times 1\) minors, \(a_1,a_2\), vanish. Vanishing of \(3\times 3\) minors of \(A_2\) coming from \(\rk A_2 = 2\) introduces no new vanishing conditions. All \(\rk A_3^c\) conditions are vacuous. Finally \(\rk A^3 = 0\) imposes \(a_3 a_7 = a_5 a_7 = 0 \). However, \(\rk A_3 = 3\) does not allow both \(A_4,A_6 = 0\). So, upon coloning out, we end up with \(a_7 = 0\). The condition \(\rk A^2 = 1\) leaves just one \(2\times 2\) minor to consider, it is \(a_3 a_6 - a_4 a_5\). 
Note $\rk A_k^3 = 0$ for all $k$.
\end{example}

\section{\texorpdfstring{Decomposing \(\ot\cap\n\)}{Spaltenstein type decomposition}}
\label{s:decom}
\begin{lemma}
    \label{lem:impliesboxyalgconverse}
Let \(\lambda\ge\mu\) be dominant effective of size \(N\). Then 
\[
    \OO_\lambda\cap\TT_\mu\cap\n = \bigcup_{\tau} \mathring X_\tau
\]
\end{lemma}
\begin{proof}
Let \(A\in\OO_\lambda\cap\TT_\mu\cap\n\) and \(B = A_{(m,\mu_m-1)}\). Suppose first that \(\mu_m > 1\). 
It suffices to show that \(A \in \mathring X_\tau\) for some \(\tau\in\cT(\lambda)_\mu\).
As in Lemma~\ref{lem:impliesboxyalg}, write
\[
    A = 
    \left[
        \begin{BMAT}(e){ccc}{ccc}
            C & v & w \\
            0 & 0 & 1 \\ 
            0 & 0 & 0 
        \end{BMAT}
    \right] 
\] 
and let
\[
    B = 
    \left[
        \begin{BMAT}(e){cc}{cc}
            C & v \\
            0 & 0 
        \end{BMAT}
    \right]\,.
\]
Denote by \(\lambda(B)\) and \(\lambda(C)\) the Jordan types of \(B\) and \(C\) respectively. 
Suppose \(p'\) and \(p''\) are such that \(\lambda\) and \(\lambda(B)\) differ by a box in column \(p'\) and \(\lambda(B) \) and \(\lambda(C)\) differ by a box in column \(p''\). 
Let us show that \(p'' < p' \). 
Indeed 
\begin{align}
    \rk B^{p'} &= \rk A^{p'} - 1  \\
    \rk C^{p''} &= \rk B^{p''} - 1\,. \label{cpbp}
\end{align}
By Lemma~\ref{lem:impliesboxyalg}, Equation~\ref{cpbp} implies that \(\rk B^{p''}<\rk A^{p''}\). 
The assumption that \(\lambda(B)\) and \(\lambda(A)\) differ in column \(p'\) (and column \(p'\) only) means that \(\lambda(B)\) and \(\lambda(A)\) have columns of equal length to the right of column \(p'\). 
In particular, \(\rk A^{p} = \rk B^{p}\) for all \(p > p'\).
By induction on $m$, and $\mu_m$ it suffices to show 
$p'\ne p''$ when $\mu_m = 2$. 

If $m = 1$ and $\mu_m = 2$, then by definition of $\TT_\mu$ and $\n$, $A = \begin{bmatrix} 0 & 1 \\ 0 & 0 \end{bmatrix}$, so $\tau = \young(11)$, which is semistandard. 

Suppose $m > 1$, $\mu_m = 2$, and $A \in \OO_\lambda\cap \TT_\mu\cap \n$ is as above, i.e. 
\[
    A = 
    \left[
        \begin{BMAT}(e){ccc}{ccc}
            C & v & w \\
            0 & 0 & 1 \\ 
            0 & 0 & 0 
        \end{BMAT}
    \right]\,.
\] 
Recall,
\[
    A^k = 
    \left[
        \begin{BMAT}(@){ccc}{ccc}
            C^k & C^{k-1} v & C^{k-1} w + C^{k-2} v\\
            0   & 0         & 0 \\
            0   & 0         & 0 
        \end{BMAT}
    \right]\,.
\] 

Suppose (towards a contradiction) that $p' = p'' = p$, and consider \(A^p \) and \(A^{p-1}\).
Since the number of boxes that are to the right of column $p$ in each of $\lambda$, $\lambda(B)$, and $\lambda(C)$ is the same, $\rk C^p = \rk B^p = \rk A^p$.
Looking at \(A^p\), this implies in particular that $C^{p-1} w + C^{p-2} v \in \Im C^p$. So $C^{p-2} v$ is in the image of $C^{p-1}$. 

On the other hand, $\rk C^{p - 1} = \rk B^{p-1} - 1 = \rk A^{p-1} - 2$. 
So neither $C^{p-2} v$ nor $C^{p-2} w + C^{p-3} v$ can be in $\Im C^{p-1}$. This contradicts the previous conclusion that $C^{p-2} v$ is in the image of $C^{p-1}$. 
%

It follows that \(p'' < p'\). 
Let \(\tau\) be a tableau of shape \(\lambda\) and weight \(\mu\) such that \(\lambda^{(m,\mu_m-1)} = \lambda(B)\) and \(\lambda^{(m,\mu_m-2)} = \lambda(C)\). 
Iterate the argument on \(B\), filling the remaining \(\mu_m - 2\) boxes of \(\lambda\) with \(m\). 
If \(\mu_m = 1\), then there is no ambiguity in filling the box where \(\lambda(B)\) and \(\lambda\) differ with \(m\).
Repeat this (two-case) argument on \(A_{(m-1)}\) and so forth to determine \(\tau\) completely. 

To conclude, the reader can check that \(A\in\mathring X_\tau\) against the boxy description of Equation~\ref{eq:boxyXt}.
\end{proof}
\begin{theorem}
\label{thm:ad1}
The map \(\tau\mapsto X_\tau\) is a bijection \(\cT(\lambda)_\mu\to\irr(\ot \cap \n)\). 
\end{theorem}
\begin{proof}
By Corollary~\ref{cor:resmvy}, \(\ot\cap\n\cong\overline{\Gr^\lambda}\cap S^\mu_-\). 
By \cite[Theorem~3.2]{mirkovic2007geometric}, the latter intersection is pure dimensional of dimension \(\langle \lambda - \mu,\rho\rangle\). It follows that \(\ot\cap\n\) must be pure dimensional also.  

Let \(Y\in\irr(\ot\cap\n)\). Because of the pure-dimensionality $Y$ must intersect $\OO_\lambda$ and that intersection must be dense in $Y$. Together with Lemma~\ref{lem:impliesboxyalgconverse}, pure-dimensionality also implies that \(Y\cap\mathring X_\tau\) is dense in $Y$ for some \(\tau\), and so \(Y = X_\tau \). 

Alternatively, one can wield the fact that the MV cycles of type \(\lambda\) and weight \(\mu\) give a basis of \(L(\lambda)_\mu\) whose dimension is well-known to be equal to the number of semi-standard Young tableaux of shape \(\lambda\) and weight \(\mu\).
\end{proof} 

\section{Applications and relation to other work}
\label{s:appli}
\subsection{Big Spaltenstein fibres}
\label{ss:bigsp}
%
Given an ordered partition \(\nu\vdash N\), let \(P_\nu\subset GL_N\) be the parabolic subgroup of block upper triangular matrices with blocks of size \(\nu\)
and denote by \(\p_\nu\) its Lie algebra. We'll view elements of the partial flag variety \(X_\nu := GL_N/P_\nu\) interchangeably as parabolic subalgebras of \( \gl_N \) which are conjugate to \(\p_\nu \) and as flags 
\[
    0 = V_0\subset V_1 \subset \cdots \subset V_{(\nu')_1} = \CC^N 
\] 
such that \(\dim V_i/V_{i-1} = \nu_i\). 
Here \(\nu'\) denotes the conjugate partition of \(\nu\).

For \(u-1\in\OO_\lambda\) fixed, Shimomura, in \cite{shimomura1980theorem}, establishes a bijection between components of  the \new{big Spaltenstein fibre} \(X_\mu^u\) and \(\cT(\lambda)_\mu\), generalizing Spaltenstein's decomposition in \cite{spaltenstein1976fixed} for the case \(\mu = (1,\ldots,1)\) and implying that big Spaltenstein fibres also have the same number of top-dimensional components as \(\ot\cap\n\). 

We conjecture that the coincidence is evidence of a correspondence implying a bijection between generalized orbital varieties and top-dimensional irreducible components of \(\OO_\lambda\cap\p_\mu\).
More precisely, let \(\cN\) denote the set of nilpotent matrices in \(\gl_N\).
Let 
\[ 
    \widetilde{\g}_\mu = 
    \{ 
        (A,V_\bullet)\in \cN \times X_\mu \big| AV_i\subset V_i \text{ for } i = 1,\ldots, \mu'_1
    \}\,.
\]
Another description is 
\[
    \widetilde{\g}_\mu = \{(A,\mathfrak p)\in \cN \times X_\mu \big| A\in \mathfrak p\}\,.
\]
Let 
\(A = u-1 \in \OO_\lambda\) and consider the restriction of 
\(\pr_1 : \widetilde{\g}_\mu\to\cN\) defined by \(\pr_1(A,\p) = A\) to 
\(\widetilde{\g}_\mu^\lambda = \pr_1^{-1}(\OO_\lambda)\). 
We conjecture that the (resulting) diagram
\begin{equation}
    \label{eq:ginzburg}
    \begin{tikzcd}[column sep = scriptsize, row sep = scriptsize]
        & \OO_\lambda \cap \p_\mu \ar[r,mapsto] \ar[d,hook] & \{\p_\mu\} \ar[d,hook] \\ 
        X_\mu^u \ar[d] \ar[r,hook] & \widetilde{\g}_\mu^\lambda \ar[r]\ar[d] & X_\mu \\
        \{A\} \ar[r,hook] & \OO_\lambda
    \end{tikzcd}    
\end{equation}
has an orbit-fibre duality (generalizing that which is established in \cite[\S 6.5]{chriss2009representation} in the case \(\mu = (1,\ldots,1)\) and \(\OO_\lambda\cap\p_\mu = \OO_\lambda\cap\n\)) such that the maps
\[
    \OO_\lambda\cap\p_\mu \to \tilde \g^{\lambda}_\mu \leftarrow X_\mu^u
\]
give a bijection of top dimensional irreducible components,
\begin{equation}
    \label{eq:shimomura}
    \irr(\OO_\lambda\cap\p_\mu) \to \irr (X_\mu^u)\,.
\end{equation}
Combining the bijection \(\irr (X_\mu^u) \to \cT(\lambda)_\mu\) coming from Shimomura's generalized decomposition, with our bijection
$\irr (\OO_\lambda\cap\TT_\mu\cap\n)\to \cT(\lambda)_\mu$,
we get a bijection 
\begin{equation}
    \label{eq:shispa}
    \irr(\OO_\lambda\cap\p_\mu)\to \irr (\OO_\lambda\cap\TT_\mu\cap\n)\,.
\end{equation}
Since \(\OO_\lambda\cap\TT_\mu\cap\n\subset\OO_\lambda\cap\p_\mu\), this suggests a bijection of \(\irr(\OO_\lambda\cap\p_\mu)\) and \(\irr(\OO_\lambda\cap\TT_\mu\cap\n)\) such that
\[
    X \in \irr(\OO_\lambda\cap\p_\mu) \Longrightarrow X\cap\TT_\mu\cap\n \in \irr(\OO_\lambda\cap\TT_\mu\cap\n)\,.
\]
\subsection{Symplectic duality of little Spaltenstein fibres} %
\label{ss:sympl}
Recall the parabolic analogue of the Grothendieck--Springer resolution 
\begin{equation}
    \label{eq:gsp}
    T^\ast X_\mu\cong 
    \{
        (A, V_\bullet) \in \overline\OO_{\mu'} \times X_{\mu} \big| AV_i\subset V_{i-1}\text{ for all } i = 1,\ldots,(\mu')_1
    \} 
    \xrightarrow{\pi_\mu} \overline\OO_{\mu'} 
\end{equation}
given by projection onto the first component. 
Denote by \(X^{\lambda'}_\mu\) 
the preimage \(\pi^{-1}_\mu (\overline{\OO}_{\mu'}\cap \TT_{\lambda'})\) and by \(\pi^{\lambda'}_\mu \) 
the restriction \(\pi_\mu\big|_{X^{\lambda'}_\mu}\). 

Let 
\(X_0 = \ot\). %
View \(\lambda\) and \(\mu \) as weights of \(G = GL_m\), i.e.\ as partitions having \(m\) parts, and denote by \(T\subset G\) the maximal torus. 
Choose a coweight \(\rho:\CC^\times\to T\) such that
\[
    \{x\in X_0\big|\lim_{t\to 0} \rho(t)\cdot x = J_\mu\} = \ot\cap\n
\]
and denote this attracting set by \(Y_0\). 
Let \(\text{top} = \dim (\OO_\lambda\cap\TT_\mu)\).
\begin{lemma}\label{lem:blpw}
    \cite[Lemma~7.22, Example~7.25]{braden2014quantizations}
\(H_{\text{top}}(Y_0) \cong IH^{\text{top}}_T(X_0)\). 
\end{lemma}
Let 
\(F^! = (\pi^{\lambda'}_\mu)^{-1}( J_{\lambda'})\).
\begin{theorem}\label{thm:blpw}
    \cite[Theorem~10.4]{braden2014quantizations}
    \(X^{\mu}_{\lambda'}\) and \(X^{\lambda'}_\mu\) are symplectic dual.
    In particular, \( IH^{\text{top}}_T(X_0) \cong H^{\text{top}}(F^!) \).
\end{theorem}
Combining these results, we get that
\begin{equation}
    \label{eq:ft1}
    H^{\text{top}}(F^!)
    \cong 
    H_{\text{top}}(Y_0)\,.
\end{equation}
On the other hand, by \cite[Lemma 2.6]{rosso2012classic},
the irreducible components of the \new{little Spaltenstein fibre} 
\(F^!\),
which is equidimensional by \cite{ginzburg2009lectures}, are in bijection with \(\cT(\lambda)_\mu\). 
Thus, our bijection is predicted by and can be viewed as symplectic dual to that of Rosso.

\chapter{Equations for Mirkovi\'c--Vilonen cycles}
\label{ch:mvs}
%
%
Fix \(G = GL_m\). In this case, \(G^\vee = G\), so we can identify \(T^\vee = T\) and \(P^\vee = P\). 
Let \(\Gr = G(\cK)/G(\cO)\).
In this chapter we use the Mirkovi\'c--Vybornov isomorphism to show that 
MV cycles can be obtained as projective closures of generalized orbital varieties (or, homogenizations of their ideals). 
We start by showing that the decomposition of the last chapter is compatible with the crystal structure on MV cycles in the sense that the Lusztig datum of an MV cycle tells us exactly which generalized orbital variety is embedded into it under the Mirkovi\'c--Vybornov isomorphism.
We exploit the explicit description of a generalized orbital variety (coming from its tableau) to get an explicit description of the projective embedding of an MV cycle into \(\PP(L(\Lambda_0))\). 
Via the embedding, we can exhibit the equivariant multiplicity of a stable MV cycle in terms of the multidegree of a generalized orbital variety, accomplishing what we set out to do following the general setup of Proposition~\ref{prop:epX}. 
\section{Lattice model for the affine Grassmannian}
\label{s:lattices}
In type \(A\) we have at our disposal the \(\cO\)-lattice model for \(\Gr\).
Recall that \(L\) is an \(\cO\)-lattice if it is a free submodule of the vector space \(\cK^m\) with the property that \(\cK\otimes_{\cO} L\cong\cK^m\). 
Let \(L_0\) denote the standard lattice \(\cO^m\) and denote by \(\lat\) the set all \(\cO\)-lattices. 
\(\lat\) carries a natural action of \(GL_m(\cK)\) for which the stabilizer of any given lattice \(L\cong L_0\) is \(GL(L)\cong GL_m(\cO)\). We can therefore apply the orbit-stabilizer theorem to see that the map
\begin{equation}
    \label{eq:lattices}
    \Gr\to\lat \qquad gG(\cO) \mapsto gL_0
\end{equation}
is an isomorphism. 
Henceforth, we identify \(\Gr \) and \(\lat\) (abandoning the latter notation).

Recall that set-theoretically \(\Gr^T \cong P\).
Since \(P \cong \ZZ^m\), the partitions of the previous chapter 
can be viewed as weights of \(T\) and, hence, as points in \(\Gr\).
%
%
If \(\lambda\) is effective, then \(L_\lambda = t^\lambda G(\cO)\) is generated by \(t^{\lambda_i}e_i\) over \(\cO\). 
%
Henceforth, we fix \(N\ge m\) and consider effective partitions \(\lambda\in P\) of size \(N\). The size of a partition \(\lambda\) determines the connected component of \(\Gr\) to which \(L_\lambda\) belongs. 

\section{The Mirkovi\'c--Vybornov isomorphism}
\label{s:mviso}
Consider the quotient map \(G(\cK)\to\Gr\). 
Given a subgroup \(H\subset G\), we denote by \(H_1[t^{-1}]\) the kernel of the evaluation at \(t^{-1} = 0 \) map \(H[t^{-1}]\to H\). 
Given \(\mu\le\lambda\) in the dominance order on \(P\), let \(\Gr_\mu\) denote the orbit \(G_1[t^{-1}]\cdot L_\mu\) in \(\Gr\). Note, \(\lvert\mu\rvert = N\) and let \(\cN \) denote the cone of nilpotent \(N\times N\) matrices.
Up to a transpose, the following version of the Mirkovi\'c--Vybornov isomorphism is due to Cautis and Kamnitzer. 
\begin{theorem}
\label{thm:mvy}
\cite[Theorem 3.2]{cautis2018categorical}, \cite{mirkovic2019comparison}
The map 
\begin{equation}
    \label{eq:mvy1}
    \tilde\Phi : \TT_\mu\cap\cN\to G_1[t^{-1}]t^\mu
\end{equation}
defined by 
\[
\begin{aligned}
    \tilde\Phi(A) &= t^\mu + a(t) \\
    a_{ij}(t) &= - \sum A_{ij}^k t^{k-1} \\
    A_{ij}^k &= k\textrm{th entry from the left of the }\mu_j\times\mu_i\textrm{ block of }A
\end{aligned}
\]
gives the Mirkovi\'c--Vybornov isomorphism of type \(\lambda\) and weight \(\mu\)
\begin{equation}
    \label{eq:mvy2}
    \Phi : \ot\to \overline{\Gr^\lambda}\cap\Gr_\mu\qquad A \mapsto \tilde\Phi(A)L_0\,.
\end{equation}
\end{theorem}
\begin{proof}
Let \((e_it^j\big| 1\le i\le m, j\ge 0)\) 
be a basis of \(\CC^m\otimes\CC[t]\). It induces a basis in any finite-dimensional quotient of \(L_0\). In particular, it induces a basis in the \(N\)-dimensional quotient \(L_0/gL_0  \).

One can check that the map taking \(gG(\cO)\) to the transpose of the matrix of left multiplication by \(t\) on the quotient \(L_0/gL_0\), 
\[
    gG(\cO) \longmapsto 
    \left[
        t\big|_{L_0/gL_0}
    \right]'\,,
\]
in the induced basis \(([e_1],[te_1],\dots,[t^{\mu_1-1}e_m],\dots,[e_m],[te_m],\dots,[t^{\mu_m-1}e_m])\) is a two-sided inverse of \(\Phi\). 
\end{proof}
\begin{example}
\label{eg:mvyimage}
Let's see what happens to a generalized orbital variety under this isomorphism. Take
\[
    \tau = \young(112,23)
\] 
as in Example~\ref{eg:mvyslice423}. 
Then
\[
    \Phi(X_\tau) = 
    \left\{ gL_0 \big| g =
        \left[
            \begin{BMAT}(@){ccc}{ccc}
                t^2   & 0   & 0 \\
                -a-bt & t^2 & 0 \\
                -c    & -d  & t
            \end{BMAT}
        \right]
        \text{for some } a,b,c,d\in\CC\text{ with } 
        a,d = 0
    \right\}\,.
\]
\end{example}
We denote by \(\phi\) the restriction of \(\Phi\) to \(\ot\cap\n\). 
\begin{corollary}
\label{cor:resmvy}
Let \(\mu\) be dominant. Then \(\phi \) is an isomorphism of \(\ot\cap \n\) and \(\overline{\Gr^\lambda}\cap S^\mu_-\).
\end{corollary}
\begin{proof}
Let \(A\in \TT_\mu \cap\n\). 
Then \( \tilde\Phi(A) \in (U_-)_1[t^{-1}]t^\mu\).
Since \((U_-)_1[t^{-1}]t^\mu\subset U_-(\cK)t^\mu\) it follows that \(\Phi(A) \in S^\mu_-\). 
If also \(A\in\overline{\OO}_\lambda\), then \(\Phi(A) \in \overline{\Gr^\lambda}\cap S^\mu_-\).

Conversely, let \(L\in S^\mu_-\), and take \(g\in U_-(\cK)t^\mu\) such that \(L = g L_0 \).
Since \(U_-\) is unipotent we can write \(U_-(\cK)\) as the product \((U_-)_1[t^{-1}]U(\cO)\). Suppose \(g = g_1g_2t^\mu \) for some \(g_1 \in (U_-)_1[t^{-1}]\) and some \(g_2\in U(\cO)\). Now, since \(\mu\) is dominant, there exists \(g_3\in U(\cO)\) such that \(g_2t^\mu = t^\mu g_3\). 
So, \(g L_0 = g_1 t^\mu L_0\), and we can trade in our representative \(g\) for \(g_1t^\mu \). Since \((U_-)_1[t^{-1}]t^\mu \subset G_1[t^{-1}]t^\mu\) it follows that \(S^{\mu}_-\subset\Gr_\mu\), so \(\overline{\Gr^\lambda}\cap S^\mu_-\subset \overline{\Gr^\lambda}\cap \Gr_\mu\).
Since \(\Phi\) is onto, so is \(\phi\). 
%
%
\end{proof}
%
%
%
\section{Lusztig and MV cycles}
\label{s:ldmv}
%
MV polytopes for \(GL_m\) can be defined via BZ data, just as in Chapter~\ref{ch:bases}, Definition~\ref{def:bzdata}. We need only to enlarge the set of chamber weights for \(GL_m\) to include the ``all 1s'' weight \(\omega_m = (1,1,\dots,1)\). 
%
Since the vertices \((\mu_w)_W\) of the MV polytope \(P(M_\bullet)\) associated to a BZ datum of weight \((\lambda,\mu)\) by Equation~\ref{eq:polbz} satisfy the ``additional'' constraint \(\langle\mu_w,\omega_m\rangle = M_{\omega_m}\), \(P(M_\bullet)\) is contained in a hyperplane of the weight lattice. 

Hyperplanes correspond to connected components, and the MV polytope of a an MV cycle of weight \((\lambda,\mu)\) will have \(M_{\omega_m} = \lvert \lambda \rvert = \lvert \mu \rvert = N\). 
The MV polytopes for \(GL_m\) satisfying \(M_{\omega_m} = 0\) are thus the MV polytopes for \(SL_m\), and we can always translate an MV polytope for \(GL_m\) to make this equation hold. 
Moreover, the stable MV polytopes are exactly the same for \(SL_m\), \(PGL_m\), and \(GL_m\).

We can define the Lusztig datum of an MV cycle in terms of the BZ datum of its MV polytope. 
Or, we can define the Lusztig datum directly using Kamnitzer's constructible functions. 
Let \(V\) be a finite dimensional complex vector space, and consider  
the valuation map
\begin{equation}
  \label{eq:val}
  \begin{split}
\val: V\otimes\cK&\to \ZZ \\
v&\mapsto k \text{ if }v\in V\otimes t^k \cO \setminus V\otimes t^{k+1} \cO\,.
  \end{split}
\end{equation}
When \(V = \CC\) we see that \(\val \) records the lowest power on an element of \(\cK\) which enjoys the filtration 
\[
  \cdots \subset t^2\CC\xt \subset t\CC\xt \subset \CC\xt \subset t^{-1}\CC\xt\subset\cdots
\]
by lowest degree of indeterminate \(t\). 
%
Note that if \(V \) is a representation of \(G\) then \(V\otimes\cK \) is naturally a representation of \(G(\cK)\). 

Fix a chamber weight \(\gamma = w \omega_i \) for some \(1\le i\le m\) and for some \(w\in W\). Let \(v_{\omega_i}\) be a highest weight vector in the \(i\)th fundamental irreducible representation \(L(\omega_i) = \bigwedge^i\CC^m\) of \(G\), (e.g.\ \(v_{\omega_i} = e_1\wedge\cdots\wedge e_i\) will do,) and consider the map 
\begin{equation}
  \label{eq:dg}
\begin{split}
    D_{w\omega_i} : \Gr &\to \ZZ \\
    gG(\cO) &\mapsto \val (g 'v_{w_0w^{-1}\omega_i})\,. 
\end{split}
\end{equation}
Let \(J \) denote the subset of \([m] = \{1,2,\dots,m\}\) corresponding to the nonzero entries of \(w_0w^{-1}\omega_i \). Then \(g'v_{w_0w^{-1}\omega_i}\) is the wedge product of those columns of \(g'\) which are indexed by \(J\), and 
\begin{equation}
    \label{eq:valdet}
    \val( g' v_{w_0w^{-1}\omega_i}) = \min_I \val \det_{I\times J} (g')\,.
\end{equation}
The minimum is taken over all \(i\) element subsets of \([m]\) and \(\det_{I\times J}\) denotes the \(i\times i\) minor using rows \(I\) and columns \(J\). 
%
%
\begin{lemma}[\cite{kamnitzer2010mirkovic} Lemma 2.4]
\label{lem:kamlem}
\(D_\gamma\)
is constructible and
\begin{equation}
    \label{eq:unflipS}
    S^\mu_w = \{gG(\cO) \in \Gr \big| D_{w\omega_i}(gG(\cO)) = \langle \mu,ww_0\omega_i\rangle \text{ for all } i  \}\,. 
\end{equation}
\end{lemma}
Piecing together Kamnitzer's results, we get the following classification theorem. 
\begin{theorem}
  \label{thm:mvcmvp}
  \cite[Proposition 2.2, Theorem 3.1]{kamnitzer2010mirkovic}
Let \(M_\bullet\) be a BZ datum of weight \((\lambda,\mu)\), and set 
\begin{equation}
  \label{eq:muw}
  \mu_w = \sum_{i = 1}^{m-1} M_{w\omega_i} w\alpha_i\,.
\end{equation}
Then
\begin{equation}
  \label{eq:mvc}
  Z = \overline{\{
    L\in\Gr\big| D_{w\omega_i}(L) = M_{ww_0\omega_i} \text{ for all }w\omega_i\in\Gamma
  \}}
\end{equation}
is an MV cycle of type \(\lambda\) and weight \(\mu\).
%
Moreover, 
\begin{equation}
  \label{eq:mvp}
  \Phi(Z) = 
\bigcap_{w\in W} \Phi(\overline{S^{\mu_w}_w}) = P(M_\bullet)\,,
\end{equation}
with 
\begin{equation}
  \label{eq:Mw}
M_{ww_0\omega_i} = \langle\mu_w,ww_0\omega_i\rangle\,.
\end{equation}
All MV cycles (MV polytopes) arise via Equation~\ref{eq:mvc} (Equation~\ref{eq:mvp}). 
\end{theorem}
%

Recall that a choice of reduced word allows us to draw a path in the 1-skeleton of an MV polytope, 
enumerate a subset of the vertices \((\mu_w)_W\) as 
\(\mu_k = \mu_{w_k}\), always starting with \(\mu_0 = \mu_e\) and ending with \(\mu_\ell = \mu_{w_0}\). 
A subset of the fixed points of an MV cycle of type \(\lambda\) and weight \(\mu\) thus acquires an order with the first fixed point given by \(\mu_e = \lambda \) and the last given by \(\mu_{w_0} = \mu\). In the associated polytope, the order on the vertices is reversed, with the lowest vertex \(\mu^{\Pol}_0 = \mu\) and the highest vertex \(\mu^{\Pol}_\ell = \lambda\). This is due to the presence of the \(w_0\) in the definition of our \(D_\gamma\) and the description of our \(S^\mu_w\). In \cite{kamnitzer2010mirkovic}, Kamnitzer works in the left quotient \(G(\cO)\backslash G(\cK)\), and the \(w_0\) is needed to transport the classification established in the left quotient to the right quotient. 

We now zoom in on the behaviour of $D_{\omega_i}$ the constructible functions corresponding to the lowest weight \(\lambda\). 
%
For any \(m\ge b\ge 1\), for any \(g\in G(\cK)\), we denote by \(g_b\) the restrction \(g\big|_{F_b}\) where \(F_b = \Sp_{\cO}(e_1,\dots,e_b)\). When \(g\) is upper triangular, we can (and do) identify \(g_b\) with the \(b\times b\) upper left submatrix of \(g\), an element of \(G_b(\cK) = GL_b(\cK)\). 
Thus, if \(g\in U_-(\cK)\), then for any \(1\le i\le b\le m\),
\begin{equation}
    \label{eq:nonameii}
    \begin{split}
        D_{\omega_i}(g_b G_b(\cO)) &= \val (g_b' v_{w_0^b\omega_i})\\
        &= \val(g' v_{[b-i+1\,b]}) \\
        &= \min_I \val \det_{I\times [b-i+1\, b] }(g') \\ 
        &= \min_K \val \det_{[b-i+1\, b]\times K }(g)\,. 
    \end{split}
\end{equation}
Again, the minimums are taken over all \(i\) element subsets of \([m] = \{1,2,\dots,m\}\), \(w_0^b\) denotes the longest element in the Weyl group of \(G_b = GL_b\), and \([b-i+1\, b]\) denotes the \(m\)-tuple that looks like a string of \(i\) 1s in positions \(b-i+1 \) through \(b\) sandwiched by 0s. 
We denote the function \(\min_K\val\det_{[b-i+1\,b]\times K}(g)\) by \(\Delta_{[b-i+1\,b]}(gG(\cO))\). 

We can define the Lusztig datum of an MV cycle \(Z\) in terms of BZ data as in \cite[Equation~17]{kamnitzer2010mirkovic}. 
\begin{equation}
    \label{eq:ldglbz}
    n_k^{\uvi}(Z) \equiv n_{(a,b)} (Z) = M_{[a\,b]}- M_{[a\,b-1]} - M_{[a+1\,b]} + M_{[a+1\,b-1]} 
\end{equation}
whenever \(\beta_k = \alpha_a + \cdots + \alpha_b\). 
In terms of the constructible functions 
\(\Delta_{[a\,b]}\) this is
\begin{equation}
    \label{eq:ldgldgamma}
        n_{(a,b)} (L) = \Delta_{[a\,b]}(L) - \Delta_{[a\, b-1]}(L) -\Delta_{[a+1\,b]}(L) +  \Delta_{[a+1\,b-1]}(L)
\end{equation}
on generic points \(L = g G(\cO)\) of \(Z\). 
\section{Equal Lusztig data}
\label{s:equal}
%
%
%
Note that 
\begin{equation}
    \label{eq:and}
    \overline{\Gr^\lambda}\cap \overline{S^\mu_-} = \overline{\Gr^\lambda\cap S^\mu_-}\cup\text{ lower dimensional components}\,.
\end{equation}
This follows from the decompositions of \(\overline{\Gr^\lambda}\) and \(\overline{S^\mu_\pm}\) afforded by \cite{mirkovic2007geometric}. In particular, an MV cycle of type \(\lambda\) and weight \(\mu \) is an irreducible component of \(\overline{\Gr^\lambda}\cap \overline{S^\mu_-}\) of dimension \(\langle\rho,\lambda-\mu\rangle\). We need this fact to ensure that the projective closure of the image of a generalized orbital variety is an MV cycle. Indeed, the restricted Mirkovi\'c--Vybornov isomorphism has \(\overline{\Gr^\lambda}\cap S^\mu_- \) as its codomain. 

Now, let us describe more precisely which MV cycle a particular generalized orbital variety will be mapped into. For the remainder of this chapter we fix a tableau \(\tau\in\cT(\lambda)_\mu\). Our first result is that a dense open subset of the generalized orbital variety labeled by \(\tau\) will be mapped to an MV cycle of type \(\lambda \) and weight \(\mu\). 
\begin{lemma}
\label{lem:whereisg} 
For a dense subset of \(A\in X_\tau\), \(\phi(A) \in S_+^\lambda\cap S^\mu_{-}\).
\end{lemma}
\begin{proof}
Let \(\mathring Z_\tau = \phi(X_\tau)\). 
By Corollary~\ref{cor:resmvy}, \(\mathring Z_\tau\)
is an irreducible component of \(\overline{\Gr^\lambda}\cap S^\mu_-\). 
In view of the decomposition recalled in Equation~\ref{eq:and} this means that
\(\mathring Z_\tau \) 
is dense in an irreducible component of \(\overline{\Gr^\lambda\cap S^\mu_-}\), i.e.\ an MV cycle. 
Indeed, \( \dim X_\tau = \langle\rho,\lambda-\mu\rangle\) is also the pure dimension of \(\overline{\Gr^\lambda\cap S^\mu_-}\) by \cite[Theorem 3.2]{mirkovic2007geometric}. 
Taking the projective closure \(\overline{\mathring Z_\tau}\) in \(\overline{\Gr^\lambda\cap S^\mu_-}\) we obtain the corresponding MV cycle, which we denote \(Z_\tau\).

Let \(U_\tau = \phi^{-1}(S^\lambda_+)\cap X_\tau \). Since \(Z_\tau \cap S^\lambda_+ \) is dense and constructible in \(Z_\tau\), \(U_\tau\) is dense and constructible in \(X_\tau\). Clearly \(\phi(A) \in S^\lambda_+\cap S^\mu_-\) for all \(A\in U_\tau\). 
%
\end{proof}
%
It remains to show that \(Z_\tau = \overline{\phi(X_\tau)}\) has Lusztig datum \(n(\tau)\). We begin by decorating the maps involved in the construction of the Mirkovi\'c--Vybornov isomorphism by a number \(1\le b\le m \) to distinguish the underlying group \(GL_b\). 
%
\begin{gather*}
    \tilde\Phi_b : \TT_{\mu^{(b)}}\cap\cN_b\to (G_b)_1 [t^{-1}]t^{\mu^{(b)}} \\
    \Phi_b : \overline{\OO}_{\lambda^{(b)}}\cap\TT_{\mu^{(b)}}\to\overline{\Gr^{\lambda^{(b)}}}\cap \Gr_{\mu^{(b)}} \\
    \phi_b: \overline{\OO}_{\lambda^{(b)}}\cap\TT_{\mu^{(b)}}\cap\n_b\to \overline{\Gr^{\lambda^{(b)}}}\cap S^{\mu^{(b)}}_- 
\end{gather*}    
Here \(\n_b \) and \(\cN_b\) denote upper triangular and nilpotent matrices of size \(\lvert \mu^{(b)}\rvert\times \lvert \mu^{(b)}\rvert\).
Furthermore, we denote by \(\pi_b\) the birational map
\[
    \ot\cap\n\to \overline{\OO}_{\lambda^{(b)}}\cap\TT_{\mu^{(b)}}\cap\n_b \qquad A\mapsto A_{(b)}\,.
\] 
Let \(U^{(b)} = U_{\tau^{(b)}} = \phi_b^{-1}(S^{\lambda^{(b)}}_+)\cap X_{\tau^{(b)}}\) as in the proof of the previous Lemma~\ref{lem:whereisg}. \(U^{(b)}\) is dense and constructible in \(X_{\tau^{(b)}}\). 
Since \(\pi_b : X_\tau \to X_{\tau^{(b)}}\) is surjective, \(U = \pi_b^{-1}(U^{(b)})\) is dense and constructible in \(X_\tau\). 
%

%
%
\begin{lemma}
    \label{lem:dabgb}
The generic value of \(\Delta_{[b-i+1\,b]}\) on \(Z_\tau\)
is \(\langle\lambda^{(b)}, w_0^b\omega_{i}\rangle\). 
\end{lemma}
\begin{proof}
    Given \(A \in U\), let \(g=\tilde\Phi(A)\) and \(g_b = \tilde\Phi_b(A_{(b)})\). 
    Since \(A_{(b)}\in U_{\tau^{(b)}}\), \(g_bG_b(\cO) \in S^{\lambda^{(b)}}_+ \) and the result follows by 
    Lemma~\ref{lem:kamlem} and the calculation in Equation~\ref{eq:nonameii}.
\end{proof}
%
We are now ready to verify our main first main result.
\begin{theorem}
\label{thm:ad2} 
\(\overline{\phi(X_\tau)}\)
is an MV cycle of type \(\lambda\) and weight \(\mu\) having Lusztig datum \(n(\tau)\). 
\end{theorem}
\begin{proof}
%
%
Given \(A\in U\), let \(g = \tilde\Phi(A)\). Then, combining Equation~\ref{eq:ldgldgamma} and Lemma~\ref{lem:dabgb}, we find that for any \(1\le a \le b\le m\),
\[
\begin{split}
    n(gG(\cO))_{(a,b)} 
        &= \langle\lambda^{(b)}, w_0^b\omega_{b-a+1}\rangle - \langle\lambda^{(b)},w_0^b\omega_{b-a}\rangle \\
          &\quad - \langle\lambda^{(b-1)},w_0^{b-1}\omega_{b-a}\rangle + \langle\lambda^{(b-1)},w_0^{b-1}\omega_{b-a-1}\rangle \\
          &= \langle\lambda^{(b)},w_0^b e_{b-a+1}\rangle - \langle\lambda^{(b-1)},w_0^{b-1}e_{b-a}\rangle \\
          &= \lambda^{(b)}_a - \lambda^{(b-1)}_a\,.
\end{split}
\]
%
%
Since \(n(\tau)_{(a,b)} =\lambda^{(b)}_a - \lambda^{(b-1)}_a\) by definition, this checks that the generic value of \(n\) on \(Z_\tau\) is \(n(\tau)\) as desired.
%
\end{proof}
%
%
Let \(Z\) be a stable MV cycle of weight \(\nu\). If \(\lambda\) is the \textit{smallest} dominant effective weight such that \(\mu := -\nu + \lambda\) is dominant effective with \(\mu_m \ne 0\),
then \(t^\mu Z\) is an MV cycle of type \(\lambda\) and weight \( \mu\). 
Let \(\tau\) be the Young tableau of shape \(\lambda\) and weight \(\mu\) whose Lusztig datum is equal to that of \(Z\). Then \(\phi(X_\tau)\) is open and dense in \(t^\mu Z\). 
%
%
\begin{proposition}\label{prop:barDeqm}
    \cite[Proposition A.5]{baumann2019mirkovic}
Let \(Z\in \cZ(\infty)_\nu\) be a stable MV cycle of Lusztig datum \(n\), and let \(X_\tau\) be as above. Then 
    \begin{equation}
        \label{eq:barDeqm}
        \barD(b_{Z}) = \varepsilon^T_{L_0}(Z) = \varepsilon^T_{L_\mu} (t^\mu Z) = \frac{\mdeg_{\TT_\mu\cap\n}(X_\tau)}{\mdeg_{\TT_\mu\cap\n}(\{J_\mu\})} 
    \end{equation}
    where the denominator is given by \(\prod_{\Delta_+}\beta\) by definition of \(\TT_\mu\).
\end{proposition} 
%
%
\section{Pl\"ucker embedding in general}
\label{s:pluckembedgral}
In this section we detail how equations for a generalized orbital variety \(X_\tau\) lead to equations for the corresponding MV cycle \(Z_\tau\). 
We begin by recalling the embedding \(\Upsilon\) of the entire affine Grassmannian based on the Kac--Moody construction of Chapter~\ref{ch:gr} and proceed to give an explicit description of \(\Upsilon\) on the subvarieties \(\overline{\Gr^\lambda}\) using the generalized orbital varieties of Chapter~\ref{ch:govs}.
%
\begin{proposition}
    \label{prop:grembed}
    The map \(\Upsilon : \Gr \to \PP(L(\Lambda_0))\) defined by \(gG^\vee(\cO)\mapsto [g \cdot v_{\Lambda_0}]\) is a \(G^\vee(\cK)\)-equivariant embedding.
\end{proposition}
\begin{proof}[Sketch of proof.]
  %
  To see that \(\Upsilon\) is well-defined we need to show that $G^\vee(\cO)$ fixes $v_{\Lambda_0}$. This is a consequence of the fact that \(\n_-\) and the Iwahori $\fb^\vee(\cO) + t\n^\vee_-(\cO)$ which together generate $\g^\vee(\cO)$ both annihilate \(v_{\Lambda_0}\).
\end{proof}
First let's recall a basic construction. Let \(\Gr_kV \) be the Grassmannian of \(k\)-planes in a finite dimensional complex vector space \(V\). The dual determinant bundle \(\Det^\ast \to\Gr_kV\) is the line bundle whose fibre over \(W\in\Gr_kV\) is equal to \((\bigwedge^kW)^* \cong \bigwedge^k V/W\).
Note that \(\Det\), unlike the tautological bundle, has nonzero holomorphic global sections. In fact, it is ample. 

To describe the orbits \(\Gr^\lambda\) in terms of lattices there is no information lost in restricting attention to those lattices which are contained in the standard lattice \(L_0\). This is because when \(\lambda\) is effective dominant, i.e.\ \(\lambda = (\lambda_1\ge\lambda_2\ge\cdots\ge\lambda_m\ge 0)\), \(gG(\cO)\in\overline{\Gr^\lambda}\) defines \(gL_0\subset L_0\). 
Moreover,
\[
    \{
        L\in\Gr \big| L\subset L_0\text{ and } \begin{bmatrix}t\big|_{L_0/L}\end{bmatrix} \text{ has Jordan type } \lambda
    \} = \Gr^\lambda \,.
\]

Let \(L \in \overline{\Gr^\lambda}\) and take \(p = \lambda_1\). Then \(t^p L_0\subset L\).
Denote by $S$ the set
\[
    ((1,0),(1,1),\dots,(1,p-1),\dots,(m,0),(m,1),\dots,(m,p-1))\,.
\]
We use it to index the standard basis
\(\{v_s\}_{s\in S}\) for \(V = L_0/t^p L_0\) so that \(v_{(i,j)} = [e_it^j]\).
%
%
\begin{lemma}
    \label{lem:todo}
Let \(L\in\overline{\Gr^\lambda}\cap\Gr_\mu\) and choose a representative \(g\in G(\cK)\) such that \(L = gL_0\). Then, as a subspace of \(V\), \( L/t^pL_0 \) has basis
\begin{equation}
    \label{eq:tpL0}
    ( [ge_1], [tge_1],\dots, [t^{p - \mu_1 - 1}g e_1], \dots, [ge_m], [tge_m],\dots,[t^{p - \mu_m - 1} ge_m])\,.
\end{equation}
\end{lemma}
Certainly, \(t^kge_i\in L = gL_0\) for all \(1\le i\le m\) and \(k\ge 0\). The purported spanning set also has the right cardinality, \( mp - N = \dim L/t^pL_0 \).
Shortly, we will see that it is in fact linearly independent.
\begin{proposition}\label{prop:upsilon}
\cite[Proposition A.6]{baumann2019mirkovic}, \cite[Section~2.6]{zhu2016introduction}
Let \(k = mp - N\). 
The map 
\begin{equation}\label{eq:upsilon}
    \upsilon : \overline{\Gr^\lambda} \to \Gr_k V \qquad L\mapsto L/t^p L_0
\end{equation}
is a closed embedding. Moreover, the dual determinant bundle on \(\Gr_k V\) pulls back to the line bundle \(\scL = \Upsilon^\ast \cO(1)\). 
\end{proposition}
\begin{proof}[Proof sketch.]
Let \(L \in \overline{\Gr^\lambda}\). Since \(t^p L_0 \subset L \subset L_0 \), \(L\) is determined by \(L/t^pL_0\). 
To see that \(\upsilon^\ast\Det^\ast = \Upsilon^\ast \cO(1)\), we start by checking that 
\[
    \scL_L \cong \bigwedge^k (L_0/t^pL_0)/(L/t^pL_0) \cong \bigwedge^k L_0/L\,. 
\]

Let \(L = g L_0\) and write \(g = u t^\nu h\) for some \(u \in U(\cK)\), \(\nu\le\lambda \) and \(h\in G(\cO)\). Note \(hL_0 = L_0\).
The natural equivariant structure on \(\Det^\ast\) as well as 
\(G(\cK)\)-equivariance of \(\Upsilon\) implies there is a natural equivariant structure on \(\scL\). In particular, multiplication by \(u\) gives an isomorphism of \(\scL_{L_\nu}\) and \( \scL_L\).

Now \(L_0/L_\nu \) has basis \(([e_it^k]\) for \(1\le i\le m\) and \(1\le k\le\nu_{i-1})\). Wedging these together in degree \(\lvert\nu\rvert\) gives the eigenvector
\[
    v_\nu = [e_1]\wedge [e_1t]\wedge \cdots \wedge [e_1t^{\nu_1-1}]\wedge\cdots\wedge [e_m]\wedge [e_mt]\wedge\cdots\wedge[e_mt^{\nu_m - 1}] 
    \in \bigwedge^{\lvert\nu\rvert} L_0/L_\nu 
\]
which has eigenvalue 
\(
    t_1^{\nu_1}t_2^{\nu_2}\cdots t_m^{\nu_m}    
\)
and, hence, weight \(\nu\). 

On the other hand, \(t^\nu v_{\Lambda_0}\) is contained in the \( \Lambda_0 - \nu + (\nu, \nu) \delta^\vee\) weight space of 
\(L(\Lambda_0)\). 
Since minus the \(T\) weight of the line through \(t^\nu v_{\Lambda_0}\) is \(-(-\nu) = \nu\), we also have that \(\Phi_T(L_\nu) = \nu\), and

Equivariant allows us to conclude that \(\Det^\ast_L\cong\scL_L \) for all \(L\) as desired. 
\end{proof}
We now return to Lemma~\ref{lem:todo}. If we expand the elements in Equation~\ref{eq:tpL0} in the standard basis of \(V\) we obtain a block matrix 
\(B\) of size \((mp - N)\times mp\) 
whose \(i\)th row block has \(r\)th row equal to the vector of coefficients on \([t^{r-1}ge_i]\) expanded in powers of \(t\).
%
The determinant map 
\[
    \bigwedge^{mp - N} : \Gr_{mp - N}(V)\to \bigwedge^{mp - N}V
\]
associates to \(L/t^pL_0\) a line in \(\bigwedge^{mp - N}V\).
The latter space has basis indexed by \((mp-N)\)-subsets \(C\) of \(S\). 
Dually,
the coordinate of the line \( \bigwedge^{mp - N} (L/t^pL_0)\) in \(\PP(\bigwedge^{mp - N}V)\) under the usual Pl\"ucker embedding is given by minors \(\Delta_C(B)\) of \(B\)
using the columns \(C\in\binom{S}{mp - N}\). 

\begin{proof}[Proof of \ref{lem:todo}]
In the event that \(gL_0 = \phi(A)\) for some \(A \in \ot\cap\n\), we can say precisely that 
\[
    ge_i = t^{\mu_i} e_i + \sum_{j > i} (-\sum_k A_{ji}^k t^{k-1})e_j
\]
where, as in the statement of Theorem~\ref{thm:mvy}, \(A_{ji}^k\) denotes the \(k\)th entry from the left of the \(\mu_i\times\mu_j\) block of \(A\). In particular, the highest power occuring in a coefficient of \(e_j \) is at most \(\mu_j -1 \) which is less than \(\mu_i\) when \(j > i\) if we assume that \(\mu\) is dominant. It follows that \(B\) is block ``upper triangular''.
Moreover, the leftmost 1 in the \(r\)th row of the \(i\)th block of \(p-\mu_i\) rows occurs in column \(\mu_i + r\) of the \(i\)th block of \(p\) columns.
This works out to a nonzero pivot entry (equal to 1) in every row, so that the set in Equation~\ref{eq:tpL0} is a basis.
\end{proof}
\begin{proposition}\label{prop:upschain}
\cite[Proposition A.7]{baumann2019mirkovic}
Let \(k = mp -N\).
Under the chain of maps 
\[
     \ot\cap\n \xrightarrow{\phi} \overline{\Gr^\lambda}\cap S^\mu_-\xrightarrow{\upsilon} \Gr_kV \xrightarrow{\psi} \PP(\bigwedge^kV)
\]
\(A\) is sent to the point \(\psi\circ\upsilon\circ\phi (A) = [\Delta_C(B)]_{C\in\binom{S}{k}}\) where \(B\) is the matrix of the basis of \(\phi(A)/t^p L_0\) supplied by Lemma~\ref{lem:todo}
expanded in the standard basis of \(V\). 
\end{proposition}
Let \(C_0\) be the subset of columns containing leading 1s. Then 
\(\Delta_{C_0}(B) = 1\) 
and we can consider the affine space \(\AA^{\binom{mp}{mp-N}-1} \subset \PP(\bigwedge^{mp - N}\CC^{mp})\) defined by the condition \(\Delta_{C_0}\ne 0\). 
Since \(\ot\cap\n\) is mapped into this affine space, the image of a generalized orbital variety \(\mathring Z\) is the intersection of the corresponding MV cycle \(Z\) with this open affine space. 

\begin{definition}\label{def:projcl}
    \cite[Chapter 8]{cox2013ideals}
    The homogenization of an ideal \(K\) is denoted \(K^h \) and defined to be the ideal \((f^h\big|f\in K)\) where
    \begin{equation}
        \label{eq:cox}
        f^h(x_0,x_1,\dots,x_d) = \sum_{i=0}^{\deg f} f_i(x_1,\dots,x_d) x_0^{\deg f - i}
    \end{equation} 
    if \(f = \sum_0^{\deg f} f_i\) is the expansion of \(f\) as the sum of its homogeneous components \(f_i\) of \(\deg f_i = i\).
\end{definition}
\begin{corollary}\label{cor:zaskh}
\cite[Corollary A.8]{baumann2019mirkovic}
The ideal of \(Z\) in \(\PP(\bigwedge^k V)\) equals the homogenization of the kernel of the map 
\begin{equation}\label{eq:zaskh}
    \CC[\{\Delta_C\big| C\in\binom{S}{k}\}]\to\CC[\mathring Z]\,. 
\end{equation}
\end{corollary}
\begin{example}
    \label{eg:ups}
Let \(\lambda = (3,2)\) and \(\mu = (2,2,1)\), and take \(p = 3\). Note \(mp - N = 3\cdot 3 - 5 = 4\).
Recall that \(A\in \ot\cap\n\) takes the form 
\[
    \left[
        \begin{BMAT}(e)[3pt]{cc.cc.c}{cc.cc.c}
            0 & 1 & 0 & 0 & 0 \\
            & 0 & a & b & c \\
            &   & 0 & 1 & 0 \\
            &   &   & 0 & d \\
            &   &   &   & 0 
        \end{BMAT}
    \right]
\]
so \(\phi(A) \) takes the form
\[
\left[
    \begin{BMAT}{c.c.c}{c.c.c}
        t^2 & & \\
        -a - bt & t^2 & \\ 
        -c & -d & t 
    \end{BMAT}
\right]
\]
and
\[
    L_A = \phi(A)L_0 = \Sp_{\cO}(t^2 e_1 - (a+bt) e_2 - ce_3, t^2 e_2 - de_3, te_3)\,. 
\]
To go back, we take the transpose of \(\begin{bmatrix}t\big|_{L_0/L_A}\end{bmatrix}\) in the basis \(([e_1],[te_1], [e_2],[te_2],[e_3])\).

Now 
\[
    L_A/t^3L_0 = \Sp_\CC([t^2 e_1] - a[e_2] - b[te_2] - c[e_3],[t^2e_2] - d[e_3], [te_3], [t^2e_3])    
\]
so
in the basis \(([e_1],[te_1],[t^2e_1],[e_2],[te_2],[t^2e_2],[e_3],[te_3],[t^2e_3])\)
\[
    B =    
    \left[
        \begin{BMAT}{ccc.ccc.ccc}{c.c.cc}
            0 & 0 & 1 &-a &-b & 0 &-c & 0 & 0 \\
            0 & 0 & 0 & 0 & 0 & 1 &-d & 0 & 0 \\
            0 & 0 & 0 & 0 & 0 & 0 & 0 & 1 & 0 \\
            0 & 0 & 0 & 0 & 0 & 0 & 0 & 0 & 1 
        \end{BMAT}
    \right]\,.
\]
Relabeling 
\[
\begin{aligned}
    S &= ((1,0),(1,1),(1,2),(2,0),(2,1),(2,2),(3,0),(3,1),(3,2)) \\
    &= (1,2,3,\bar 1, \bar 2,\bar 3, \b 1, \b 2,\b 3)
\end{aligned}    
\]
we find that \(\Delta_{3\bar 3\scriptsize\b 2\b 3} = 1\).
\end{example}
\subsection{The case \(\mu = \omega_m\)}
\label{ss:pluckembspec}
Let \(N = m\), \(\mu = (1,1,\dots,1)\) and \(p = 2\), so that \(mp - N = m\). 
In this case the Mirkovi\'c--Vybornov isomorphism assumes the very easy form 
\(\Phi(A) = (t^\mu - A')G(\cO)\) 
and the corresponding lattice is 
\begin{equation}
    \label{eq:LA}
    L_A = \Sp_\cO(\{e_it - A' e_i\}_{1\le i\le m})\,. 
\end{equation}
We re-index the basis elements \([e_it^j]\) of \(V = L_0/t^2 L_0 \) as follows.
\begin{equation}
    \label{eq:plbasis}
    (v_1,\dots,v_m, v_{\bar 1},\dots,v_{\bar m}) = ([e_1t],\dots,[e_mt],\dots,[e_1],\dots,[e_m])\,.
\end{equation}
In this basis, \(\upsilon(L_A) \) is the row space of the \(m\times 2m \) matrix \(B = \begin{bmatrix} I & -A \end{bmatrix}\). 
    
As above, subsets \(C \subset S\) of size \(m\) index a basis of \(\bigwedge^m V\), with \(S = (1,\dots,m,\bar 1,\dots,\bar m)\) denoting the special index set. For any \(C\subset S\) of size \(m\), we can consider the minor \(\Delta_C(B)\) using the columns \(C\). 

The results of the last section combine to form a chain of \(T\)-equivariant maps 
\begin{equation}
    \label{eq:onupschain}
    \overline{\OO}_\lambda\cap\n\to\overline{\Gr^\lambda}\cap S^0_-\to\Gr(m,V) \to \PP(\bigwedge^m V)
\end{equation}
under which \(A \) is sent to the line \([\Delta_C(B)]_{C\in\binom{S}{m}}\). 
In particular, the weight of \(v_i\) is equal to the weight of \(v_{\bar i}\) is equal to the standard basis vector \(\varepsilon_i\) of \(\t^\ast\) such that \(\alpha_i = \varepsilon_i - \varepsilon_{i+1}\). 

From the simple form of \(B\), it is immediate that \(\overline{\OO}_\lambda\cap\n\) is mapped into the affine space \(\AA^{\binom{2m}{m} - 1}\subset  \PP(\bigwedge^m V)\) cut out by \(\Delta_{12345}\ne 0\). 

\chapter{Case studies}
\label{ch:calculs}
%
%
Fix \(G = PGL_m\) and \(\uvi = (1\,2\,\dots\,m-1\,\dots\,1\,2\,1)\).
%
In this chapter we present examples of equality and non-equality between ``corresponding'' MV basis and dual semicanonical basis vectors, as well as evidence of ``extra-compatibility'' (both terms recalled below).

We espouse the \textit{ex crystal} point of view, where the basic input is \(\uvi\)-Lusztig data, even though we always start with a tableau.\footnote{Under the rug we are applying the appropriate partial inverse \(B(\infty)\to\cT(\lambda)\) supplied by Theorem~\ref{thm:tabmseg} and described in Chapter~\ref{ch:bases}.} 
The output of primary interest to us is the ideal of a generalized orbital variety. 
Its multidegree and Hilbert series aid and abet us in establishing non-equality and in studying extra-compatibility, respectively.
\section{Varying degrees of equality}
\label{s:degofeq}
Recall the question we set out to answer.
\begin{question*}
    \label{q1}
    If \(M\) and \(Z\) define the same polytope, \(P\), do they define the same measure, \(D\), on \(P\)?
\end{question*}
Now that we have met the structures, the models, and the methods, we are in a position to carry out some computations.
%
Figure~\ref{fig:flow} below depicts the functional flowchart by which we process the data of our examples. 
\begin{figure}[ht!]
  \centering
  $$
\begin{tikzcd}
  &\irr\Lambda(\nu) \ar[dr,"Y\mapsto c_Y"] \\
  \NN^\ell \ar[ur,"n\mapsto Y"]\ar[dr,"n\mapsto Z"'] & & \CC[U] \ar[rr,bend left,"\barD"] \ar[r,"D"] & \cPP \ar[r,"FT"] & \CC[\treg]\\ 
  &\irr\overline{S^\nu_+\cap S^0_-} \ar[ur,"Z\mapsto b_Z"] \ar[rr,"\phi(X) = t^\mu Z\cap\Gr_\mu"',"Z\mapsto X"] & & \irr\ot\cap\n  \ar[ur,"\frac{\mdeg}{p(\mu)}"']
\end{tikzcd}
$$
  \caption{Flowchart}
  \label{fig:flow}
\end{figure}
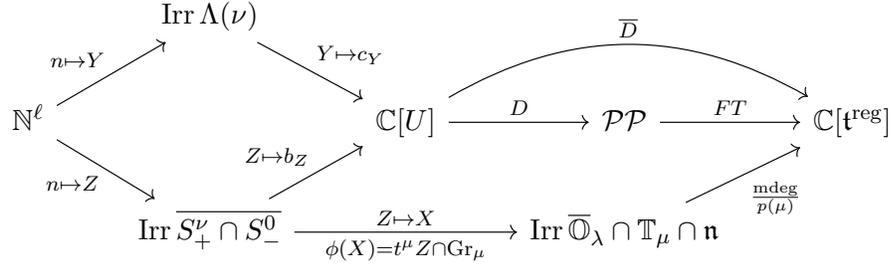

We will say that \(Y\in\irr\Lambda(\nu)\) and \(Z\in\cZ(\infty)_\nu\) \new{correspond}
if \(n(Y) = n(Z)\). 
By way of Theorem~\ref{thm:tabmseg}, we can associate to a Lusztig datum of weight \(\nu\), a tableau of weight \(\mu\), where \(\mu \) is the ``smallest'' possible weight such that \(\nu + \mu \) is a partition. 
\begin{theorem}
  \label{thm:a9}
  (\cite[c.f.\ Corollary A.10]{baumann2019mirkovic})
  If \(\tau \in\cT(\nu + \mu )_\mu\) is such that \(n(\tau) = n(Z) \), then \(t^\mu Z \) is the MV cycle corresponding to the generalized orbital variety labeled by \(\tau\). 
\end{theorem}
%
%
%
This theorem is actually a corollary of our first two main results, Theorems~\ref{mthm:t1v1} and~\ref{mthm:t2v1}:
the generalized orbital variety labeled by a semistandard Young tableau \(\tau\in\cT(\lambda)_\mu\) is mapped under the Mirkovi\'c--Vybornov isomorphism to the MV cycle of weight \((\lambda,\mu)\) and Lusztig datum \(n(\tau)\); therefore, \(t^{-\mu}Z\) is a stable MV cycle of weight \(\lambda - \mu\) and Lusztig datum \(n(\tau)\). 

Using the the fact that equality of polytopes, \(\Pol(Y) = \Pol(Z)\), is the same as equality of \(\uvi\)-Lusztig data, \(n(Z) = n(Y) \), we can actually avoid polytopes in answering the opening question.
Identify \(P^\vee = \ZZ^m/\ZZ(1,1,\dots,1)\) and let \(M \) be a general point in \(Y\). Recall 
that \(D(c_Y)\) (Theorem~\ref{thm:DDh}) and \(D(b_Z)\) (Theorem~\ref{thm:DHFn}) are limits in \(n\) of (scalings of) distributions on \(\t^\ast\) which are obtained from the representation classes 
\begin{equation}
  \label{eq:dsaslim}
  [H^\bullet(F_n(M))], [H^0(Z,\scL^{\otimes n})] \in R(T^\vee) \qquad (n \in \NN)
\end{equation}
via the map that sends a \(T^\vee\)-representation \(V\) to the measure \(\sum\dim V_\nu \delta_\nu\). 
The torus weight spaces are given by \(H^\bullet(F_{n,\nu}(M))\) and \( H^0(Z,\scL^{\otimes n})_{\nu}\) respectively, and
%
corresponding \(Z\) and \(Y\) are said to be \new{extra-compatible} if for all \(n\in\NN \) and for all \(\mu \in Q_+\) 
\begin{equation}
  \label{eq:nmu}
  \chi(F_{n,\nu}(M)) = \dim H^0(Z,\scL^{\otimes n})_{\nu}\,. 
\end{equation}

Note, extra-compatibility is a priori stronger than equality of measures \(D\) or \(\barD\). Conversely, it is not known whether equality of basis vectors implies extra-compatibility. 
All in all, we have the following relations. 
\begin{proposition}\label{prop:allthecomps}
  (cf.\ \cite[Proposition 12.2]{baumann2019mirkovic})
Consider the five statements below.
\begin{enumerate}[label=\roman*)]
    \item \(Y\) and \(Z\) are extra-compatible.
    \item \(c_Y = b_Z\)
    \item \(D(c_Y) = D(b_Z)\)
    \item \(\barD(c_Y) = \barD(b_Z)\)
    \item \(\Pol(Y) = \Pol(Z)\)
\end{enumerate}
They are related as follows.
%
%
\[
\begin{tikzcd}[column sep = huge]
  (i) \ar[blue,dr,Rightarrow,"\text{Def.~\ref{def:dhxav}}+\text{Thm.~\ref{thm:DHFn}}+\text{Thm.~\ref{thm:DDh}}" description]\ar[dd, bend left,Rightarrow,"?" marking] & & (v) \ar[dd,"\text{Thm.~\ref{thm:adcmf}}", Rightarrow, "/" marking]\\
  & (iii) \ar[red,ur,Rightarrow,"\text{Thm.~\ref{thm:DDh}}"']\ar[dr,Rightarrow,"\text{Def. of }\barD"' description] \\
  (ii) \ar[ur,Rightarrow,"\text{Def. of }D"' description]\ar[uu,bend left,Rightarrow,"?" marking] & & (iv)
\end{tikzcd}
\] 
%
%
\end{proposition}

\section{The string order on \(B(\infty)\)}
To check equality of basis vectors (in our first two examples) we rely on the string order on \(B(\infty)\) developed by Baumann in \cite{baumann2011canonical}. 
%
Given \((b',b'')\) and \((c',c'')\) in \(B(\infty)\times B(\infty)\), we write \((b',b'')\approx(c',c'')\) if one of the three conditions below is satisfied.
\begin{itemize}
  \item There is \(i\in I\) such that \(\varphi_i(b') = \varphi_i(b'')\) and \((c',c'') = (\tilde e_i b',\tilde e_i b'')\). 
  \item There is \(i\in I\) such that \(\varphi_i (b') -\varphi_i(b'') > 0 \) and \((c',c'') = (\tilde f_i b',\tilde f_i b'')\). 
  \item \((c',c'') = (\sigma b', \sigma b'')\).
\end{itemize}
Here \(\sigma \) is is an antiautomorphism which we do not define---we will only make use of the first condition.

Given \((b',b'')\in B(\infty)\times B(\infty)\) we write \(b'\le_{\text{str}} b''\) if \(b'\) and \(b''\) have the same weight and if for any finite sequence of elementary moves 
\begin{equation}
  \label{eq:moves}
  (b',b'') = (b_0',b_0'') \approx (b_1',b_1'') \approx \cdots \approx (b_{n}',b_{n}'') 
\end{equation}
one has \(\varphi_i(b_n') \le \varphi_i(b_n'') \) for all \(i\in I\). 

\begin{proposition}
  \cite[Proposition 2.6]{baumann2011canonical}
  \begin{enumerate}
    \item The relation \(\le_{\text{str}}\) is an order on \(B(\infty)\).
    \item The transition matrix between two perfect bases is upper triangular with respect to the order \(\le_{\text{str}}\).
  \end{enumerate}
\end{proposition}
This means that if \(Z\) and \(Y\) correspond, then 
\begin{equation}
  \label{eq:chgofbasis}
  b_Z = \sum_{b(Y')\le_{\text{str}}b(Z)} m(Z,{Y'}) c_{Y'}
\end{equation}
with \(b(Z) = b(Y)\), so \(m(Z,Y) = 1\). Here \(b(?) \) is an element of the \(B(\infty)\) crystal (c.f.\ Section~\ref{s:bmodels}), not to be confused with \(b_?\) (or \(c_?\)) which denotes an element of \(\CC[U]\). 
If we can show that \(b(Z)\) is minimal with respect to the string order \(\le_{\text{str}}\) then we'll have \(b_Z = c_Y\) in \(\CC[U]\).

In particular, to show that some $c_Y$ agrees with some $b_Z$ does not require any geometry, while to show that some pair of corresponding basis vectors is \textit{not} equal---as we will do in our last example---does.

\section{Case 1: Evidence for extra-compatibility in \(A_4\)}
Take
\(\lambda = (2,2,1,0,0) \) and \(\mu = (1,1,1,1,1) = 0\). Put \(\nu = \lambda - \mu =  \alpha_1 + 2\alpha_2 + 2\alpha_3 + \alpha_4\) and
%
%
%
let \(Z\in \cZ(\infty)_\nu\) and \(Y\in\irr\Lambda(\nu)\) be the corresponding pair with Lusztig datum equal to that of the tableau \(\tau\in\cT(\lambda)_\mu\) pictured below.
\[
  \young(12,34,5)
\]
%
This case is interesting because it is the ``smallest'' example not covered by the analysis of \cite[\S 12.5]{baumann2019mirkovic}. 
In particular: 
\begin{itemize}
  \item \(c_Y\) is not a ``flag minor'' (see \cite{geiss2005semicanonical});
  \item general points \(M\in Y\), while still rigid, are not Maya modules (see \cite{kamnitzer2011modules});
  \item \(Z\) is not a Schubert variety (see \cite{baumann2019mirkovic}).
\end{itemize}
%


\begin{proposition}
  \(b_{Z} = c_Y\).
\end{proposition}
\begin{proof}
We will show that \(b(\tau)\) is minimal with respect to the string order. The result then follows by applying Equation~\ref{eq:chgofbasis}.
First note that \( b(\tau) (= b(Z) = b(Y))\) has string datum 
\((3,2,1,4,2,3)\),
meaning that in terms of the highest weight vector \(b_\infty\in B(\infty)\) 
\[
  b(\tau) = f_3 f_2 f_1 f_4 f_2 f_3 b_\infty \,.  
\]
Note that under the map \(\cT(\lambda) \to B(\infty)\), \(b_\infty = b(\tau_h)\) where \(\tau_h\) is the highest weight tableau pictured below.
\[
  \young(11,22,3)
\] 
%
Since elements of different weights are not comparable in the string order, we need only consider the five elements of \(\cT(\lambda)_\mu\) and their images in \(B(\infty)\). These elements along with their \(\varphi_i\) values are listed below. Any other elements of weight \(\lambda-\mu\) in \(B(\infty)\) will necessarily not be less than any element in the image of \(\cT(\lambda)_\mu\) in the string order.
%
\begin{table}[ht]
  \centering
  \newcolumntype{Q}{>{$}l<{$}>{$}l<{$}}
  \begin{tabular}{Q} 
    b & \vec\varepsilon = (\varepsilon_1,\varepsilon_2,\varepsilon_3,\varepsilon_4) \\
    \midrule 
    b_1 & (0,0,1,0) \\
    b_2 & (0,1,0,1) \\
    b_3 & (1,0,0,1) \\ 
    b_4 & (0,1,0,0) \\
    b_5 = b(\tau) & (1,0,1,0) 
  \end{tabular}
\end{table}

\noindent From this table it follows by definition that  
\(b_i\not\le_{\text{str}}b_5\) for any \(i\ne 5\) except perhaps \(i = 1\). 
To see that 
\(b_1 \not\le_{\text{str}} b_5 \)
we need to work a little harder. 
One can check by hand or with the help of Sage that  
\[
\begin{aligned}
  \varepsilon_3(b_1) &= \varepsilon_3(b_5) = 1 \\
  \varepsilon_2(f_3b_1) &= \varepsilon_2(f_3 b_5) = 0 \\
  \varepsilon_1(f_2f_3b_1) &= \varepsilon_1(f_2f_3b_5) = 0 \\
  \varepsilon_4(f_1f_2f_3b_1) &= \varepsilon_4(f_1f_2f_3b_5) = 0 
\end{aligned}  
\]
and 
\[
  \varphi_3(e_4e_1e_2e_3 b_1) = 1 > \varphi_3(e_4e_1e_2e_3b_5) = 0\,.  
\]
We followed this particular sequence of elementary moves, because it is the start of the path back to \(b_\infty\) from \(b_5\) given by our choice of string datum. In terms of tableaux,
\[
  \tau\xrightarrow{3}\young(12,33,5)\xrightarrow{2}\young(12,23,5)\xrightarrow{1}\young(11,23,5)\xrightarrow{4}\young(11,23,4)\xrightarrow{2}\young(11,22,4)\xrightarrow{3} \tau_h 
\]
while 
\[
  \tau_1 = \young(14,25,3) \xrightarrow{3} \young(13,25,3) \xrightarrow{2}\young(12,25,3) \xrightarrow{1} \young(11,25,3) \xrightarrow{4} \young(11,24,3)\xrightarrow{2} 0 
\]
where \(\tau_1\) is such that \(b(\tau_1) = b_1\).
\end{proof}
As a corollary, it follows that \(\barD(b_{Z}) = \barD(c_Y)\). Once we have the ideal of the orbital variety which is open dense in \(t^\mu Z\) 
we can check this directly. 

For our main result, we show that \(Z\) and \(Y\) are extra-compatible to order 2. 
\begin{theorem}
\cite[Theorem A.11]{baumann2019mirkovic}.
\label{thm:aio1}
  Let \(\scL\) be the line bundle afforded by the embedding \(\Upsilon\). (See Chapter~\ref{ch:mvs}.) Let \(M\) be a general point of \(Y\). Then 
  \begin{enumerate}[label=(\roman*)]
    \item For all \(n\in\NN\), \(\dim H^0 (Z,\scL^{\otimes n}) = \dim H^\bullet(F_n(M))\) \label{thm:aio11}
    \item For \(n = 1,2\) and for all \(\nu\in Q_+\), \(\dim H^0 (Z,\scL^{\otimes n})_\nu = \dim H^\bullet(F_{n,\nu}(M))\) \label{thm:aio12}
  \end{enumerate}
\end{theorem}
To prove this theorem we begin by describing the orbital variety of \(Z\).  
\subsection{Orbital variety of \(\tau\)}
If \(A\in\TT_\mu\cap\n = \n\) has coordinate
\[
  \left[
  \begin{BMAT}(@){ccccc}{ccccc}  
  0 & a_1 & a_2 & a_3 & a_4 \\
  0 & 0   & a_5 & a_6 & a_7 \\
  0 & 0   &  0  & a_8 & a_9 \\
  0 & 0   &  0  &  0  & a_{10} \\
  0 & 0   &  0  &  0  & 0
  \end{BMAT}  
  \right]
\]
then the generalized orbital variety \(X\) which is labeled by \(\tau\) is the vanishing locus of the ideal
\[
  I = (a_5, a_{10}, a_1a_6 + a_2a_8, a_7a_8 - a_6a_9, a_1a_7 + a_2a_9)\,.
\]
One can verify either by hand or using Macaulay2 that this ideal is prime. 

Next, applying the rules for computing multidegrees from Chapter~\ref{s:mdegrules} (or using a computer) we find that 
%
%
\begin{equation}
  \label{eq:mdeg1} 
  \mdeg_{\n}(X) = 
  \alpha_2\alpha_4(\alpha_{13} + \alpha_{24} + \alpha_1\alpha_3)\,.
\end{equation}
\subsection{MV cycle}
Let \(p = \lambda_1 = 2\). Since \(X\) is contained in the subspace \(\{a_5, a_{10} = 0\}\), the description of the MV cycle \(t^\mu Z = \overline{\phi(X)}\) afforded by Corollary~\ref{cor:zaskh} can be simplified. 
By ignoring minors which are forced to vanish among the set of \(\binom{10}{5}\) possible Pl\"ucker coordinates, we can exhibit \(t^\mu Z\) as a subvariety of \(\PP^{16}\) using just the following minors. 
$$
\begin{aligned} 
  b_0 &= \Delta_{12345} & b_1 &= \Delta_{1345\overline{1}} & b_2 &= \Delta_{1245\overline{1}} & b_3 &= \Delta_{1235\overline{1}} & b_4 &= \Delta_{1234\overline{1}} & b_5 &= \Delta_{1235\overline{2}} \\
  b_6 &= \Delta_{1234\overline{2}} & b_7 &= \Delta_{1235\overline{3}} & b_8 &= \Delta_{1234\overline{3}} & b_9 &= \Delta_{123\overline{12}} &
  b_{10} &= \Delta_{124\overline{12}} & b_{11} &= \Delta_{124\overline{13}} \\ 
  b_{12} &= \Delta_{125\overline{12}} & b_{13} &= \Delta_{125\overline{13}} & b_{14} &= \Delta_{123\overline{13}} & b_{15} &= \Delta_{134\overline{13}} & b_{16} &= \Delta_{135\overline{13}} 
\end{aligned}$$
Let \(P = \mathbb{C}[b_0, b_1, \dots, b_{16}]\). Then \(t^\mu Z\) is equal to \(\proj(P/K^h)\) where \(K^h\) is the homogenization
of 
\(\Ker(P \to \CC[\n]/I)\). 
Explicitly \(K = K_1 + K_2\) where  
$$
\begin{aligned}
K_1 &= (b_9 - b_3b_6 + b_4b_5, b_{10} - b_2b_6, b_{11} - b_2b_8, \\ 
&\qquad \qquad b_{12} - b_2b_5,b_{13} - b_2b_7, b_{14} - b_3b_8 + b_4b_7, b_{15} - b_1b_8, b_{16} - b_1b_7) \\
K_2 &= (b_1b_5 + b_2b_7, b_6b_7 - b_5b_8, b_1b_6 + b_2b_8)
\end{aligned}
$$
are
the ideals coming from \(\CC[\n]\) and 
\(I\) respectively. 

If we're careful, \(K^h\) is actually the ideal of \(\psi\circ\upsilon(t^\mu Z) \cong t^\mu Z\), but the distinction is moot thanks to Propositions~\ref{prop:upsilon} and~\ref{prop:upschain} which guarantee that 
\begin{equation}
  H^0(\phi\circ \upsilon(t^\mu Z), \cO(1) ) = H^0(\upsilon(t^\mu Z), \Det^\ast) = H^0(t^\mu Z, \scL) \,.
\end{equation}
since \(\psi^\ast\cO(1) = \Det^{-1}\) and \(\upsilon^\ast\Det^{-1} = \scL\).

Since \(H^0(t^\mu Z, \scL) = H^0 (Z,\scL) \), it suffices to compute the Hilbert series of \(K^h\). We do so with the help of Macaulay2. 
First we find that the Hilbert polynomial of \(K^h\) is 
\[
  P_Z = -2P_3 + 18 P_4 -40 P_5 + 25 P_6 
\]
where $P_k$ denotes the Hilbert polynomial of \(\PP^k\) which is given by \( n\mapsto \binom{n + k}{k}\). 

Next we verify that (i.e.\ that \(Z\) is projectively normal)
\[
  P/K^h\cong \bigoplus H^0(Z, \cO_Z(n))  
\]
by checking that \(\depth K^h = 7\). (See \cite[Corollary A.1.12 and Corollary A.1.13]{eisenbud2005geometry}.)
%
We conclude that 
\[
  \dim H^0 (Z, \scL^{\otimes n}) = P_Z(n) 
\]
for all \(n\).

%
%
\subsection{Preprojective algebra module}
Using Theorem~\ref{thm:modld} we find that the (abridged---0 summands omitted) HN table for the Lusztig datum \(n(\tau)\),
\[
  \begin{BMAT}{c}{c|c|c|c}
    1\\
    2\\
    2\leftarrow 3 \\
    3\leftarrow 4
  \end{BMAT}\,,
\]
determines a general point \(M\in Y\) of the form 
\[
  \begin{tikzpicture}
    \begin{scope}[xshift=6cm]
    \node (02) at (0,2){$1$};  
    \node (22) at (2,2){$3$};
    \node (11) at (1,1){$22$};
    \node (31) at (3,1){$4$};
    \node (20) at (2,0){$3$}; 
    \draw[->] (02) -- (11); 
    \draw[->] (22) -- (11); 
    \draw[->] (22) -- (31); 
    \draw[->] (11) -- (20); 
    \draw[->] (31) -- (20); 
    \end{scope}
  \end{tikzpicture}
\]
with the maps chosen such that \(\Ker(M_2\to M_3) \), \(\Im (M_3\to M_2)\) and \(\Im (M_1 \to M_2)\) are distinct. 
%

This module has the property that any two submodules of a given dimension vector are isomorphic, and this property is exploited in \cite[Appendix]{baumann2019mirkovic} to compute \(\dim H^\bullet(F_n(M))\) recursively as follows. 
Given \(\nu\in Q_+\) we denote by \(F_n(M)^\nu\) the component of \(F_n(M)\) consisting of those \((n+1)\)-step flags that have dimension \(\nu\) at step \(n\). 
Since all submodules of dimension \(\nu\) are isomorphic, there exists \(N\subset M\) such that \(F_n(M)^\nu\cong F_{n-1}(N)\times F_{1,\nu}(M)\). 
The resulting recurrence 
\[
  \dim H^\bullet(F_n(M)^\nu) = \dim H^\bullet(F_{n-1}(N)) + \dim H^\bullet (F_{1,\nu}(M))
\] 
is helpful for computing the right-hand side of the formula
\[
   \dim H^\bullet(F_n(M)) = \sum_{\nu\in Q_+} \dim H^\bullet(F_n(M)^\nu)\,.
\]
By induction on \(n\) one finds that 
\[
  \dim H^\bullet(F_n(M)) = \frac{(n+1)^2(n+2)^2(n + 3)(5n + 12)}{144}  
\]
which is equal to \(P_Z(n)\). 
See Table~\ref{tab:A4exgr1} below. This shows Theorem~\ref{thm:aio1}~\ref{thm:aio11}.
\begin{table}[ht!]
\centering
\begin{tabular}{@{}llll@{}} \toprule
$\nu$ & $F_1(M)^\nu$ & $\dim H^\bullet(F_{1,\nu}(M))$ & $\dim H^\bullet(F_n(M)^\nu)$\\ 
\midrule
$(0, 0, 0, 0)$ & Point & 1 & $1$\\
$(0, 1, 0, 0)$ & Point & 1 & $n$\\
$(0, 0, 1, 0)$ & Point & 1 & $n$\\
$(0, 0, 1, 1)$ & Point & 1 & $\frac{1}{2}n(n+1)$\\
$(0, 1, 1, 0)$ & $\mathbb{P}^1$ & 2 & $\frac{1}{2}n(3n+1)$\\
$(0, 1, 1, 1)$ & $\mathbb{P}^1$ & 2 & $\frac{1}{6}n(n+1)(5n+1)$\\
$(0, 2, 1, 0)$ & Point & 1 & $\frac{1}{2}n^2(n+1)$\\
$(1, 1, 1, 0)$ & Point & 1 & $\frac{1}{6}n(n+1)(n+2)$\\
$(0, 1, 2, 1)$ & Point & 1 & $\frac{1}{12}n(n+1)^2(n+2)$\\
$(0, 2, 1, 1)$ & Point & 1 & $\frac{1}{6}n^2(n+1)(2n+1)$\\
$(1, 1, 1, 1)$ & Point & 1 & $\frac{1}{24}n(n+1)(n+2)(3n+1)$\\
$(1, 2, 1, 0)$ & Point & 1 & $\frac{1}{6}n^2(n+1)(n+2)$\\
$(0, 2, 2, 1)$ & Point & 1 & $\frac{1}{12}n^2(n+1)^2(n+2)$\\
$(1, 2, 1, 1)$ & Point & 1 & $\frac{1}{24}n^2(n+1)(n+2)(3n+1)$\\
$(1, 2, 2, 1)$ & Point & 1 & $\frac{1}{144}n^2(n+1)^2(n+2)(5n+7)$\\
\bottomrule
\end{tabular}

\hfill
\caption{Flags of submodules of $M$.\label{tab:A4exgr1}}
\end{table}
\subsection{Extra-compatibility}
Recall that $F_{n, \nu}(M)$ denotes the variety of $(n + 1)$-step flags of submodules of $M$ whose dimension vectors sum to a fixed dimension vector $\nu$. 
%
%
%
In Theorem~\ref{thm:aio1}~\ref{thm:aio12} we wish to show that 
\begin{equation}
  \label{eq:eca4}
  \dim H^0(t^\mu Z,\cO(n))_\nu = \dim H^\bullet(F_{n,\nu}(M))
\end{equation}
for \(n = 1,2\) and for all \(\mu\in P\). We can describe \(F_{n,\mu}(M)\) by hand. 
\begin{table}[ht!]
  \begin{longtable}{MZM}
      \toprule
  \nu & (P/J^h)_{\nu + \mu} & H^\bullet(F_1(M)^\nu) \\
  \midrule
  0 & b_0 & 1 \\
  \alpha_1 & b_1 = a_1 & 1 \\
  \alpha_1 + \alpha_2 & b_2 = a_2 & 1 \\
  \alpha_1 + \alpha_2 + \alpha_3 & b_3 = a_3, b_{13} = a_2a_8 & 2 \\ 
  \alpha_1 + \alpha_2 + \alpha_3 + \alpha_4 & b_4 = a_4,b_{11} = a_2a_9 & 2\\
  \alpha_2 + \alpha_3 & b_5 = a_6 & 1 \\
  \alpha_2 + \alpha_3 + \alpha_4 & b_6 = a_7 & 1 \\
  \alpha_ 3 & b_7 = a_8 & 1 \\
  \alpha_3 + \alpha_4 & b_8 = a_9 & 1 \\
  \alpha_1 + \alpha_3 & b_{16} = a_1a_8 & 1 \\ 
  \alpha_1 + \alpha_3 + \alpha_4 & b_{15} = a_1a_9 & 1 \\
  \alpha_1 + 2\alpha_2 + \alpha_3 & b_{12} = a_2a_6 & 1 \\
  \alpha_1 + 2\alpha_2 + \alpha_3 + \alpha_4 & b_{10} = a_2a_7 & 1 \\
  \alpha_1 + \alpha_2 + 2\alpha_3 + \alpha_4 & b_{14} = a_3a_9 - a_4a_8 & 1 \\
  \alpha_1 + 2\alpha_2 + 2\alpha_3 + \alpha_4 & b_9 = a_3a_7 - a_4a_6 & 1 \\
  \midrule
  & & 17 \\
  \bottomrule \\
  \caption{
    $F_1(M)$ and $Z$ as torus representations
      \label{a4onestep}
  }
  \end{longtable}
\end{table}

\newpage To describe the left-hand side we use Macaulay2, computing a truncated Hilbert series of \(J^h\) in the polynomial ring \(P\) which is multigraded by the usual degree plus the torus action. The latter grading is inherited from \(\CC[\n]\). 
Indeed, the information of the first and second columns of the above table is contained in the following truncated Hilbert series, with the formal variable \(u_i = t_i/t_{i+1}\) capturing weight \(\alpha_i\) piece of \(I\subset\CC[\n]\). 
The torus graded Hilbert series of the 17 dimensional representation \(H^0(Z,\cO(1))\) is 

  $
      1+u_1+u_3+u_1u_2+u_1u_3+u_2u_3+u_3u_4+2u_1u_2u_3+u_1u_3u_4+u_2u_3u_4+u_1u_2^2u_3+2u_1u_2u_3u_4+u_1u_2^2u_3u_4+u_1u_2u_3^2u_4+u_1u_2^2u_3^2u_4
  $.

Remembering some more terms we can see the correspondence on the level of torus representations between \(F_2(M)\) and \(Z\). The torus graded Hilbert series of the 110 dimensional representation \(H^0(Z,\cO(2))\) is:

  $
  1+u_1+u_3+u_1^2+u_1u_2+u_1u_3+u_2u_3+u_3^2+u_3u_4+u_1^2u_2+u_1^2u_3+2u_1u_2u_3+u_1u_3^2+u_2u_3^2+u_1u_3u_4+u_2u_3u_4+u_3^2u_4+u_1^2u_2^2+2u_1^2u_2u_3+u_1u_2^2u_3+u_1^2u_3^2+2u_1u_2u_3^2+u_2^2u_3^2+u_1^2u_3u_4+2u_1u_2u_3u_4+u_1u_3^2u_4+u_2u_3^2u_4+u_3^2u_4^2+2u_1^2u_2^2u_3+2u_1^2u_2u_3^2+2u_1u_2^2u_3^2+2u_1^2u_2u_3u_4+u_1u_2^2u_3u_4+u_1^2u_3^2u_4+3u_1u_2u_3^2u_4+u_2^2u_3^2u_4+u_1u_3^2u_4^2+u_2u_3^2u_4^2+u_1^2u_2^3u_3+3u_1^2u_2^2u_3^2+u_1u_2^3u_3^2+2u_1^2u_2^2u_3u_4+3u_1^2u_2u_3^2u_4+3u_1u_2^2u_3^2u_4+u_1u_2u_3^3u_4+u_1^2u_3^2u_4^2+2u_1u_2u_3^2u_4^2+u_2^2u_3^2u_4^2+2u_1^2u_2^3u_3^2+u_1^2u_2^3u_3u_4+4u_1^2u_2^2u_3^2u_4+u_1u_2^3u_3^2u_4+u_1^2u_2u_3^3u_4+u_1u_2^2u_3^3u_4+2u_1^2u_2u_3^2u_4^2+2u_1u_2^2u_3^2u_4^2+u_1u_2u_3^3u_4^2+u_1^2u_2^4u_3^2+3u_1^2u_2^3u_3^2u_4+2u_1^2u_2^2u_3^3u_4+u_1u_2^3u_3^3u_4+3u_1^2u_2^2u_3^2u_4^2+u_1u_2^3u_3^2u_4^2+u_1^2u_2u_3^3u_4^2+u_1u_2^2u_3^3u_4^2+u_1^2u_2^4u_3^2u_4+2u_1^2u_2^3u_3^3u_4+2u_1^2u_2^3u_3^2u_4^2+2u_1^2u_2^2u_3^3u_4^2+u_1u_2^3u_3^3u_4^2+u_1^2u_2^4u_3^3u_4+u_1^2u_2^4u_3^2u_4^2+2u_1^2u_2^3u_3^3u_4^2+u_1^2u_2^2u_3^4u_4^2+u_1^2u_2^4u_3^3u_4^2+u_1^2u_2^3u_3^4u_4^2+u_1^2u_2^4u_3^4u_4^2
  $.
%

In Table~\ref{tab:a4twostep} below we see that \(\dim H^\bullet(F_{2,\nu}(M))\) and the coefficients on the Hilbert series of \(H^0(Z,\cO(2))\), the dimensions of \(H^0(Z,\cO(2))_\nu = H^0(t^\mu Z,\cO(n))_{\nu+2\mu} \), coincide. 
%
\begin{center}
  \begin{longtable}{HHcHHc}
  \toprule
Dimension vector & $\nu$ & $\nu + 2\mu$ & $(P/J^h)_{\nu + 2\mu}$ & Basis size & $H^\bullet(F_{2,\nu}(M))$\\
\midrule
$(0, 0, 0, 0)$ & $(0, 0, 0, 0, 0)$ & $(2, 2, 2, 2, 2)$ & $a_{0}a_{0}$ & $1$ & $1$ \\
$(1, 0, 0, 0)$ & $(1, -1, 0, 0, 0)$ & $(3, 1, 2, 2, 2)$ & $a_{0}a_{1}$ & $1$ & $1$ \\
$(1, 1, 0, 0)$ & $(1, 0, -1, 0, 0)$ & $(3, 2, 1, 2, 2)$ & $a_{0}a_{2}$ & $1$ & $1$ \\
$(1, 1, 1, 0)$ & $(1, 0, 0, -1, 0)$ & $(3, 2, 2, 1, 2)$ & $a_{0}a_{3}, a_{0}a_{28}$ & $2$ & $2$ \\
$(1, 1, 1, 1)$ & $(1, 0, 0, 0, -1)$ & $(3, 2, 2, 2, 1)$ & $a_{0}a_{4}, a_{0}a_{29}$ & $2$ & $2$ \\
$(0, 1, 1, 0)$ & $(0, 1, 0, -1, 0)$ & $(2, 3, 2, 1, 2)$ & $a_{0}a_{6}$ & $1$ & $1$ \\
$(0, 1, 1, 1)$ & $(0, 1, 0, 0, -1)$ & $(2, 3, 2, 2, 1)$ & $a_{0}a_{7}$ & $1$ & $1$ \\
$(0, 0, 1, 0)$ & $(0, 0, 1, -1, 0)$ & $(2, 2, 3, 1, 2)$ & $a_{0}a_{8}$ & $1$ & $1$ \\
$(0, 0, 1, 1)$ & $(0, 0, 1, 0, -1)$ & $(2, 2, 3, 2, 1)$ & $a_{0}a_{9}$ & $1$ & $1$ \\
$(1, 2, 1, 1)$ & $(1, 1, -1, 0, -1)$ & $(3, 3, 1, 2, 1)$ & $a_{0}a_{27}$ & $1$ & $1$ \\
$(1, 2, 1, 0)$ & $(1, 1, -1, -1, 0)$ & $(3, 3, 1, 1, 2)$ & $a_{0}a_{26}$ & $1$ & $1$ \\
$(1, 0, 1, 1)$ & $(1, -1, 1, 0, -1)$ & $(3, 1, 3, 2, 1)$ & $a_{0}a_{19}$ & $1$ & $1$ \\
$(1, 0, 1, 0)$ & $(1, -1, 1, -1, 0)$ & $(3, 1, 3, 1, 2)$ & $a_{0}a_{18}$ & $1$ & $1$ \\
$(2, 0, 0, 0)$ & $(2, -2, 0, 0, 0)$ & $(4, 0, 2, 2, 2)$ & $a_{1}a_{1}$ & $1$ & $1$ \\
$(2, 1, 0, 0)$ & $(2, -1, -1, 0, 0)$ & $(4, 1, 1, 2, 2)$ & $a_{1}a_{2}$ & $1$ & $1$ \\
$(2, 1, 1, 0)$ & $(2, -1, 0, -1, 0)$ & $(4, 1, 2, 1, 2)$ & $a_{1}a_{3}, a_{1}a_{28}$ & $2$ & $2$ \\
$(2, 1, 1, 1)$ & $(2, -1, 0, 0, -1)$ & $(4, 1, 2, 2, 1)$ & $a_{1}a_{4}, a_{1}a_{29}$ & $2$ & $2$ \\
$(2, 2, 1, 1)$ & $(2, 0, -1, 0, -1)$ & $(4, 2, 1, 2, 1)$ & $a_{1}a_{27}, a_{2}a_{4}$ & $2$ & $2$ \\
$(2, 2, 1, 0)$ & $(2, 0, -1, -1, 0)$ & $(4, 2, 1, 1, 2)$ & $a_{1}a_{26}, a_{2}a_{3}$ & $2$ & $2$ \\
$(2, 0, 1, 1)$ & $(2, -2, 1, 0, -1)$ & $(4, 0, 3, 2, 1)$ & $a_{1}a_{19}$ & $1$ & $1$ \\
$(2, 0, 1, 0)$ & $(2, -2, 1, -1, 0)$ & $(4, 0, 3, 1, 2)$ & $a_{1}a_{18}$ & $1$ & $1$ \\
$(2, 2, 0, 0)$ & $(2, 0, -2, 0, 0)$ & $(4, 2, 0, 2, 2)$ & $a_{2}a_{2}$ & $1$ & $1$ \\
$(2, 3, 1, 1)$ & $(2, 1, -2, 0, -1)$ & $(4, 3, 0, 2, 1)$ & $a_{2}a_{27}$ & $1$ & $1$ \\
$(2, 3, 1, 0)$ & $(2, 1, -2, -1, 0)$ & $(4, 3, 0, 1, 2)$ & $a_{2}a_{26}$ & $1$ & $1$ \\
$(2, 2, 2, 0)$ & $(2, 0, 0, -2, 0)$ & $(4, 2, 2, 0, 2)$ & $a_{3}a_{3}, a_{3}a_{28}, a_{28}a_{28}$ & $3$ & $3$ \\
$(2, 2, 2, 1)$ & $(2, 0, 0, -1, -1)$ & $(4, 2, 2, 1, 1)$ & $a_{3}a_{4}, a_{3}a_{29}, a_{4}a_{28}, a_{29}a_{28}$ & $4$ & $4$ \\
$(1, 2, 2, 0)$ & $(1, 1, 0, -2, 0)$ & $(3, 3, 2, 0, 2)$ & $a_{3}a_{6}, a_{6}a_{28}$ & $2$ & $2$ \\
$(1, 2, 2, 1)$ & $(1, 1, 0, -1, -1)$ & $(3, 3, 2, 1, 1)$ & $a_{3}a_{7}, a_{4}a_{6}, a_{6}a_{29}$ & $3$ & $3$ \\
$(1, 1, 2, 0)$ & $(1, 0, 1, -2, 0)$ & $(3, 2, 3, 0, 2)$ & $a_{3}a_{8}, a_{6}a_{18}$ & $2$ & $2$ \\
$(1, 1, 2, 1)$ & $(1, 0, 1, -1, -1)$ & $(3, 2, 3, 1, 1)$ & $a_{3}a_{9}, a_{4}a_{8}, a_{6}a_{19}$ & $3$ & $3$ \\
$(2, 3, 3, 1)$ & $(2, 1, 0, -2, -1)$ & $(4, 3, 2, 0, 1)$ & $a_{3}a_{3746}, a_{3746}a_{28}$ & $2$ & $2$ \\
$(2, 3, 2, 1)$ & $(2, 1, -1, -1, -1)$ & $(4, 3, 1, 1, 1)$ & $a_{3}a_{27}, a_{4}a_{26}, a_{27}a_{28}$ & $3$ & $3$ \\
$(2, 3, 2, 0)$ & $(2, 1, -1, -2, 0)$ & $(4, 3, 1, 0, 2)$ & $a_{3}a_{26}, a_{26}a_{28}$ & $2$ & $2$ \\
$(2, 2, 3, 1)$ & $(2, 0, 1, -2, -1)$ & $(4, 2, 3, 0, 1)$ & $a_{3}a_{3948}, a_{3746}a_{18}$ & $2$ & $2$ \\
$(2, 1, 2, 1)$ & $(2, -1, 1, -1, -1)$ & $(4, 1, 3, 1, 1)$ & $a_{3}a_{19}, a_{4}a_{18}, a_{29}a_{18}$ & $3$ & $3$ \\
$(2, 1, 2, 0)$ & $(2, -1, 1, -2, 0)$ & $(4, 1, 3, 0, 2)$ & $a_{3}a_{18}, a_{28}a_{18}$ & $2$ & $2$ \\
$(2, 2, 2, 2)$ & $(2, 0, 0, 0, -2)$ & $(4, 2, 2, 2, 0)$ & $a_{4}a_{4}, a_{4}a_{29}, a_{29}a_{29}$ & $3$ & $3$ \\
$(1, 2, 2, 2)$ & $(1, 1, 0, 0, -2)$ & $(3, 3, 2, 2, 0)$ & $a_{4}a_{7}, a_{7}a_{29}$ & $2$ & $2$ \\
$(1, 1, 2, 2)$ & $(1, 0, 1, 0, -2)$ & $(3, 2, 3, 2, 0)$ & $a_{4}a_{9}, a_{7}a_{19}$ & $2$ & $2$ \\
$(2, 3, 3, 2)$ & $(2, 1, 0, -1, -2)$ & $(4, 3, 2, 1, 0)$ & $a_{4}a_{3746}, a_{3746}a_{29}$ & $2$ & $2$ \\
$(2, 3, 2, 2)$ & $(2, 1, -1, 0, -2)$ & $(4, 3, 1, 2, 0)$ & $a_{4}a_{27}, a_{27}a_{29}$ & $2$ & $2$ \\
$(2, 2, 3, 2)$ & $(2, 0, 1, -1, -2)$ & $(4, 2, 3, 1, 0)$ & $a_{4}a_{3948}, a_{3746}a_{19}$ & $2$ & $2$ \\
$(2, 1, 2, 2)$ & $(2, -1, 1, 0, -2)$ & $(4, 1, 3, 2, 0)$ & $a_{4}a_{19}, a_{29}a_{19}$ & $2$ & $2$ \\
$(0, 2, 2, 0)$ & $(0, 2, 0, -2, 0)$ & $(2, 4, 2, 0, 2)$ & $a_{6}a_{6}$ & $1$ & $1$ \\
$(0, 2, 2, 1)$ & $(0, 2, 0, -1, -1)$ & $(2, 4, 2, 1, 1)$ & $a_{6}a_{7}$ & $1$ & $1$ \\
$(0, 1, 2, 0)$ & $(0, 1, 1, -2, 0)$ & $(2, 3, 3, 0, 2)$ & $a_{6}a_{8}$ & $1$ & $1$ \\
$(0, 1, 2, 1)$ & $(0, 1, 1, -1, -1)$ & $(2, 3, 3, 1, 1)$ & $a_{6}a_{9}$ & $1$ & $1$ \\
$(1, 3, 3, 1)$ & $(1, 2, 0, -2, -1)$ & $(3, 4, 2, 0, 1)$ & $a_{6}a_{3746}$ & $1$ & $1$ \\
$(1, 3, 2, 1)$ & $(1, 2, -1, -1, -1)$ & $(3, 4, 1, 1, 1)$ & $a_{6}a_{27}$ & $1$ & $1$ \\
$(1, 3, 2, 0)$ & $(1, 2, -1, -2, 0)$ & $(3, 4, 1, 0, 2)$ & $a_{6}a_{26}$ & $1$ & $1$ \\
$(1, 2, 3, 1)$ & $(1, 1, 1, -2, -1)$ & $(3, 3, 3, 0, 1)$ & $a_{6}a_{3948}$ & $1$ & $1$ \\
$(0, 2, 2, 2)$ & $(0, 2, 0, 0, -2)$ & $(2, 4, 2, 2, 0)$ & $a_{7}a_{7}$ & $1$ & $1$ \\
$(0, 1, 2, 2)$ & $(0, 1, 1, 0, -2)$ & $(2, 3, 3, 2, 0)$ & $a_{7}a_{9}$ & $1$ & $1$ \\
$(1, 3, 3, 2)$ & $(1, 2, 0, -1, -2)$ & $(3, 4, 2, 1, 0)$ & $a_{7}a_{3746}$ & $1$ & $1$ \\
$(1, 3, 2, 2)$ & $(1, 2, -1, 0, -2)$ & $(3, 4, 1, 2, 0)$ & $a_{7}a_{27}$ & $1$ & $1$ \\
$(1, 2, 3, 2)$ & $(1, 1, 1, -1, -2)$ & $(3, 3, 3, 1, 0)$ & $a_{7}a_{3948}$ & $1$ & $1$ \\
$(0, 0, 2, 0)$ & $(0, 0, 2, -2, 0)$ & $(2, 2, 4, 0, 2)$ & $a_{8}a_{8}$ & $1$ & $1$ \\
$(0, 0, 2, 1)$ & $(0, 0, 2, -1, -1)$ & $(2, 2, 4, 1, 1)$ & $a_{8}a_{9}$ & $1$ & $1$ \\
$(1, 1, 3, 1)$ & $(1, 0, 2, -2, -1)$ & $(3, 2, 4, 0, 1)$ & $a_{8}a_{3948}$ & $1$ & $1$ \\
$(1, 0, 2, 1)$ & $(1, -1, 2, -1, -1)$ & $(3, 1, 4, 1, 1)$ & $a_{8}a_{19}$ & $1$ & $1$ \\
$(1, 0, 2, 0)$ & $(1, -1, 2, -2, 0)$ & $(3, 1, 4, 0, 2)$ & $a_{8}a_{18}$ & $1$ & $1$ \\
$(0, 0, 2, 2)$ & $(0, 0, 2, 0, -2)$ & $(2, 2, 4, 2, 0)$ & $a_{9}a_{9}$ & $1$ & $1$ \\
$(1, 1, 3, 2)$ & $(1, 0, 2, -1, -2)$ & $(3, 2, 4, 1, 0)$ & $a_{9}a_{3948}$ & $1$ & $1$ \\
$(1, 0, 2, 2)$ & $(1, -1, 2, 0, -2)$ & $(3, 1, 4, 2, 0)$ & $a_{9}a_{19}$ & $1$ & $1$ \\
$(2, 4, 4, 2)$ & $(2, 2, 0, -2, -2)$ & $(4, 4, 2, 0, 0)$ & $a_{3746}a_{3746}$ & $1$ & $1$ \\
$(2, 4, 3, 2)$ & $(2, 2, -1, -1, -2)$ & $(4, 4, 1, 1, 0)$ & $a_{3746}a_{27}$ & $1$ & $1$ \\
$(2, 4, 3, 1)$ & $(2, 2, -1, -2, -1)$ & $(4, 4, 1, 0, 1)$ & $a_{3746}a_{26}$ & $1$ & $1$ \\
$(2, 3, 4, 2)$ & $(2, 1, 1, -2, -2)$ & $(4, 3, 3, 0, 0)$ & $a_{3746}a_{3948}$ & $1$ & $1$ \\
$(2, 4, 2, 2)$ & $(2, 2, -2, 0, -2)$ & $(4, 4, 0, 2, 0)$ & $a_{27}a_{27}$ & $1$ & $1$ \\
$(2, 4, 2, 1)$ & $(2, 2, -2, -1, -1)$ & $(4, 4, 0, 1, 1)$ & $a_{27}a_{26}$ & $1$ & $1$ \\
$(2, 4, 2, 0)$ & $(2, 2, -2, -2, 0)$ & $(4, 4, 0, 0, 2)$ & $a_{26}a_{26}$ & $1$ & $1$ \\
$(2, 2, 4, 2)$ & $(2, 0, 2, -2, -2)$ & $(4, 2, 4, 0, 0)$ & $a_{3948}a_{3948}$ & $1$ & $1$ \\
$(2, 1, 3, 2)$ & $(2, -1, 2, -1, -2)$ & $(4, 1, 4, 1, 0)$ & $a_{3948}a_{19}$ & $1$ & $1$ \\
$(2, 1, 3, 1)$ & $(2, -1, 2, -2, -1)$ & $(4, 1, 4, 0, 1)$ & $a_{3948}a_{18}$ & $1$ & $1$ \\
$(2, 0, 2, 2)$ & $(2, -2, 2, 0, -2)$ & $(4, 0, 4, 2, 0)$ & $a_{19}a_{19}$ & $1$ & $1$ \\
$(2, 0, 2, 1)$ & $(2, -2, 2, -1, -1)$ & $(4, 0, 4, 1, 1)$ & $a_{19}a_{18}$ & $1$ & $1$ \\
$(2, 0, 2, 0)$ & $(2, -2, 2, -2, 0)$ & $(4, 0, 4, 0, 2)$ & $a_{18}a_{18}$ & $1$ & $1$ \\ \midrule
& & & & $110$ & $110$ \\ \bottomrule \\
\caption{$F_2(M)$ and $Z$ as torus representations.\label{tab:a4twostep}}
\end{longtable}
\end{center}
\vspace*{-1.5cm}
This check completes Theorem~\ref{thm:aio1}~\ref{thm:aio12}.
\subsection{Equal \(\barD\)}
We now check directly that \(\barD(b_Z) = \barD(c_Y)\). 
To compute the flag function of \(Y\) we analyze the geometry of composition series \(F_\vi(M)\). 

If \(\nu\) belongs to the subset
\[
\begin{aligned}
  (1, 2, 3, 2, 4, 3) & & (1, 2, 3, 4, 2, 3) & & (1, 3, 2, 4, 3, 2) & & 
  (3, 4, 1, 2, 3, 2) \\
  (1, 3, 4, 2, 3, 2) & & (3, 1, 2, 4, 3, 2) & & (3, 1, 4, 2, 3, 2)  & & (3, 4, 2, 1, 2, 3) \\
  (3, 2, 1, 2, 4, 3) & & (3, 2, 1, 4, 2, 3) & & (3, 2, 4, 1, 2, 3) 
\end{aligned}  
\] 
of $\Seq(\nu)$
then $F_{\vi}(M)$ is a point, so $\chi(F_{\vi}(M)) = 1$. 

If \(\nu\) belongs to the subset 
\[
\begin{aligned}
  (1, 3, 2, 2, 4, 3) & & (1, 3, 2, 4, 2, 3) & & (1, 3, 4, 2, 2, 3) & & (3, 1, 2, 2, 4, 3) \\
  (3, 1, 2, 4, 2, 3) & & (3, 1, 4, 2, 2, 3) & & (3, 4, 1, 2, 2, 3) 
\end{aligned}  
\] 
of $\Seq(\nu)$
then $F_{\vi}(M)\cong\mathbb P ^1$, so $\chi(F_{\vi}(M)) = 2$.
For all other values of $\vi$ in $\Seq(\nu)$, the variety is empty.

The flag function is a rational function, but we can use $p(\mu) = \mdeg_{\n}(0)$ to clear the denominator. By direct computation, one obtains that the flag function is given by $\barD(c_Y) p(\mu)$. Comparing with $\mdeg_{\n}(X)$ and applying Proposition~\ref{prop:barDeqm} proves the expected equality.
\section{Case 2: Weak evidence of extra-compatibility in \(A_5\)}
This time, we take \(\lambda = (2,2,1,1) \) and \(\mu = (1,1,1,1,1,1)\).
Put \(\nu = \lambda - \mu =  \alpha_1 + 2\alpha_2 + 2\alpha_3 + 2\alpha_4 + \alpha_5\) and
%
%
%
let \(Z\in \cZ(\infty)_\nu\) and \(Y\in\irr\Lambda(\nu)\) be the corresponding pair with Lusztig datum equal to that of the tableau \(\tau\in\cT(\lambda)_\mu\) pictured below.
%
%
%
\[
  \tau = \young(13,25,4,6)
\]
%
%
%
This case is interesting because it is the simplest example for which \(Y\) is a component whose general point is \textit{not} rigid.
%
%
\begin{proposition}
  \(b_Z = c_Y\)
\end{proposition}
\begin{proof}
  First note that \(b(\tau) = b(Z) = b(Y)\) has string data \((2,1,3,4,3,2,5,4)\) meaning that in terms of the highest weight vector \(b_\infty\in B(\infty)\)
  \[
    b(\tau) = f_2f_1f_3f_4f_3f_2f_5f_4 b_\infty\,. 
  \]
  Note that we are once again identifying \(B(\lambda) \) with a subset of \(B(\infty)\) such that the highest weight tableau 
  \[
    \young(13,25,4,6)
  \]
  goes to the highest weight vector \(b_\infty\). 

  Since \(\wt(b(\tau)) = 0 \) we need only consider the zero weight space of \(B(\omega_2 + \omega_4)\subset B(\infty)\):
  %
  \begin{table}[ht]
    \centering
    \newcolumntype{Q}{>{$}l<{$}>{$}l<{$}}
    \begin{tabular}{Q} 
      b & \vec\varepsilon = (\varepsilon_1,\varepsilon_2,\varepsilon_3,\varepsilon_4,\varepsilon_5)\\
      \midrule 
      b_1 & (0,0,0,1,0) \\
      b_2 & (0,0,1,0,1) \\
      b_3 & (0,1,0,0,1) \\ 
      b_4 & (1,0,0,0,1) \\
      b_5 & (0,0,1,0,0) \\
      b_6 = b(\tau) & (0,1,0,1,0) \\
      b_7 & (1,0,0,1,0) \\
      b_8 & (0,1,0,0,0) \\
      b_9 & (1,0,1,0,0) 
    \end{tabular}
  \end{table}
From this table it is immediate that \(b_i \not \le_{\text{str}} b_6  \) for any \(i\ne 6\) except perhaps \(i = 1\) or \(i = 8\). To see that \(b_1\not\le_{\text{str}} b_6 \) and \(b_8\not\le_{\text{str}} b_6 \) we use Sage to determine \(f\) strings for each of \(b_1,b_8\).

Since 
\[
  \begin{aligned}
    \varepsilon_4(b_1) &= \varepsilon_4(b_6) = 1 \\
    \varepsilon_3(f_4b_1) &= \varepsilon_3(f_4 b_6) = 0 \\
    \varepsilon_2(f_3f_4b_1) &= \varepsilon_2(f_3f_4 b_6) = 0 \\
    \varepsilon_1(f_2f_3f_4b_1) &= \varepsilon_1(f_2f_3f_4b_6) = 0 \\
    \varepsilon_5(f_1f_2f_3f_4b_1) &= \varepsilon_5(f_1f_2f_3f_4b_6) = 0 
  \end{aligned}
\]
and 
\[
  \varepsilon_4(f_5f_1f_2f_3f_4b_1) = 1 \ge \varepsilon_4(f_5f_1f_2f_3f_4b_6) = 0 
\]
it follows that \(m(b_1,b_6) = 0\).

Similarly, since
\[
  \begin{aligned}
  \varepsilon_2(b_6) &= \varepsilon_2(b_8) = 1 \\
  \varepsilon_1(f_2b_6) &= \varepsilon_1(f_2b_8) = 0 \\
  \varepsilon_3(f_1f_2b_6) &= \varepsilon_3(f_1f_2b_8) = 0 
  \end{aligned}
\]
while 
\[
  \varepsilon_2(f_3f_1f_2b_6) = 0 < \varepsilon_2(f_3f_1f_2b_8) = 1 
\]
it follows that \(m(b_6,b_8) = 0 \). 
%
\end{proof}
%
%
%
As a corollary, it follows that \(\barD(b_{Z}) = \barD(c_Y)\). Once we have the ideal of the orbital variety which is open dense in \(t^\mu Z = Z\)
we can check this directly. 

For our main result, we show that \(Z\) and \(Y\) are extra-compatible to order 1. 
\begin{theorem}\label{thm:a12}
  \cite[Theorem A.12]{baumann2019mirkovic}
Let \(M \) be a general point of \(Y\). 
  For all \(\nu\in Q_+\), 
  \begin{equation}
    \label{thm:121}
    \dim H^0(Z,\scL)_{\nu} = \dim H^\bullet(F_{1,\nu}(M))\,.
  \end{equation}
\end{theorem}
%
To prove this theorem we again begin by describing the generalized orbital variety of \(Z\). 
%
%
%
\subsection{Orbital variety}
If \(A\in\TT_\mu\cap\n = \n\) has coordinate
$$
A = 
\left[\begin{BMAT}{cccccc}{cccccc}
  0 & a_1 & a_2 & a_3 & a_4 & a_5\\
  0 & 0 & a_6 & a_7 & a_8 & a_9\\
  0 & 0 & 0 & a_{10} & a_{11} & a_{12}\\
  0 & 0 & 0 & 0 & a_{13} & a_{14}\\
  0 & 0 & 0 & 0 & 0 & a_{15}\\
  0 & 0 & 0 & 0 & 0 & 0 
\end{BMAT}\right]
$$
then the generalized orbital variety \(X\) which is labeled by \(\tau\) is the vanishing locus of the prime ideal \(I(X)\) computed and verified with the aid of Macaulay2. 

\begin{table}[ht]
  \centering
  \newcolumntype{Q}{>{$}l<{$}}
  \begin{tabular}{Q} 
  {a}_{15}\\
  {a}_{10}\\
  {a}_{1}\\
  {a}_{6}{a}_{11}+{a}_{7}{a}_{13}\\
  {a}_{12}{a}_{13}-{a}_{11}{a}_{14}\\
  {a}_{6}{a}_{12}+{a}_{7}{a}_{14}\\
  {a}_{2}{a}_{11}+{a}_{3}{a}_{13}\\
  {a}_{3}{a}_{6}-{a}_{2}{a}_{7}\\
  {a}_{2}{a}_{12}+{a}_{3}{a}_{14}\\
  {a}_{5}{a}_{6}{a}_{13}-{a}_{2}{a}_{9}{a}_{13}-{a}_{4}{a}_{6}{a}_{14}+{a}_{2}{a}_{8}{a}_{14}\\
  {a}_{5}{a}_{7}{a}_{13}-{a}_{3}{a}_{9}{a}_{13}-{a}_{4}{a}_{7}{a}_{14}+{a}_{3}{a}_{8}{a}_{14}\\
  {a}_{5}{a}_{7}{a}_{11}-{a}_{3}{a}_{9}{a}_{11}-{a}_{4}{a}_{7}{a}_{12}+{a}_{3}{a}_{8}{a}_{12}\\
  \end{tabular}
  \caption{Generators of \(I(X)\)}
\end{table} 
One can compute the Hilbert series of \(I(X)\) as with the previous example.
\begin{table}[ht!]
  \centering
  \ra{1.1}
  \begin{longtable}{@{}lHlH@{}}
  \toprule
  $C$ & Count & $\nu$  & $\dimvec$ \\
  \midrule
  $(1, 2, 3, 4, 5, 6)$ & $(1, 1, 1, 1, 1, 1)$ & $(0, 0, 0, 0, 0, 0)$ & $(0, 0, 0, 0, 0)$ \\
  $(1, 2, 3, 4, 5, \bar 1)$ & $(2, 1, 1, 1, 1, 0)$ & $(1, 0, 0, 0, 0, -1)$ & $(1, 1, 1, 1, 1)$ \\
  $(1, 2, 3, 4, 6, \bar 1)$ & $(2, 1, 1, 1, 0, 1)$ & $(1, 0, 0, 0, -1, 0)$ & $(1, 1, 1, 1, 0)$ \\
  $(1, 2, 3, 5, 6, \bar 1)$ & $(2, 1, 1, 0, 1, 1)$ & $(1, 0, 0, -1, 0, 0)$ & $(1, 1, 1, 0, 0)$ \\
  $(1, 2, 4, 5, 6, \bar 1)$ & $(2, 1, 0, 1, 1, 1)$ & $(1, 0, -1, 0, 0, 0)$ & $(1, 1, 0, 0, 0)$ \\
  $(1, 2, 3, 4, 5, \bar 2)$ & $(1, 2, 1, 1, 1, 0)$ & $(0, 1, 0, 0, 0, -1)$ & $(0, 1, 1, 1, 1)$ \\
  $(1, 2, 3, 4, 6, \bar 2)$ & $(1, 2, 1, 1, 0, 1)$ & $(0, 1, 0, 0, -1, 0)$ & $(0, 1, 1, 1, 0)$ \\
  $(1, 2, 3, 5, 6, \bar 2)$ & $(1, 2, 1, 0, 1, 1)$ & $(0, 1, 0, -1, 0, 0)$ & $(0, 1, 1, 0, 0)$ \\
  $(1, 2, 4, 5, 6, \bar 2)$ & $(1, 2, 0, 1, 1, 1)$ & $(0, 1, -1, 0, 0, 0)$ & $(0, 1, 0, 0, 0)$ \\
  $(1, 2, 3, 4, \bar 1, \bar 2)$ & $(2, 2, 1, 1, 0, 0)$ & $(1, 1, 0, 0, -1, -1)$ & $(1, 2, 2, 2, 1)$ \\
  $(1, 2, 3, 5, \bar 1, \bar 2)$ & $(2, 2, 1, 0, 1, 0)$ & $(1, 1, 0, -1, 0, -1)$ & $(1, 2, 2, 1, 1)$ \\
  $(1, 2, 4, 5, \bar 1, \bar 2)$ & $(2, 2, 0, 1, 1, 0)$ & $(1, 1, -1, 0, 0, -1)$ & $(1, 2, 1, 1, 1)$ \\
  $(1, 2, 3, 6, \bar 1, \bar 2)$ & $(2, 2, 1, 0, 0, 1)$ & $(1, 1, 0, -1, -1, 0)$ & $(1, 2, 2, 1, 0)$ \\
  $(1, 2, 4, 6, \bar 1, \bar 2)$ & $(2, 2, 0, 1, 0, 1)$ & $(1, 1, -1, 0, -1, 0)$ & $(1, 2, 1, 1, 0)$ \\
  $(1, 2, 3, 4, 5, \bar 3)$ & $(1, 1, 2, 1, 1, 0)$ & $(0, 0, 1, 0, 0, -1)$ & $(0, 0, 1, 1, 1)$ \\
  $(1, 2, 3, 4, 6, \bar 3)$ & $(1, 1, 2, 1, 0, 1)$ & $(0, 0, 1, 0, -1, 0)$ & $(0, 0, 1, 1, 0)$ \\
  $(1, 2, 3, 4, \bar 1, \bar 3)$ & $(2, 1, 2, 1, 0, 0)$ & $(1, 0, 1, 0, -1, -1)$ & $(1, 1, 2, 2, 1)$ \\
  $(1, 2, 3, 5, \bar 1, \bar 3)$ & $(2, 1, 2, 0, 1, 0)$ & $(1, 0, 1, -1, 0, -1)$ & $(1, 1, 2, 1, 1)$ \\
  $(1, 2, 3, 5, \bar 1, 4)$ & $(2, 1, 1, 1, 1, 0)$ & $(1, 0, 0, 0, 0, -1)$ & $(1, 1, 1, 1, 1)$ \\
  $(1, 2, 3, 6, \bar 1, \bar 3)$ & $(2, 1, 2, 0, 0, 1)$ & $(1, 0, 1, -1, -1, 0)$ & $(1, 1, 2, 1, 0)$ \\
  $(1, 2, 3, 6, \bar 1, 4)$ & $(2, 1, 1, 1, 0, 1)$ & $(1, 0, 0, 0, -1, 0)$ & $(1, 1, 1, 1, 0)$ \\
  $(1, 2, 3, 4, \bar 2, \bar 3)$ & $(1, 2, 2, 1, 0, 0)$ & $(0, 1, 1, 0, -1, -1)$ & $(0, 1, 2, 2, 1)$ \\
  $(1, 2, 3, 5, \bar 2, \bar 3)$ & $(1, 2, 2, 0, 1, 0)$ & $(0, 1, 1, -1, 0, -1)$ & $(0, 1, 2, 1, 1)$ \\
  $(1, 2, 3, 5, 2, 4)$ & $(1, 2, 1, 1, 1, 0)$ & $(0, 1, 0, 0, 0, -1)$ & $(0, 1, 1, 1, 1)$ \\
  $(1, 2, 3, 6, \bar 2, \bar 3)$ & $(1, 2, 2, 0, 0, 1)$ & $(0, 1, 1, -1, -1, 0)$ & $(0, 1, 2, 1, 0)$ \\
  $(1, 2, 3, 6, \bar 2, \bar 4)$ & $(1, 2, 1, 1, 0, 1)$ & $(0, 1, 0, 0, -1, 0)$ & $(0, 1, 1, 1, 0)$ \\
  $(1, 2, 3, 4, 5, \bar 4)$ & $(1, 1, 1, 2, 1, 0)$ & $(0, 0, 0, 1, 0, -1)$ & $(0, 0, 0, 1, 1)$ \\
  $(1, 2, 3, 4, 6, \bar 4)$ & $(1, 1, 1, 2, 0, 1)$ & $(0, 0, 0, 1, -1, 0)$ & $(0, 0, 0, 1, 0)$ \\
  $(1, 2, 3, 4, \bar 1, \bar 4)$ & $(2, 1, 1, 2, 0, 0)$ & $(1, 0, 0, 1, -1, -1)$ & $(1, 1, 1, 2, 1)$ \\
  $(1, 2, 4, 5, \bar 1, \bar 4)$ & $(2, 1, 0, 2, 1, 0)$ & $(1, 0, -1, 1, 0, -1)$ & $(1, 1, 0, 1, 1)$ \\
  $(1, 2, 4, 6, \bar 1, \bar 4)$ & $(2, 1, 0, 2, 0, 1)$ & $(1, 0, -1, 1, -1, 0)$ & $(1, 1, 0, 1, 0)$ \\
  $(1, 2, 3, 4, \bar 2, \bar 4)$ & $(1, 2, 1, 2, 0, 0)$ & $(0, 1, 0, 1, -1, -1)$ & $(0, 1, 1, 2, 1)$ \\
  $(1, 2, 4, 5, \bar 2, \bar 4)$ & $(1, 2, 0, 2, 1, 0)$ & $(0, 1, -1, 1, 0, -1)$ & $(0, 1, 0, 1, 1)$ \\
  $(1, 2, 4, 6, \bar 2, \bar 4)$ & $(1, 2, 0, 2, 0, 1)$ & $(0, 1, -1, 1, -1, 0)$ & $(0, 1, 0, 1, 0)$ \\
  $(1, 2, 3, 4, \bar 3, \bar 4)$ & $(1, 1, 2, 2, 0, 0)$ & $(0, 0, 1, 1, -1, -1)$ & $(0, 0, 1, 2, 1)$ \\
  \midrule
  & & 35 & \\ 
  \bottomrule
  \end{longtable}
  \caption{\(H^0(Z,\scL)\)
    \label{tab:indexbasis}}
  \end{table}

%
%

\begin{table}[ht!]
  \centering
  \ra{1.1}
  \begin{longtable}{@{}lll@{}}
  \toprule
  \(\nu\) & $F_1(M)^\nu$ & \(\dim H^\bullet(F_1(M)^\nu)\) \\
  \midrule
  $(1, 2, 2, 1, 0)$ & Point & $1$ \\
  $(1, 1, 0, 1, 1)$ & Point & $1$ \\
  $(0, 1, 1, 2, 1)$ & Point & $1$ \\
  $(1, 1, 2, 1, 1)$ & Point & $1$ \\
  $(1, 1, 1, 2, 1)$ & Point & $1$ \\
  $(0, 0, 1, 1, 1)$ & Point & $1$ \\
  $(0, 0, 1, 2, 1)$ & Point & $1$ \\
  $(0, 0, 0, 1, 1)$ & Point & $1$ \\
  $(1, 1, 1, 1, 1)$ & $\mathbb{P}^1$ & $2$ \\
  $(0, 1, 1, 1, 0)$ & $\mathbb{P}^1$ & $2$ \\
  $(0, 1, 0, 0, 0)$ & Point & $1$ \\
  $(0, 1, 1, 0, 0)$ & Point & $1$ \\
  $(1, 2, 2, 1, 1)$ & Point & $1$ \\
  $(0, 1, 0, 1, 0)$ & Point & $1$ \\
  $(1, 2, 2, 2, 1)$ & Point & $1$ \\
  $(1, 2, 1, 1, 1)$ & Point & $1$ \\
  $(1, 1, 1, 0, 0)$ & Point & $1$ \\
  $(1, 1, 0, 1, 0)$ & Point & $1$ \\
  $(0, 1, 2, 1, 0)$ & Point & $1$ \\
  $(1, 1, 1, 1, 0)$ & $\mathbb{P}^1$ & $2$ \\
  $(1, 1, 0, 0, 0)$ & Point & $1$ \\
  $(1, 1, 2, 1, 0)$ & Point & $1$ \\
  $(0, 0, 0, 1, 0)$ & Point & $1$ \\
  $(0, 0, 0, 0, 0)$ & Point & $1$ \\
  $(0, 0, 1, 1, 0)$ & Point & $1$ \\
  $(0, 1, 1, 1, 1)$ & $\mathbb{P}^1$ & $2$ \\
  $(0, 1, 0, 1, 1)$ & Point & $1$ \\
  $(1, 2, 1, 1, 0)$ & Point & $1$ \\
  $(0, 1, 2, 1, 1)$ & Point & $1$ \\
  $(1, 1, 2, 2, 1)$ & Point & $1$ \\
  $(0, 1, 2, 2, 1)$ & Point & $1$\\
  \midrule
  & & $35$\\
  \bottomrule\\
  \caption{
      Components of the quiver Grassmannian of $Q_\xi$
  \label{tab:qlambdacomponents}
  }
  \end{longtable}
\end{table}
\subsection{Preprojective algebra module}
Next, a general point \(M_a\in Y\) takes the form 
\begin{equation}
  \label{eq:a4m2}
  \begin{tikzpicture}
    \begin{scope}[scale=1.3,xshift=6cm]
    \node (10) at (1,0){$2$};
    \node (30) at (3,0){$4$};
    \node (01) at (0,1){$1$};
    \node (21) at (2,1){$33$};
    \node (41) at (4,1){$5$};
    \node (12) at (1,2){$2$};
    \node (32) at (3,2){$4$};
    \draw[->] (12) to node[pos=.3,left]{$\scriptstyle1$} (01);
    \draw[->] (12) to node[pos=.6,left]{$\left[\bsm1\\0\esm\right]$} (21);
    \draw[->] (32) to node[pos=.6,right]{$\left[\bsm0\\1\esm\right]$} (21);
    \draw[->] (32) to node[pos=.3,right]{$\scriptstyle1$} (41);
    \draw[->] (01) to node[pos=.6,left]{$\scriptstyle-a$} (10);
    \draw[->] (21) to node[pos=.3,left]{$\left[\bsm a&1\esm\right]$} (10);
    \draw[->] (21) to node[pos=.3,right]{$\left[\bsm1&1\esm\right]$} (30);
    \draw[->] (41) to node[pos=.6,right]{$\scriptstyle-1$} (30);
    \end{scope}
  \end{tikzpicture}
\end{equation}
It is straightforward to find all submodules of \(M_a\) and thus to describe the space \(F_1(M_a)\). Comparing with the computation of \(H^0(Z,\scL)\) via Table~\ref{tab:indexbasis} yields the proof of Theorem~\ref{thm:a12}.

\subsection{Equal $\barD$}
%
%
To compute \(\barD(c_Y)\) we again enumerate composition series. We find 148 sequences \(\vi\) such that \(F_\vi(M_a)\cong\PP^1\) and 104 sequences \(\vi\) such that \(F_\vi(M_a)\) is a point. Using a computer, we check that the resulting polynomial \(p(\mu) \barD(c_Y)\) agrees with the multidegree of \(X\),
\[
  \begin{aligned}
    \mdeg_{\n}(X) &= {\alpha}_{1}{\alpha}_{3}{\alpha}_{5}
    \left(
      2\alpha_{1,4}\alpha_{1,5}\alpha_2\alpha_{3,5} +{\alpha}_{1,3}{\alpha}_{1,5}{\alpha}_{3,4}{\alpha}_{3,5}+{\alpha}_{1,2}{\alpha}_{15}{\alpha}_{2}{\alpha}_{4}+{\alpha}_{1,2}{\alpha}_{2}{\alpha}_{2,3}{\alpha}_{4}\right. \\
      & \qquad \left. + {\alpha}_{1,2}{\alpha}_{2}{\alpha}_{3,4}{\alpha}_{4}+{\alpha}_{1,3}{\alpha}_{34}{\alpha}_{3,5}{\alpha}_{4}+{\alpha}_{2}{\alpha}_{3,4}{\alpha}_{3,5}{\alpha}_{4}
    \right)
  \end{aligned}
\]
%
%

%
\section{Case 3: Nonequality of basis vectors in \texorpdfstring{\(A_5\)}{A_5}}
%
Take 
\(\lambda = (4,4,2,2)\) and \(\mu = (2,2,2,2,2,2)\). Consider the ``double'' of the tableau from the previous case.
\[
  \tau = \young(1133,2255,44,66)
\]

Let \(Y\in\irr\Lambda\) have Lusztig data equal to \(n(\tau)\).
A general point in \(Y\) is of the form \(M_a\oplus M_{a'}\) with \(a\ne a' \in\CC\). Let \(I(\omega_2 + \omega_4) = I(\omega_2) \oplus I(\omega_4)\) where \(I(\omega_2)\) and \(I(\omega_4)\) are the standard injective \(\cA\)-modules (see Figure~\ref{fig:iw24}).
\begin{figure}[ht!]
  \centering
  \begin{tikzpicture}
    \begin{scope}[scale=1,xshift=6cm]
      \node (03) at (0,3) {$1$};
      \node (12) at (1,2) {$2$};
      \node (21) at (2,1) {$3$}; 
      \node (30) at (3,0) {$4$}; 
      \node (14) at (1,4) {$2$}; 
      \node (23) at (2,3) {$3$}; 
      \node (32) at (3,2) {$4$}; 
      \node (41) at (4,1) {$5$}; 
      \draw[->] (03) -- (12); 
      \draw[->] (12) -- (21); 
      \draw[->] (21) -- (30); 
      \draw[->] (14) -- (23);
      \draw[->] (23) -- (32);
      \draw[->] (32) -- (41); 
      \draw[->] (14) -- (03); \draw[->] (23) -- (12); \draw[->] (32) -- (21); \draw[->] (41) -- (30); 
      \node at (2,-1) {$I(\omega_2)$};
    \end{scope}
  \end{tikzpicture} 
  \qquad
  \begin{tikzpicture}
    \begin{scope}[scale=1,xshift=6cm] 
      \node (01) at (0,1) {$1$}; 
      \node (12) at (1,2) {$2$}; 
      \node (23) at (2,3) {$3$}; 
      \node (34) at (3,4) {$4$}; 
      \draw[->] (34) -- (23);
      \draw[->] (23) -- (12);
      \draw[->] (12) -- (01);
      \node (10) at (1,0) {$2$};
      \node (21) at (2,1) {$3$}; 
      \node (32) at (3,2) {$4$}; 
      \node (43) at (4,3) {$5$};
      \draw[->] (43) -- (32); 
      \draw[->] (32) -- (21);
      \draw[->] (21) -- (10); 
      \draw[->] (34) -- (43); \draw[->] (23) -- (32); \draw[->] (12) -- (21); \draw[->] (01) -- (10); 
      \node at (2,-1) {$I(\omega_4)$};
    \end{scope}
  \end{tikzpicture}
  \caption{Standard injective \(\cA\)-modules
  \label{fig:iw24}}
\end{figure}
\begin{theorem}
  \label{thm:adcmf}
  \cite[Theorem A.13]{baumann2019mirkovic}
Let \(Y'\in\irr\Lambda(\nu/2)\) be the component from Equation~\ref{eq:a4m2}, and let \(Y''\in\irr\Lambda(\nu)\) be the component whose general point is \(I(\omega_2 + \omega_4)\). 
Then 
\[
  \barD (b_Z) = \barD(c_{Y'})^2 - 2\barD(c_{Y''})\,.
\]
In particular, \(\barD(b_Z) \ne \barD( c_Y)\), and therefore \(b_Z\ne c_Y\). 
\end{theorem}
\begin{proof}
By Lemma~\ref{lem:glshom} and multiplicativity of \(\barD\) the computation of the right-hand side is reduced to the previous section (and the easy computation of \(\barD\) with respect to \(I(\omega_i)\)). 

For the left-hand side, we use Equation~\ref{eq:boxyXt} to give a description of the generalized orbital variety \(X_\tau\). 
With the help of a computer we find that it is the vanishing locus of a prime ideal \(I\) seen in Figure~\ref{yogi} of dimension \(16\) inside a polynomial ring with \(24\) generators. 
%

From there, it is easy to compute the multidegree of \(X_\tau\) and thus \(\barD(b_Z)\). 
\end{proof}
\begin{figure}[ht]
    \centering
    {\small
\hspace*{-1.3cm}$\bgroup\begin{matrix}
{b}_{15},{b}_{10},{b}_{1},{a}_{15},{a}_{10},{a}_{1}\\
{a}_{14} {b}_{11}-{a}_{13} {b}_{12}-{a}_{12} {b}_{13}+{a}_{11} {b}_{14}\\
{a}_{12} {b}_{6}+{a}_{14} {b}_{7}+{a}_{6} {b}_{12}+{a}_{7} {b}_{14}\\
{a}_{11} {b}_{6}+{a}_{13} {b}_{7}+{a}_{6} {b}_{11}+{a}_{7} {b}_{13}\\
{a}_{12} {b}_{2}+{a}_{14} {b}_{3}+{a}_{2} {b}_{12}+{a}_{3} {b}_{14}\\
{a}_{11} {b}_{2}+{a}_{13} {b}_{3}+{a}_{2} {b}_{11}+{a}_{3} {b}_{13}\\
{a}_{7} {b}_{2}-{a}_{6} {b}_{3}-{a}_{3} {b}_{6}+{a}_{2} {b}_{7}\\
{a}_{12} {a}_{13}-{a}_{11} {a}_{14}\\
{a}_{6} {a}_{12}+{a}_{7} {a}_{14}\\
{a}_{2} {a}_{12}+{a}_{3} {a}_{14}\\
{a}_{6} {a}_{11}+{a}_{7} {a}_{13}\\
{a}_{2} {a}_{11}+{a}_{3} {a}_{13}\\
{a}_{3} {a}_{6}-{a}_{2} {a}_{7}\\
{a}_{5} {a}_{6} {b}_{11}-{a}_{2} {a}_{9} {b}_{11}-{a}_{4} {a}_{6} {b}_{12}+{a}_{2} {a}_{8} {b}_{12}+{a}_{5} {a}_{7} {b}_{13}-{a}_{3} {a}_{9} {b}_{13}-{a}_{4} {a}_{7} {b}_{14}+{a}_{3} {a}_{8} {b}_{14}\\
{a}_{9} {a}_{13} {b}_{3}-{a}_{8} {a}_{14} {b}_{3}+{a}_{7} {a}_{14} {b}_{4}-{a}_{7} {a}_{13} {b}_{5}-{a}_{5} {a}_{13} {b}_{7}+{a}_{4} {a}_{14} {b}_{7}-{a}_{3} {a}_{14} {b}_{8}+{a}_{3} {a}_{13} {b}_{9}-{a}_{5} {a}_{7} {b}_{13}+{a}_{3} {a}_{9} {b}_{13}+{a}_{4} {a}_{7} {b}_{14}-{a}_{3} {a}_{8} {b}_{14}\\
{a}_{9} {a}_{11} {b}_{3}-{a}_{8} {a}_{12} {b}_{3}+{a}_{7} {a}_{12} {b}_{4}-{a}_{7} {a}_{11} {b}_{5}-{a}_{5} {a}_{11} {b}_{7}+{a}_{4} {a}_{12} {b}_{7}-{a}_{3} {a}_{12} {b}_{8}+{a}_{3} {a}_{11} {b}_{9}-{a}_{5} {a}_{7} {b}_{11}+{a}_{3} {a}_{9} {b}_{11}+{a}_{4} {a}_{7} {b}_{12}-{a}_{3} {a}_{8} {b}_{12}\\
{a}_{9} {a}_{13} {b}_{2}-{a}_{8} {a}_{14} {b}_{2}+{a}_{6} {a}_{14} {b}_{4}-{a}_{6} {a}_{13} {b}_{5}-{a}_{5} {a}_{13} {b}_{6}+{a}_{4} {a}_{14} {b}_{6}-{a}_{2} {a}_{14} {b}_{8}+{a}_{2} {a}_{13} {b}_{9}-{a}_{5} {a}_{6} {b}_{13}+{a}_{2} {a}_{9} {b}_{13}+{a}_{4} {a}_{6} {b}_{14}-{a}_{2} {a}_{8} {b}_{14}\\
{a}_{5} {a}_{7} {a}_{13}-{a}_{3} {a}_{9} {a}_{13}-{a}_{4} {a}_{7} {a}_{14}+{a}_{3} {a}_{8} {a}_{14}\\
{a}_{5} {a}_{6} {a}_{13}-{a}_{2} {a}_{9} {a}_{13}-{a}_{4} {a}_{6} {a}_{14}+{a}_{2} {a}_{8} {a}_{14}\\
{a}_{5} {a}_{7} {a}_{11}-{a}_{3} {a}_{9} {a}_{11}-{a}_{4} {a}_{7} {a}_{12}+{a}_{3} {a}_{8} {a}_{12}
\end{matrix}\egroup$}

\centering
$
{a}_{9}^{2} {b}_{2}^{2} {b}_{11}^{2}-{a}_{6}^{2} {b}_{5}^{2} {b}_{11}^{2}-2 {a}_{5} {a}_{9} {b}_{2} {b}_{6} {b}_{11}^{2}+{a}_{5}^{2} {b}_{6}^{2} {b}_{11}^{2}+2 {a}_{2} {a}_{6} {b}_{5} {b}_{9} {b}_{11}^{2}-{a}_{2}^{2} {b}_{9}^{2} {b}_{11}^{2}-2 {a}_{8} {a}_{9} {b}_{2}^{2} {b}_{11} {b}_{12}+
2 {a}_{6}^{2} {b}_{4} {b}_{5} {b}_{11} {b}_{12}+
2 {a}_{5} {a}_{8} {b}_{2} {b}_{6} {b}_{11} {b}_{12}+2 {a}_{4} {a}_{9} {b}_{2} {b}_{6} {b}_{11} {b}_{12}-2 {a}_{4} {a}_{5} {b}_{6}^{2} {b}_{11} {b}_{12}-2 {a}_{2} {a}_{6} {b}_{5} {b}_{8} {b}_{11} {b}_{12}-2 {a}_{2} {a}_{6} {b}_{4} {b}_{9} {b}_{11} {b}_{12}+2 {a}_{2}^{2} {b}_{8} {b}_{9} {b}_{11} {b}_{12}+{a}_{8}^{2} {b}_{2}^{2} {b}_{12}^{2}-{a}_{6}^{2} {b}_{4}^{2} {b}_{12}^{2}-2 {a}_{4} {a}_{8} {b}_{2} {b}_{6} {b}_{12}^{2}+{a}_{4}^{2} {b}_{6}^{2} {b}_{12}^{2}+2 {a}_{2} {a}_{6} {b}_{4} {b}_{8} {b}_{12}^{2}-{a}_{2}^{2} {b}_{8}^{2} {b}_{12}^{2}+2 {a}_{9}^{2} {b}_{2} {b}_{3} {b}_{11} {b}_{13}-2 {a}_{6} {a}_{7} {b}_{5}^{2} {b}_{11} {b}_{13}-2 {a}_{5} {a}_{9} {b}_{3} {b}_{6} {b}_{11} {b}_{13}-2 {a}_{5} {a}_{9} {b}_{2} {b}_{7} {b}_{11} {b}_{13}+2 {a}_{5}^{2} {b}_{6} {b}_{7} {b}_{11} {b}_{13}+4 {a}_{2} {a}_{7} {b}_{5} {b}_{9} {b}_{11} {b}_{13}-2 {a}_{2} {a}_{3} {b}_{9}^{2} {b}_{11} {b}_{13}-2 {a}_{8} {a}_{9} {b}_{2} {b}_{3} {b}_{12} {b}_{13}+2 {a}_{6} {a}_{7} {b}_{4} {b}_{5} {b}_{12} {b}_{13}+4 {a}_{5} {a}_{8} {b}_{3} {b}_{6} {b}_{12} {b}_{13}-2 {a}_{4} {a}_{9} {b}_{3} {b}_{6} {b}_{12} {b}_{13}-2 {a}_{5} {a}_{8} {b}_{2} {b}_{7} {b}_{12} {b}_{13}+4 {a}_{4} {a}_{9} {b}_{2} {b}_{7} {b}_{12} {b}_{13}-2 {a}_{4} {a}_{5} {b}_{6} {b}_{7} {b}_{12} {b}_{13}-2 {a}_{2} {a}_{7} {b}_{5} {b}_{8} {b}_{12} {b}_{13}-2 {a}_{2} {a}_{7} {b}_{4} {b}_{9} {b}_{12} {b}_{13}+2 {a}_{2} {a}_{3} {b}_{8} {b}_{9} {b}_{12} {b}_{13}+{a}_{9}^{2} {b}_{3}^{2} {b}_{13}^{2}-{a}_{7}^{2} {b}_{5}^{2} {b}_{13}^{2}-2 {a}_{5} {a}_{9} {b}_{3} {b}_{7} {b}_{13}^{2}+{a}_{5}^{2} {b}_{7}^{2} {b}_{13}^{2}+2 {a}_{3} {a}_{7} {b}_{5} {b}_{9} {b}_{13}^{2}-{a}_{3}^{2} {b}_{9}^{2} {b}_{13}^{2}-2 {a}_{8} {a}_{9} {b}_{2} {b}_{3} {b}_{11} {b}_{14}+2 {a}_{6} {a}_{7} {b}_{4} {b}_{5} {b}_{11} {b}_{14}-2 {a}_{5} {a}_{8} {b}_{3} {b}_{6} {b}_{11} {b}_{14}+4 {a}_{4} {a}_{9} {b}_{3} {b}_{6} {b}_{11} {b}_{14}+4 {a}_{5} {a}_{8} {b}_{2} {b}_{7} {b}_{11} {b}_{14}-2 {a}_{4} {a}_{9} {b}_{2} {b}_{7} {b}_{11} {b}_{14}-2 {a}_{4} {a}_{5} {b}_{6} {b}_{7} {b}_{11} {b}_{14}-2 {a}_{2} {a}_{7} {b}_{5} {b}_{8} {b}_{11} {b}_{14}-2 {a}_{2} {a}_{7} {b}_{4} {b}_{9} {b}_{11} {b}_{14}+2 {a}_{2} {a}_{3} {b}_{8} {b}_{9} {b}_{11} {b}_{14}+2 {a}_{8}^{2} {b}_{2} {b}_{3} {b}_{12} {b}_{14}-2 {a}_{6} {a}_{7} {b}_{4}^{2} {b}_{12} {b}_{14}-2 {a}_{4} {a}_{8} {b}_{3} {b}_{6} {b}_{12} {b}_{14}-2 {a}_{4} {a}_{8} {b}_{2} {b}_{7} {b}_{12} {b}_{14}+2 {a}_{4}^{2} {b}_{6} {b}_{7} {b}_{12} {b}_{14}+4 {a}_{2} {a}_{7} {b}_{4} {b}_{8} {b}_{12} {b}_{14}-2 {a}_{2} {a}_{3} {b}_{8}^{2} {b}_{12} {b}_{14}-2 {a}_{8} {a}_{9} {b}_{3}^{2} {b}_{13} {b}_{14}+2 {a}_{7}^{2} {b}_{4} {b}_{5} {b}_{13} {b}_{14}+2 {a}_{5} {a}_{8} {b}_{3} {b}_{7} {b}_{13} {b}_{14}+2 {a}_{4} {a}_{9} {b}_{3} {b}_{7} {b}_{13} {b}_{14}-2 {a}_{4} {a}_{5} {b}_{7}^{2} {b}_{13} {b}_{14}-2 {a}_{3} {a}_{7} {b}_{5} {b}_{8} {b}_{13} {b}_{14}-2 {a}_{3} {a}_{7} {b}_{4} {b}_{9} {b}_{13} {b}_{14}+2 {a}_{3}^{2} {b}_{8} {b}_{9} {b}_{13} {b}_{14}+{a}_{8}^{2} {b}_{3}^{2} {b}_{14}^{2}-{a}_{7}^{2} {b}_{4}^{2} {b}_{14}^{2}-2 {a}_{4} {a}_{8} {b}_{3} {b}_{7} {b}_{14}^{2}+{a}_{4}^{2} {b}_{7}^{2} {b}_{14}^{2}+2 {a}_{3} {a}_{7} {b}_{4} {b}_{8} {b}_{14}^{2}-{a}_{3}^{2} {b}_{8}^{2} {b}_{14}^{2}
$
    \caption{Ideal of \(X_\tau\)}
    \label{yogi}
\end{figure}

%
\begin{figure}[ht]
    \centering
    \small $({{\alpha}_{5}})^{2} ({{\alpha}_{3}})^{2} ({{\alpha}_{1}})^{2} ({\alpha}_{1}^{4} {\alpha}_{2}^{2} {\alpha}_{3}^{2}+6 {\alpha}_{1}^{3} {\alpha}_{2}^{3} {\alpha}_{3}^{2}+13 {\alpha}_{1}^{2} {\alpha}_{2}^{4} {\alpha}_{3}^{2}+12 {\alpha}_{1} {\alpha}_{2}^{5} {\alpha}_{3}^{2}+4 {\alpha}_{2}^{6} {\alpha}_{3}^{2}+2 {\alpha}_{1}^{4} {\alpha}_{2} {\alpha}_{3}^{3}+16 {\alpha}_{1}^{3} {\alpha}_{2}^{2} {\alpha}_{3}^{3}+44 {\alpha}_{1}^{2} {\alpha}_{2}^{3} {\alpha}_{3}^{3}+50 {\alpha}_{1} {\alpha}_{2}^{4} {\alpha}_{3}^{3}+20 {\alpha}_{2}^{5} {\alpha}_{3}^{3}+{\alpha}_{1}^{4} {\alpha}_{3}^{4}+14 {\alpha}_{1}^{3} {\alpha}_{2} {\alpha}_{3}^{4}+55 {\alpha}_{1}^{2} {\alpha}_{2}^{2} {\alpha}_{3}^{4}+82 {\alpha}_{1} {\alpha}_{2}^{3} {\alpha}_{3}^{4}+41 {\alpha}_{2}^{4} {\alpha}_{3}^{4}+4 {\alpha}_{1}^{3} {\alpha}_{3}^{5}+30 {\alpha}_{1}^{2} {\alpha}_{2} {\alpha}_{3}^{5}+66 {\alpha}_{1} {\alpha}_{2}^{2} {\alpha}_{3}^{5}+44 {\alpha}_{2}^{3} {\alpha}_{3}^{5}+6 {\alpha}_{1}^{2} {\alpha}_{3}^{6}+26 {\alpha}_{1} {\alpha}_{2} {\alpha}_{3}^{6}+26 {\alpha}_{2}^{2} {\alpha}_{3}^{6}+4 {\alpha}_{1} {\alpha}_{3}^{7}+8 {\alpha}_{2} {\alpha}_{3}^{7}+{\alpha}_{3}^{8}+4 {\alpha}_{1}^{4} {\alpha}_{2}^{2} {\alpha}_{3} {\alpha}_{4}+24 {\alpha}_{1}^{3} {\alpha}_{2}^{3} {\alpha}_{3} {\alpha}_{4}+52 {\alpha}_{1}^{2} {\alpha}_{2}^{4} {\alpha}_{3} {\alpha}_{4}+48 {\alpha}_{1} {\alpha}_{2}^{5} {\alpha}_{3} {\alpha}_{4}+16 {\alpha}_{2}^{6} {\alpha}_{3} {\alpha}_{4}+8 {\alpha}_{1}^{4} {\alpha}_{2} {\alpha}_{3}^{2} {\alpha}_{4}+68 {\alpha}_{1}^{3} {\alpha}_{2}^{2} {\alpha}_{3}^{2} {\alpha}_{4}+192 {\alpha}_{1}^{2} {\alpha}_{2}^{3} {\alpha}_{3}^{2} {\alpha}_{4}+220 {\alpha}_{1} {\alpha}_{2}^{4} {\alpha}_{3}^{2} {\alpha}_{4}+88 {\alpha}_{2}^{5} {\alpha}_{3}^{2} {\alpha}_{4}+4 {\alpha}_{1}^{4} {\alpha}_{3}^{3} {\alpha}_{4}+64 {\alpha}_{1}^{3} {\alpha}_{2} {\alpha}_{3}^{3} {\alpha}_{4}+264 {\alpha}_{1}^{2} {\alpha}_{2}^{2} {\alpha}_{3}^{3} {\alpha}_{4}+400 {\alpha}_{1} {\alpha}_{2}^{3} {\alpha}_{3}^{3} {\alpha}_{4}+200 {\alpha}_{2}^{4} {\alpha}_{3}^{3} {\alpha}_{4}+20 {\alpha}_{1}^{3} {\alpha}_{3}^{4} {\alpha}_{4}+160 {\alpha}_{1}^{2} {\alpha}_{2} {\alpha}_{3}^{4} {\alpha}_{4}+360 {\alpha}_{1} {\alpha}_{2}^{2} {\alpha}_{3}^{4} {\alpha}_{4}+240 {\alpha}_{2}^{3} {\alpha}_{3}^{4} {\alpha}_{4}+36 {\alpha}_{1}^{2} {\alpha}_{3}^{5} {\alpha}_{4}+160 {\alpha}_{1} {\alpha}_{2} {\alpha}_{3}^{5} {\alpha}_{4}+160 {\alpha}_{2}^{2} {\alpha}_{3}^{5} {\alpha}_{4}+28 {\alpha}_{1} {\alpha}_{3}^{6} {\alpha}_{4}+56 {\alpha}_{2} {\alpha}_{3}^{6} {\alpha}_{4}+8 {\alpha}_{3}^{7} {\alpha}_{4}+4 {\alpha}_{1}^{4} {\alpha}_{2}^{2} {\alpha}_{4}^{2}+24 {\alpha}_{1}^{3} {\alpha}_{2}^{3} {\alpha}_{4}^{2}+52 {\alpha}_{1}^{2} {\alpha}_{2}^{4} {\alpha}_{4}^{2}+48 {\alpha}_{1} {\alpha}_{2}^{5} {\alpha}_{4}^{2}+16 {\alpha}_{2}^{6} {\alpha}_{4}^{2}+10 {\alpha}_{1}^{4} {\alpha}_{2} {\alpha}_{3} {\alpha}_{4}^{2}+90 {\alpha}_{1}^{3} {\alpha}_{2}^{2} {\alpha}_{3} {\alpha}_{4}^{2}+260 {\alpha}_{1}^{2} {\alpha}_{2}^{3} {\alpha}_{3} {\alpha}_{4}^{2}+300 {\alpha}_{1} {\alpha}_{2}^{4} {\alpha}_{3} {\alpha}_{4}^{2}+120 {\alpha}_{2}^{5} {\alpha}_{3} {\alpha}_{4}^{2}+6 {\alpha}_{1}^{4} {\alpha}_{3}^{2} {\alpha}_{4}^{2}+106 {\alpha}_{1}^{3} {\alpha}_{2} {\alpha}_{3}^{2} {\alpha}_{4}^{2}+452 {\alpha}_{1}^{2} {\alpha}_{2}^{2} {\alpha}_{3}^{2} {\alpha}_{4}^{2}+692 {\alpha}_{1} {\alpha}_{2}^{3} {\alpha}_{3}^{2} {\alpha}_{4}^{2}+346 {\alpha}_{2}^{4} {\alpha}_{3}^{2} {\alpha}_{4}^{2}+40 {\alpha}_{1}^{3} {\alpha}_{3}^{3} {\alpha}_{4}^{2}+332 {\alpha}_{1}^{2} {\alpha}_{2} {\alpha}_{3}^{3} {\alpha}_{4}^{2}+756 {\alpha}_{1} {\alpha}_{2}^{2} {\alpha}_{3}^{3} {\alpha}_{4}^{2}+504 {\alpha}_{2}^{3} {\alpha}_{3}^{3} {\alpha}_{4}^{2}+88 {\alpha}_{1}^{2} {\alpha}_{3}^{4} {\alpha}_{4}^{2}+396 {\alpha}_{1} {\alpha}_{2} {\alpha}_{3}^{4} {\alpha}_{4}^{2}+396 {\alpha}_{2}^{2} {\alpha}_{3}^{4} {\alpha}_{4}^{2}+80 {\alpha}_{1} {\alpha}_{3}^{5} {\alpha}_{4}^{2}+160 {\alpha}_{2} {\alpha}_{3}^{5} {\alpha}_{4}^{2}+26 {\alpha}_{3}^{6} {\alpha}_{4}^{2}+4 {\alpha}_{1}^{4} {\alpha}_{2} {\alpha}_{4}^{3}+36 {\alpha}_{1}^{3} {\alpha}_{2}^{2} {\alpha}_{4}^{3}+104 {\alpha}_{1}^{2} {\alpha}_{2}^{3} {\alpha}_{4}^{3}+120 {\alpha}_{1} {\alpha}_{2}^{4} {\alpha}_{4}^{3}+48 {\alpha}_{2}^{5} {\alpha}_{4}^{3}+4 {\alpha}_{1}^{4} {\alpha}_{3} {\alpha}_{4}^{3}+76 {\alpha}_{1}^{3} {\alpha}_{2} {\alpha}_{3} {\alpha}_{4}^{3}+328 {\alpha}_{1}^{2} {\alpha}_{2}^{2} {\alpha}_{3} {\alpha}_{4}^{3}+504 {\alpha}_{1} {\alpha}_{2}^{3} {\alpha}_{3} {\alpha}_{4}^{3}+252 {\alpha}_{2}^{4} {\alpha}_{3} {\alpha}_{4}^{3}+40 {\alpha}_{1}^{3} {\alpha}_{3}^{2} {\alpha}_{4}^{3}+336 {\alpha}_{1}^{2} {\alpha}_{2} {\alpha}_{3}^{2} {\alpha}_{4}^{3}+768 {\alpha}_{1} {\alpha}_{2}^{2} {\alpha}_{3}^{2} {\alpha}_{4}^{3}+512 {\alpha}_{2}^{3} {\alpha}_{3}^{2} {\alpha}_{4}^{3}+112 {\alpha}_{1}^{2} {\alpha}_{3}^{3} {\alpha}_{4}^{3}+504 {\alpha}_{1} {\alpha}_{2} {\alpha}_{3}^{3} {\alpha}_{4}^{3}+504 {\alpha}_{2}^{2} {\alpha}_{3}^{3} {\alpha}_{4}^{3}+120 {\alpha}_{1} {\alpha}_{3}^{4} {\alpha}_{4}^{3}+240 {\alpha}_{2} {\alpha}_{3}^{4} {\alpha}_{4}^{3}+44 {\alpha}_{3}^{5} {\alpha}_{4}^{3}+{\alpha}_{1}^{4} {\alpha}_{4}^{4}+20 {\alpha}_{1}^{3} {\alpha}_{2} {\alpha}_{4}^{4}+86 {\alpha}_{1}^{2} {\alpha}_{2}^{2} {\alpha}_{4}^{4}+132 {\alpha}_{1} {\alpha}_{2}^{3} {\alpha}_{4}^{4}+66 {\alpha}_{2}^{4} {\alpha}_{4}^{4}+20 {\alpha}_{1}^{3} {\alpha}_{3} {\alpha}_{4}^{4}+166 {\alpha}_{1}^{2} {\alpha}_{2} {\alpha}_{3} {\alpha}_{4}^{4}+378 {\alpha}_{1} {\alpha}_{2}^{2} {\alpha}_{3} {\alpha}_{4}^{4}+252 {\alpha}_{2}^{3} {\alpha}_{3} {\alpha}_{4}^{4}+78 {\alpha}_{1}^{2} {\alpha}_{3}^{2} {\alpha}_{4}^{4}+346 {\alpha}_{1} {\alpha}_{2} {\alpha}_{3}^{2} {\alpha}_{4}^{4}+346 {\alpha}_{2}^{2} {\alpha}_{3}^{2} {\alpha}_{4}^{4}+100 {\alpha}_{1} {\alpha}_{3}^{3} {\alpha}_{4}^{4}+200 {\alpha}_{2} {\alpha}_{3}^{3} {\alpha}_{4}^{4}+41 {\alpha}_{3}^{4} {\alpha}_{4}^{4}+4 {\alpha}_{1}^{3} {\alpha}_{4}^{5}+32 {\alpha}_{1}^{2} {\alpha}_{2} {\alpha}_{4}^{5}+72 {\alpha}_{1} {\alpha}_{2}^{2} {\alpha}_{4}^{5}+48 {\alpha}_{2}^{3} {\alpha}_{4}^{5}+28 {\alpha}_{1}^{2} {\alpha}_{3} {\alpha}_{4}^{5}+120 {\alpha}_{1} {\alpha}_{2} {\alpha}_{3} {\alpha}_{4}^{5}+120 {\alpha}_{2}^{2} {\alpha}_{3} {\alpha}_{4}^{5}+44 {\alpha}_{1} {\alpha}_{3}^{2} {\alpha}_{4}^{5}+88 {\alpha}_{2} {\alpha}_{3}^{2} {\alpha}_{4}^{5}+20 {\alpha}_{3}^{3} {\alpha}_{4}^{5}+4 {\alpha}_{1}^{2} {\alpha}_{4}^{6}+16 {\alpha}_{1} {\alpha}_{2} {\alpha}_{4}^{6}+16 {\alpha}_{2}^{2} {\alpha}_{4}^{6}+8 {\alpha}_{1} {\alpha}_{3} {\alpha}_{4}^{6}+16 {\alpha}_{2} {\alpha}_{3} {\alpha}_{4}^{6}+4 {\alpha}_{3}^{2} {\alpha}_{4}^{6}+2 {\alpha}_{1}^{4} {\alpha}_{2}^{2} {\alpha}_{3} {\alpha}_{5}+12 {\alpha}_{1}^{3} {\alpha}_{2}^{3} {\alpha}_{3} {\alpha}_{5}+26 {\alpha}_{1}^{2} {\alpha}_{2}^{4} {\alpha}_{3} {\alpha}_{5}+24 {\alpha}_{1} {\alpha}_{2}^{5} {\alpha}_{3} {\alpha}_{5}+8 {\alpha}_{2}^{6} {\alpha}_{3} {\alpha}_{5}+4 {\alpha}_{1}^{4} {\alpha}_{2} {\alpha}_{3}^{2} {\alpha}_{5}+34 {\alpha}_{1}^{3} {\alpha}_{2}^{2} {\alpha}_{3}^{2} {\alpha}_{5}+96 {\alpha}_{1}^{2} {\alpha}_{2}^{3} {\alpha}_{3}^{2} {\alpha}_{5}+110 {\alpha}_{1} {\alpha}_{2}^{4} {\alpha}_{3}^{2} {\alpha}_{5}+44 {\alpha}_{2}^{5} {\alpha}_{3}^{2} {\alpha}_{5}+2 {\alpha}_{1}^{4} {\alpha}_{3}^{3} {\alpha}_{5}+32 {\alpha}_{1}^{3} {\alpha}_{2} {\alpha}_{3}^{3} {\alpha}_{5}+132 {\alpha}_{1}^{2} {\alpha}_{2}^{2} {\alpha}_{3}^{3} {\alpha}_{5}+200 {\alpha}_{1} {\alpha}_{2}^{3} {\alpha}_{3}^{3} {\alpha}_{5}+100 {\alpha}_{2}^{4} {\alpha}_{3}^{3} {\alpha}_{5}+10 {\alpha}_{1}^{3} {\alpha}_{3}^{4} {\alpha}_{5}+80 {\alpha}_{1}^{2} {\alpha}_{2} {\alpha}_{3}^{4} {\alpha}_{5}+180 {\alpha}_{1} {\alpha}_{2}^{2} {\alpha}_{3}^{4} {\alpha}_{5}+120 {\alpha}_{2}^{3} {\alpha}_{3}^{4} {\alpha}_{5}+18 {\alpha}_{1}^{2} {\alpha}_{3}^{5} {\alpha}_{5}+80 {\alpha}_{1} {\alpha}_{2} {\alpha}_{3}^{5} {\alpha}_{5}+80 {\alpha}_{2}^{2} {\alpha}_{3}^{5} {\alpha}_{5}+14 {\alpha}_{1} {\alpha}_{3}^{6} {\alpha}_{5}+28 {\alpha}_{2} {\alpha}_{3}^{6} {\alpha}_{5}+4 {\alpha}_{3}^{7} {\alpha}_{5}+4 {\alpha}_{1}^{4} {\alpha}_{2}^{2} {\alpha}_{4} {\alpha}_{5}+24 {\alpha}_{1}^{3} {\alpha}_{2}^{3} {\alpha}_{4} {\alpha}_{5}+52 {\alpha}_{1}^{2} {\alpha}_{2}^{4} {\alpha}_{4} {\alpha}_{5}+48 {\alpha}_{1} {\alpha}_{2}^{5} {\alpha}_{4} {\alpha}_{5}+16 {\alpha}_{2}^{6} {\alpha}_{4} {\alpha}_{5}+10 {\alpha}_{1}^{4} {\alpha}_{2} {\alpha}_{3} {\alpha}_{4} {\alpha}_{5}+90 {\alpha}_{1}^{3} {\alpha}_{2}^{2} {\alpha}_{3} {\alpha}_{4} {\alpha}_{5}+260 {\alpha}_{1}^{2} {\alpha}_{2}^{3} {\alpha}_{3} {\alpha}_{4} {\alpha}_{5}+300 {\alpha}_{1} {\alpha}_{2}^{4} {\alpha}_{3} {\alpha}_{4} {\alpha}_{5}+120 {\alpha}_{2}^{5} {\alpha}_{3} {\alpha}_{4} {\alpha}_{5}+6 {\alpha}_{1}^{4} {\alpha}_{3}^{2} {\alpha}_{4} {\alpha}_{5}+106 {\alpha}_{1}^{3} {\alpha}_{2} {\alpha}_{3}^{2} {\alpha}_{4} {\alpha}_{5}+452 {\alpha}_{1}^{2} {\alpha}_{2}^{2} {\alpha}_{3}^{2} {\alpha}_{4} {\alpha}_{5}+692 {\alpha}_{1} {\alpha}_{2}^{3} {\alpha}_{3}^{2} {\alpha}_{4} {\alpha}_{5}+346 {\alpha}_{2}^{4} {\alpha}_{3}^{2} {\alpha}_{4} {\alpha}_{5}+40 {\alpha}_{1}^{3} {\alpha}_{3}^{3} {\alpha}_{4} {\alpha}_{5}+332 {\alpha}_{1}^{2} {\alpha}_{2} {\alpha}_{3}^{3} {\alpha}_{4} {\alpha}_{5}+756 {\alpha}_{1} {\alpha}_{2}^{2} {\alpha}_{3}^{3} {\alpha}_{4} {\alpha}_{5}+504 {\alpha}_{2}^{3} {\alpha}_{3}^{3} {\alpha}_{4} {\alpha}_{5}+88 {\alpha}_{1}^{2} {\alpha}_{3}^{4} {\alpha}_{4} {\alpha}_{5}+396 {\alpha}_{1} {\alpha}_{2} {\alpha}_{3}^{4} {\alpha}_{4} {\alpha}_{5}+396 {\alpha}_{2}^{2} {\alpha}_{3}^{4} {\alpha}_{4} {\alpha}_{5}+80 {\alpha}_{1} {\alpha}_{3}^{5} {\alpha}_{4} {\alpha}_{5}+160 {\alpha}_{2} {\alpha}_{3}^{5} {\alpha}_{4} {\alpha}_{5}+26 {\alpha}_{3}^{6} {\alpha}_{4} {\alpha}_{5}+6 {\alpha}_{1}^{4} {\alpha}_{2} {\alpha}_{4}^{2} {\alpha}_{5}+54 {\alpha}_{1}^{3} {\alpha}_{2}^{2} {\alpha}_{4}^{2} {\alpha}_{5}+156 {\alpha}_{1}^{2} {\alpha}_{2}^{3} {\alpha}_{4}^{2} {\alpha}_{5}+180 {\alpha}_{1} {\alpha}_{2}^{4} {\alpha}_{4}^{2} {\alpha}_{5}+72 {\alpha}_{2}^{5} {\alpha}_{4}^{2} {\alpha}_{5}+6 {\alpha}_{1}^{4} {\alpha}_{3} {\alpha}_{4}^{2} {\alpha}_{5}+114 {\alpha}_{1}^{3} {\alpha}_{2} {\alpha}_{3} {\alpha}_{4}^{2} {\alpha}_{5}+492 {\alpha}_{1}^{2} {\alpha}_{2}^{2} {\alpha}_{3} {\alpha}_{4}^{2} {\alpha}_{5}+756 {\alpha}_{1} {\alpha}_{2}^{3} {\alpha}_{3} {\alpha}_{4}^{2} {\alpha}_{5}+378 {\alpha}_{2}^{4} {\alpha}_{3} {\alpha}_{4}^{2} {\alpha}_{5}+60 {\alpha}_{1}^{3} {\alpha}_{3}^{2} {\alpha}_{4}^{2} {\alpha}_{5}+504 {\alpha}_{1}^{2} {\alpha}_{2} {\alpha}_{3}^{2} {\alpha}_{4}^{2} {\alpha}_{5}+1152 {\alpha}_{1} {\alpha}_{2}^{2} {\alpha}_{3}^{2} {\alpha}_{4}^{2} {\alpha}_{5}+768 {\alpha}_{2}^{3} {\alpha}_{3}^{2} {\alpha}_{4}^{2} {\alpha}_{5}+168 {\alpha}_{1}^{2} {\alpha}_{3}^{3} {\alpha}_{4}^{2} {\alpha}_{5}+756 {\alpha}_{1} {\alpha}_{2} {\alpha}_{3}^{3} {\alpha}_{4}^{2} {\alpha}_{5}+756 {\alpha}_{2}^{2} {\alpha}_{3}^{3} {\alpha}_{4}^{2} {\alpha}_{5}+180 {\alpha}_{1} {\alpha}_{3}^{4} {\alpha}_{4}^{2} {\alpha}_{5}+360 {\alpha}_{2} {\alpha}_{3}^{4} {\alpha}_{4}^{2} {\alpha}_{5}+66 {\alpha}_{3}^{5} {\alpha}_{4}^{2} {\alpha}_{5}+2 {\alpha}_{1}^{4} {\alpha}_{4}^{3} {\alpha}_{5}+40 {\alpha}_{1}^{3} {\alpha}_{2} {\alpha}_{4}^{3} {\alpha}_{5}+172 {\alpha}_{1}^{2} {\alpha}_{2}^{2} {\alpha}_{4}^{3} {\alpha}_{5}+264 {\alpha}_{1} {\alpha}_{2}^{3} {\alpha}_{4}^{3} {\alpha}_{5}+132 {\alpha}_{2}^{4} {\alpha}_{4}^{3} {\alpha}_{5}+40 {\alpha}_{1}^{3} {\alpha}_{3} {\alpha}_{4}^{3} {\alpha}_{5}+332 {\alpha}_{1}^{2} {\alpha}_{2} {\alpha}_{3} {\alpha}_{4}^{3} {\alpha}_{5}+756 {\alpha}_{1} {\alpha}_{2}^{2} {\alpha}_{3} {\alpha}_{4}^{3} {\alpha}_{5}+504 {\alpha}_{2}^{3} {\alpha}_{3} {\alpha}_{4}^{3} {\alpha}_{5}+156 {\alpha}_{1}^{2} {\alpha}_{3}^{2} {\alpha}_{4}^{3} {\alpha}_{5}+692 {\alpha}_{1} {\alpha}_{2} {\alpha}_{3}^{2} {\alpha}_{4}^{3} {\alpha}_{5}+692 {\alpha}_{2}^{2} {\alpha}_{3}^{2} {\alpha}_{4}^{3} {\alpha}_{5}+200 {\alpha}_{1} {\alpha}_{3}^{3} {\alpha}_{4}^{3} {\alpha}_{5}+400 {\alpha}_{2} {\alpha}_{3}^{3} {\alpha}_{4}^{3} {\alpha}_{5}+82 {\alpha}_{3}^{4} {\alpha}_{4}^{3} {\alpha}_{5}+10 {\alpha}_{1}^{3} {\alpha}_{4}^{4} {\alpha}_{5}+80 {\alpha}_{1}^{2} {\alpha}_{2} {\alpha}_{4}^{4} {\alpha}_{5}+180 {\alpha}_{1} {\alpha}_{2}^{2} {\alpha}_{4}^{4} {\alpha}_{5}+120 {\alpha}_{2}^{3} {\alpha}_{4}^{4} {\alpha}_{5}+70 {\alpha}_{1}^{2} {\alpha}_{3} {\alpha}_{4}^{4} {\alpha}_{5}+300 {\alpha}_{1} {\alpha}_{2} {\alpha}_{3} {\alpha}_{4}^{4} {\alpha}_{5}+300 {\alpha}_{2}^{2} {\alpha}_{3} {\alpha}_{4}^{4} {\alpha}_{5}+110 {\alpha}_{1} {\alpha}_{3}^{2} {\alpha}_{4}^{4} {\alpha}_{5}+220 {\alpha}_{2} {\alpha}_{3}^{2} {\alpha}_{4}^{4} {\alpha}_{5}+50 {\alpha}_{3}^{3} {\alpha}_{4}^{4} {\alpha}_{5}+12 {\alpha}_{1}^{2} {\alpha}_{4}^{5} {\alpha}_{5}+48 {\alpha}_{1} {\alpha}_{2} {\alpha}_{4}^{5} {\alpha}_{5}+48 {\alpha}_{2}^{2} {\alpha}_{4}^{5} {\alpha}_{5}+24 {\alpha}_{1} {\alpha}_{3} {\alpha}_{4}^{5} {\alpha}_{5}+48 {\alpha}_{2} {\alpha}_{3} {\alpha}_{4}^{5} {\alpha}_{5}+12 {\alpha}_{3}^{2} {\alpha}_{4}^{5} {\alpha}_{5}+{\alpha}_{1}^{4} {\alpha}_{2}^{2} {\alpha}_{5}^{2}+6 {\alpha}_{1}^{3} {\alpha}_{2}^{3} {\alpha}_{5}^{2}+13 {\alpha}_{1}^{2} {\alpha}_{2}^{4} {\alpha}_{5}^{2}+12 {\alpha}_{1} {\alpha}_{2}^{5} {\alpha}_{5}^{2}+4 {\alpha}_{2}^{6} {\alpha}_{5}^{2}+2 {\alpha}_{1}^{4} {\alpha}_{2} {\alpha}_{3} {\alpha}_{5}^{2}+20 {\alpha}_{1}^{3} {\alpha}_{2}^{2} {\alpha}_{3} {\alpha}_{5}^{2}+60 {\alpha}_{1}^{2} {\alpha}_{2}^{3} {\alpha}_{3} {\alpha}_{5}^{2}+70 {\alpha}_{1} {\alpha}_{2}^{4} {\alpha}_{3} {\alpha}_{5}^{2}+28 {\alpha}_{2}^{5} {\alpha}_{3} {\alpha}_{5}^{2}+{\alpha}_{1}^{4} {\alpha}_{3}^{2} {\alpha}_{5}^{2}+22 {\alpha}_{1}^{3} {\alpha}_{2} {\alpha}_{3}^{2} {\alpha}_{5}^{2}+100 {\alpha}_{1}^{2} {\alpha}_{2}^{2} {\alpha}_{3}^{2} {\alpha}_{5}^{2}+156 {\alpha}_{1} {\alpha}_{2}^{3} {\alpha}_{3}^{2} {\alpha}_{5}^{2}+78 {\alpha}_{2}^{4} {\alpha}_{3}^{2} {\alpha}_{5}^{2}+8 {\alpha}_{1}^{3} {\alpha}_{3}^{3} {\alpha}_{5}^{2}+72 {\alpha}_{1}^{2} {\alpha}_{2} {\alpha}_{3}^{3} {\alpha}_{5}^{2}+168 {\alpha}_{1} {\alpha}_{2}^{2} {\alpha}_{3}^{3} {\alpha}_{5}^{2}+112 {\alpha}_{2}^{3} {\alpha}_{3}^{3} {\alpha}_{5}^{2}+19 {\alpha}_{1}^{2} {\alpha}_{3}^{4} {\alpha}_{5}^{2}+88 {\alpha}_{1} {\alpha}_{2} {\alpha}_{3}^{4} {\alpha}_{5}^{2}+88 {\alpha}_{2}^{2} {\alpha}_{3}^{4} {\alpha}_{5}^{2}+18 {\alpha}_{1} {\alpha}_{3}^{5} {\alpha}_{5}^{2}+36 {\alpha}_{2} {\alpha}_{3}^{5} {\alpha}_{5}^{2}+6 {\alpha}_{3}^{6} {\alpha}_{5}^{2}+2 {\alpha}_{1}^{4} {\alpha}_{2} {\alpha}_{4} {\alpha}_{5}^{2}+22 {\alpha}_{1}^{3} {\alpha}_{2}^{2} {\alpha}_{4} {\alpha}_{5}^{2}+68 {\alpha}_{1}^{2} {\alpha}_{2}^{3} {\alpha}_{4} {\alpha}_{5}^{2}+80 {\alpha}_{1} {\alpha}_{2}^{4} {\alpha}_{4} {\alpha}_{5}^{2}+32 {\alpha}_{2}^{5} {\alpha}_{4} {\alpha}_{5}^{2}+2 {\alpha}_{1}^{4} {\alpha}_{3} {\alpha}_{4} {\alpha}_{5}^{2}+46 {\alpha}_{1}^{3} {\alpha}_{2} {\alpha}_{3} {\alpha}_{4} {\alpha}_{5}^{2}+212 {\alpha}_{1}^{2} {\alpha}_{2}^{2} {\alpha}_{3} {\alpha}_{4} {\alpha}_{5}^{2}+332 {\alpha}_{1} {\alpha}_{2}^{3} {\alpha}_{3} {\alpha}_{4} {\alpha}_{5}^{2}+166 {\alpha}_{2}^{4} {\alpha}_{3} {\alpha}_{4} {\alpha}_{5}^{2}+24 {\alpha}_{1}^{3} {\alpha}_{3}^{2} {\alpha}_{4} {\alpha}_{5}^{2}+216 {\alpha}_{1}^{2} {\alpha}_{2} {\alpha}_{3}^{2} {\alpha}_{4} {\alpha}_{5}^{2}+504 {\alpha}_{1} {\alpha}_{2}^{2} {\alpha}_{3}^{2} {\alpha}_{4} {\alpha}_{5}^{2}+336 {\alpha}_{2}^{3} {\alpha}_{3}^{2} {\alpha}_{4} {\alpha}_{5}^{2}+72 {\alpha}_{1}^{2} {\alpha}_{3}^{3} {\alpha}_{4} {\alpha}_{5}^{2}+332 {\alpha}_{1} {\alpha}_{2} {\alpha}_{3}^{3} {\alpha}_{4} {\alpha}_{5}^{2}+332 {\alpha}_{2}^{2} {\alpha}_{3}^{3} {\alpha}_{4} {\alpha}_{5}^{2}+80 {\alpha}_{1} {\alpha}_{3}^{4} {\alpha}_{4} {\alpha}_{5}^{2}+160 {\alpha}_{2} {\alpha}_{3}^{4} {\alpha}_{4} {\alpha}_{5}^{2}+30 {\alpha}_{3}^{5} {\alpha}_{4} {\alpha}_{5}^{2}+{\alpha}_{1}^{4} {\alpha}_{4}^{2} {\alpha}_{5}^{2}+24 {\alpha}_{1}^{3} {\alpha}_{2} {\alpha}_{4}^{2} {\alpha}_{5}^{2}+110 {\alpha}_{1}^{2} {\alpha}_{2}^{2} {\alpha}_{4}^{2} {\alpha}_{5}^{2}+172 {\alpha}_{1} {\alpha}_{2}^{3} {\alpha}_{4}^{2} {\alpha}_{5}^{2}+86 {\alpha}_{2}^{4} {\alpha}_{4}^{2} {\alpha}_{5}^{2}+24 {\alpha}_{1}^{3} {\alpha}_{3} {\alpha}_{4}^{2} {\alpha}_{5}^{2}+212 {\alpha}_{1}^{2} {\alpha}_{2} {\alpha}_{3} {\alpha}_{4}^{2} {\alpha}_{5}^{2}+492 {\alpha}_{1} {\alpha}_{2}^{2} {\alpha}_{3} {\alpha}_{4}^{2} {\alpha}_{5}^{2}+328 {\alpha}_{2}^{3} {\alpha}_{3} {\alpha}_{4}^{2} {\alpha}_{5}^{2}+100 {\alpha}_{1}^{2} {\alpha}_{3}^{2} {\alpha}_{4}^{2} {\alpha}_{5}^{2}+452 {\alpha}_{1} {\alpha}_{2} {\alpha}_{3}^{2} {\alpha}_{4}^{2} {\alpha}_{5}^{2}+452 {\alpha}_{2}^{2} {\alpha}_{3}^{2} {\alpha}_{4}^{2} {\alpha}_{5}^{2}+132 {\alpha}_{1} {\alpha}_{3}^{3} {\alpha}_{4}^{2} {\alpha}_{5}^{2}+264 {\alpha}_{2} {\alpha}_{3}^{3} {\alpha}_{4}^{2} {\alpha}_{5}^{2}+55 {\alpha}_{3}^{4} {\alpha}_{4}^{2} {\alpha}_{5}^{2}+8 {\alpha}_{1}^{3} {\alpha}_{4}^{3} {\alpha}_{5}^{2}+68 {\alpha}_{1}^{2} {\alpha}_{2} {\alpha}_{4}^{3} {\alpha}_{5}^{2}+156 {\alpha}_{1} {\alpha}_{2}^{2} {\alpha}_{4}^{3} {\alpha}_{5}^{2}+104 {\alpha}_{2}^{3} {\alpha}_{4}^{3} {\alpha}_{5}^{2}+60 {\alpha}_{1}^{2} {\alpha}_{3} {\alpha}_{4}^{3} {\alpha}_{5}^{2}+260 {\alpha}_{1} {\alpha}_{2} {\alpha}_{3} {\alpha}_{4}^{3} {\alpha}_{5}^{2}+260 {\alpha}_{2}^{2} {\alpha}_{3} {\alpha}_{4}^{3} {\alpha}_{5}^{2}+96 {\alpha}_{1} {\alpha}_{3}^{2} {\alpha}_{4}^{3} {\alpha}_{5}^{2}+192 {\alpha}_{2} {\alpha}_{3}^{2} {\alpha}_{4}^{3} {\alpha}_{5}^{2}+44 {\alpha}_{3}^{3} {\alpha}_{4}^{3} {\alpha}_{5}^{2}+13 {\alpha}_{1}^{2} {\alpha}_{4}^{4} {\alpha}_{5}^{2}+52 {\alpha}_{1} {\alpha}_{2} {\alpha}_{4}^{4} {\alpha}_{5}^{2}+52 {\alpha}_{2}^{2} {\alpha}_{4}^{4} {\alpha}_{5}^{2}+26 {\alpha}_{1} {\alpha}_{3} {\alpha}_{4}^{4} {\alpha}_{5}^{2}+52 {\alpha}_{2} {\alpha}_{3} {\alpha}_{4}^{4} {\alpha}_{5}^{2}+13 {\alpha}_{3}^{2} {\alpha}_{4}^{4} {\alpha}_{5}^{2}+2 {\alpha}_{1}^{3} {\alpha}_{2}^{2} {\alpha}_{5}^{3}+8 {\alpha}_{1}^{2} {\alpha}_{2}^{3} {\alpha}_{5}^{3}+10 {\alpha}_{1} {\alpha}_{2}^{4} {\alpha}_{5}^{3}+4 {\alpha}_{2}^{5} {\alpha}_{5}^{3}+4 {\alpha}_{1}^{3} {\alpha}_{2} {\alpha}_{3} {\alpha}_{5}^{3}+24 {\alpha}_{1}^{2} {\alpha}_{2}^{2} {\alpha}_{3} {\alpha}_{5}^{3}+40 {\alpha}_{1} {\alpha}_{2}^{3} {\alpha}_{3} {\alpha}_{5}^{3}+20 {\alpha}_{2}^{4} {\alpha}_{3} {\alpha}_{5}^{3}+2 {\alpha}_{1}^{3} {\alpha}_{3}^{2} {\alpha}_{5}^{3}+24 {\alpha}_{1}^{2} {\alpha}_{2} {\alpha}_{3}^{2} {\alpha}_{5}^{3}+60 {\alpha}_{1} {\alpha}_{2}^{2} {\alpha}_{3}^{2} {\alpha}_{5}^{3}+40 {\alpha}_{2}^{3} {\alpha}_{3}^{2} {\alpha}_{5}^{3}+8 {\alpha}_{1}^{2} {\alpha}_{3}^{3} {\alpha}_{5}^{3}+40 {\alpha}_{1} {\alpha}_{2} {\alpha}_{3}^{3} {\alpha}_{5}^{3}+40 {\alpha}_{2}^{2} {\alpha}_{3}^{3} {\alpha}_{5}^{3}+10 {\alpha}_{1} {\alpha}_{3}^{4} {\alpha}_{5}^{3}+20 {\alpha}_{2} {\alpha}_{3}^{4} {\alpha}_{5}^{3}+4 {\alpha}_{3}^{5} {\alpha}_{5}^{3}+4 {\alpha}_{1}^{3} {\alpha}_{2} {\alpha}_{4} {\alpha}_{5}^{3}+24 {\alpha}_{1}^{2} {\alpha}_{2}^{2} {\alpha}_{4} {\alpha}_{5}^{3}+40 {\alpha}_{1} {\alpha}_{2}^{3} {\alpha}_{4} {\alpha}_{5}^{3}+20 {\alpha}_{2}^{4} {\alpha}_{4} {\alpha}_{5}^{3}+4 {\alpha}_{1}^{3} {\alpha}_{3} {\alpha}_{4} {\alpha}_{5}^{3}+46 {\alpha}_{1}^{2} {\alpha}_{2} {\alpha}_{3} {\alpha}_{4} {\alpha}_{5}^{3}+114 {\alpha}_{1} {\alpha}_{2}^{2} {\alpha}_{3} {\alpha}_{4} {\alpha}_{5}^{3}+76 {\alpha}_{2}^{3} {\alpha}_{3} {\alpha}_{4} {\alpha}_{5}^{3}+22 {\alpha}_{1}^{2} {\alpha}_{3}^{2} {\alpha}_{4} {\alpha}_{5}^{3}+106 {\alpha}_{1} {\alpha}_{2} {\alpha}_{3}^{2} {\alpha}_{4} {\alpha}_{5}^{3}+106 {\alpha}_{2}^{2} {\alpha}_{3}^{2} {\alpha}_{4} {\alpha}_{5}^{3}+32 {\alpha}_{1} {\alpha}_{3}^{3} {\alpha}_{4} {\alpha}_{5}^{3}+64 {\alpha}_{2} {\alpha}_{3}^{3} {\alpha}_{4} {\alpha}_{5}^{3}+14 {\alpha}_{3}^{4} {\alpha}_{4} {\alpha}_{5}^{3}+2 {\alpha}_{1}^{3} {\alpha}_{4}^{2} {\alpha}_{5}^{3}+22 {\alpha}_{1}^{2} {\alpha}_{2} {\alpha}_{4}^{2} {\alpha}_{5}^{3}+54 {\alpha}_{1} {\alpha}_{2}^{2} {\alpha}_{4}^{2} {\alpha}_{5}^{3}+36 {\alpha}_{2}^{3} {\alpha}_{4}^{2} {\alpha}_{5}^{3}+20 {\alpha}_{1}^{2} {\alpha}_{3} {\alpha}_{4}^{2} {\alpha}_{5}^{3}+90 {\alpha}_{1} {\alpha}_{2} {\alpha}_{3} {\alpha}_{4}^{2} {\alpha}_{5}^{3}+90 {\alpha}_{2}^{2} {\alpha}_{3} {\alpha}_{4}^{2} {\alpha}_{5}^{3}+34 {\alpha}_{1} {\alpha}_{3}^{2} {\alpha}_{4}^{2} {\alpha}_{5}^{3}+68 {\alpha}_{2} {\alpha}_{3}^{2} {\alpha}_{4}^{2} {\alpha}_{5}^{3}+16 {\alpha}_{3}^{3} {\alpha}_{4}^{2} {\alpha}_{5}^{3}+6 {\alpha}_{1}^{2} {\alpha}_{4}^{3} {\alpha}_{5}^{3}+24 {\alpha}_{1} {\alpha}_{2} {\alpha}_{4}^{3} {\alpha}_{5}^{3}+24 {\alpha}_{2}^{2} {\alpha}_{4}^{3} {\alpha}_{5}^{3}+12 {\alpha}_{1} {\alpha}_{3} {\alpha}_{4}^{3} {\alpha}_{5}^{3}+24 {\alpha}_{2} {\alpha}_{3} {\alpha}_{4}^{3} {\alpha}_{5}^{3}+6 {\alpha}_{3}^{2} {\alpha}_{4}^{3} {\alpha}_{5}^{3}+{\alpha}_{1}^{2} {\alpha}_{2}^{2} {\alpha}_{5}^{4}+2 {\alpha}_{1} {\alpha}_{2}^{3} {\alpha}_{5}^{4}+{\alpha}_{2}^{4} {\alpha}_{5}^{4}+2 {\alpha}_{1}^{2} {\alpha}_{2} {\alpha}_{3} {\alpha}_{5}^{4}+6 {\alpha}_{1} {\alpha}_{2}^{2} {\alpha}_{3} {\alpha}_{5}^{4}+4 {\alpha}_{2}^{3} {\alpha}_{3} {\alpha}_{5}^{4}+{\alpha}_{1}^{2} {\alpha}_{3}^{2} {\alpha}_{5}^{4}+6 {\alpha}_{1} {\alpha}_{2} {\alpha}_{3}^{2} {\alpha}_{5}^{4}+6 {\alpha}_{2}^{2} {\alpha}_{3}^{2} {\alpha}_{5}^{4}+2 {\alpha}_{1} {\alpha}_{3}^{3} {\alpha}_{5}^{4}+4 {\alpha}_{2} {\alpha}_{3}^{3} {\alpha}_{5}^{4}+{\alpha}_{3}^{4} {\alpha}_{5}^{4}+2 {\alpha}_{1}^{2} {\alpha}_{2} {\alpha}_{4} {\alpha}_{5}^{4}+6 {\alpha}_{1} {\alpha}_{2}^{2} {\alpha}_{4} {\alpha}_{5}^{4}+4 {\alpha}_{2}^{3} {\alpha}_{4} {\alpha}_{5}^{4}+2 {\alpha}_{1}^{2} {\alpha}_{3} {\alpha}_{4} {\alpha}_{5}^{4}+10 {\alpha}_{1} {\alpha}_{2} {\alpha}_{3} {\alpha}_{4} {\alpha}_{5}^{4}+10 {\alpha}_{2}^{2} {\alpha}_{3} {\alpha}_{4} {\alpha}_{5}^{4}+4 {\alpha}_{1} {\alpha}_{3}^{2} {\alpha}_{4} {\alpha}_{5}^{4}+8 {\alpha}_{2} {\alpha}_{3}^{2} {\alpha}_{4} {\alpha}_{5}^{4}+2 {\alpha}_{3}^{3} {\alpha}_{4} {\alpha}_{5}^{4}+{\alpha}_{1}^{2} {\alpha}_{4}^{2} {\alpha}_{5}^{4}+4 {\alpha}_{1} {\alpha}_{2} {\alpha}_{4}^{2} {\alpha}_{5}^{4}+4 {\alpha}_{2}^{2} {\alpha}_{4}^{2} {\alpha}_{5}^{4}+2 {\alpha}_{1} {\alpha}_{3} {\alpha}_{4}^{2} {\alpha}_{5}^{4}+4 {\alpha}_{2} {\alpha}_{3} {\alpha}_{4}^{2} {\alpha}_{5}^{4}+{\alpha}_{3}^{2} {\alpha}_{4}^{2} {\alpha}_{5}^{4})$
    \caption{Multidegree of \(X_\tau\)}
    \label{pierre}
\end{figure}

%
\begin{figure}[ht]
    \centering
    \begin{multicols}{2}
        \medskip
        
        \setbox\ltmcbox\vbox{
        \makeatletter\col@number\@ne
        \begin{longtable}{>{$}l<{$}}
        {\Delta}_{1\bar 1\bar 2\bar 3\bar 4}\\
            {\Delta}_{12\bar 2\bar 3\bar 4}\\
            {\Delta}_{13\bar 1\bar 3\bar 4}\\
            {\Delta}_{12\bar 1\bar 3\bar 4}\\
            {\Delta}_{123\bar 3\bar 4}\\
            {\Delta}_{14\bar 1\bar 2\bar 4}\\
            {\Delta}_{13\bar 1\bar 2\bar 4}\\
            {\Delta}_{12\bar 1\bar 2\bar 4}\\
            {\Delta}_{124\bar 2\bar 4}\\
            {\Delta}_{123\bar 2\bar 4}\\
            {\Delta}_{134\bar 1\bar 4}\\
            {\Delta}_{124\bar 1\bar 4}\\
            {\Delta}_{123\bar 1\bar 4}\\
            {\Delta}_{1234\bar 4}\\
            {\Delta}_{15\bar 1\bar 2\bar 3}\\
            {\Delta}_{14\bar 1\bar 2\bar 3}\\
            {\Delta}_{13\bar 1\bar 2\bar 3}\\
            {\Delta}_{12\bar 1\bar 2\bar 3}\\
            {\Delta}_{125\bar 2\bar 3}\\
            {\Delta}_{124\bar 2\bar 3}\\
            {\Delta}_{123\bar 2\bar 3}\\
            {\Delta}_{145\bar 1\bar 2}\\
            {\Delta}_{135\bar 1\bar 2}-{\Delta}_{125\bar 1\bar 3}\\
            {\Delta}_{134\bar 1\bar 2}-{\Delta}_{124\bar 1\bar 3}\\
            {\Delta}_{1245\bar 2}\\
            1\\
            {\Delta}_{125\bar 1\bar 3}^{3}-{\Delta}_{125\bar 1\bar 2}{\Delta}_{135\bar 1\bar 3}\\
            {\Delta}_{134\bar 1\bar 3}{\Delta}_{125\bar 1\bar 3}-{\Delta}_{124\bar 1\bar 3}{\Delta}_{135\bar 1\bar 3}\\
            {\Delta}_{124\bar 1\bar 3}{\Delta}_{125\bar 1\bar 3}-{\Delta}_{124\bar 1\bar 2}{\Delta}_{135\bar 1\bar 3}\\
            {\Delta}_{123\bar 1\bar 3}{\Delta}_{125\bar 1\bar 3}-{\Delta}_{123\bar 1\bar 2}{\Delta}_{135\bar 1\bar 3}\\
            {\Delta}_{1235\bar 3}{\Delta}_{125\bar 1\bar 3}-{\Delta}_{1235\bar 2}{\Delta}_{135\bar 1\bar 3}\\
            {\Delta}_{1234\bar 3}{\Delta}_{125\bar 1\bar 3}-{\Delta}_{1234\bar 2}{\Delta}_{135\bar 1\bar 3}\\
            {\Delta}_{1345\bar 1}{\Delta}_{125\bar 1\bar 3}-{\Delta}_{1245\bar 1}{\Delta}_{135\bar 1\bar 3}\\
            {\Delta}_{1235\bar 3}{\Delta}_{134\bar 1\bar 3}-{\Delta}_{1234\bar 3}{\Delta}_{135\bar 1\bar 3}\\
            {\Delta}_{125\bar 1\bar 2}{\Delta}_{134\bar 1\bar 3}-{\Delta}_{124\bar 1\bar 2}{\Delta}_{135\bar 1\bar 3}\\
            {\Delta}_{1235\bar 2}{\Delta}_{134\bar 1\bar 3}-{\Delta}_{1234\bar 2}{\Delta}_{135\bar 1\bar 3}\\
            {\Delta}_{124\bar 1\bar 3}^{3}-{\Delta}_{124\bar 1\bar 2}{\Delta}_{134\bar 1\bar 3}\\
            {\Delta}_{123\bar 1\bar 3}{\Delta}_{124\bar 1\bar 3}-{\Delta}_{123\bar 1\bar 2}{\Delta}_{134\bar 1\bar 3}\\
            {\Delta}_{1235\bar 3}{\Delta}_{124\bar 1\bar 3}-{\Delta}_{1234\bar 2}{\Delta}_{135\bar 1\bar 3}\\
            {\Delta}_{1234\bar 3}{\Delta}_{124\bar 1\bar 3}-{\Delta}_{1234\bar 2}{\Delta}_{134\bar 1\bar 3}\\
            {\Delta}_{125\bar 1\bar 2}{\Delta}_{124\bar 1\bar 3}-{\Delta}_{124\bar 1\bar 2}{\Delta}_{125\bar 1\bar 3}\\
            {\Delta}_{1235\bar 2}{\Delta}_{124\bar 1\bar 3}-{\Delta}_{1234\bar 2}{\Delta}_{125\bar 1\bar 3}\\
            {\Delta}_{1345\bar 1}{\Delta}_{124\bar 1\bar 3}-{\Delta}_{1245\bar 1}{\Delta}_{134\bar 1\bar 3}\\
            {\Delta}_{125\bar 1\bar 2}{\Delta}_{123\bar 1\bar 3}-{\Delta}_{123\bar 1\bar 2}{\Delta}_{125\bar 1\bar 3}\\
            {\Delta}_{124\bar 1\bar 2}{\Delta}_{123\bar 1\bar 3}-{\Delta}_{123\bar 1\bar 2}{\Delta}_{124\bar 1\bar 3}\\
            {\Delta}_{1345\bar 1}{\Delta}_{123\bar 1\bar 3}-{\Delta}_{1235\bar 1}{\Delta}_{134\bar 1\bar 3}+{\Delta}_{1234\bar 1}{\Delta}_{135\bar 1\bar 3}\\
            {\Delta}_{1245\bar 1}{\Delta}_{123\bar 1\bar 3}-{\Delta}_{1235\bar 1}{\Delta}_{124\bar 1\bar 3}+{\Delta}_{1234\bar 1}{\Delta}_{125\bar 1\bar 3}\\
            {\Delta}_{125\bar 1\bar 2}{\Delta}_{1235\bar 3}-{\Delta}_{1235\bar 2}{\Delta}_{125\bar 1\bar 3}\\
            {\Delta}_{124\bar 1\bar 2}{\Delta}_{1235\bar 3}-{\Delta}_{1234\bar 2}{\Delta}_{125\bar 1\bar 3}\\
            {\Delta}_{123\bar 1\bar 2}{\Delta}_{1235\bar 3}-{\Delta}_{1235\bar 2}{\Delta}_{123\bar 1\bar 3}\\
            {\Delta}_{1345\bar 1}{\Delta}_{1235\bar 3}+{\Delta}_{12345}{\Delta}_{135\bar 1\bar 3}\\
            {\Delta}_{1245\bar 1}{\Delta}_{1235\bar 3}+{\Delta}_{12345}{\Delta}_{125\bar 1\bar 3}\\
            {\Delta}_{125\bar 1\bar 2}{\Delta}_{1234\bar 3}-{\Delta}_{1234\bar 2}{\Delta}_{125\bar 1\bar 3}\\
            {\Delta}_{124\bar 1\bar 2}{\Delta}_{1234\bar 3}-{\Delta}_{1234\bar 2}{\Delta}_{124\bar 1\bar 3}\\
            {\Delta}_{123\bar 1\bar 2}{\Delta}_{1234\bar 3}-{\Delta}_{1234\bar 2}{\Delta}_{123\bar 1\bar 3}\\
            {\Delta}_{1235\bar 2}{\Delta}_{1234\bar 3}-{\Delta}_{1234\bar 2}{\Delta}_{1235\bar 3}\\
            {\Delta}_{1345\bar 1}{\Delta}_{1234\bar 3}+{\Delta}_{12345}{\Delta}_{134\bar 1\bar 3}\\
            {\Delta}_{1245\bar 1}{\Delta}_{1234\bar 3}+{\Delta}_{12345}{\Delta}_{124\bar 1\bar 3}\\
            {\Delta}_{1235\bar 1}{\Delta}_{1234\bar 3}-{\Delta}_{1234\bar 1}{\Delta}_{1235\bar 3}+{\Delta}_{12345}{\Delta}_{123\bar 1\bar 3}\\
            {\Delta}_{1345\bar 1}{\Delta}_{125\bar 1\bar 2}-{\Delta}_{1245\bar 1}{\Delta}_{125\bar 1\bar 3}\\
            {\Delta}_{1235\bar 2}{\Delta}_{124\bar 1\bar 2}-{\Delta}_{1234\bar 2}{\Delta}_{125\bar 1\bar 2}\\
            {\Delta}_{1345\bar 1}{\Delta}_{124\bar 1\bar 2}-{\Delta}_{1245\bar 1}{\Delta}_{124\bar 1\bar 3}\\
            {\Delta}_{1345\bar 1}{\Delta}_{123\bar 1\bar 2}-{\Delta}_{1235\bar 1}{\Delta}_{124\bar 1\bar 3}+{\Delta}_{1234\bar 1}{\Delta}_{125\bar 1\bar 3}\\
            {\Delta}_{1245\bar 1}{\Delta}_{123\bar 1\bar 2}-{\Delta}_{1235\bar 1}{\Delta}_{124\bar 1\bar 2}+{\Delta}_{1234\bar 1}{\Delta}_{125\bar 1\bar 2}\\
            {\Delta}_{1345\bar 1}{\Delta}_{1235\bar 2}+{\Delta}_{12345}{\Delta}_{125\bar 1\bar 3}\\
            {\Delta}_{1245\bar 1}{\Delta}_{1235\bar 2}+{\Delta}_{12345}{\Delta}_{125\bar 1\bar 2}\\
            {\Delta}_{1345\bar 1}{\Delta}_{1234\bar 2}+{\Delta}_{12345}{\Delta}_{124\bar 1\bar 3}\\
            {\Delta}_{1245\bar 1}{\Delta}_{1234\bar 2}+{\Delta}_{12345}{\Delta}_{124\bar 1\bar 2}\\
            {\Delta}_{1235\bar 1}{\Delta}_{1234\bar 2}-{\Delta}_{1234\bar 1}{\Delta}_{1235\bar 2}+{\Delta}_{12345}{\Delta}_{123\bar 1\bar 2}
        \end{longtable}
        \unskip
        \unpenalty
        \unpenalty}
        
        \unvbox\ltmcbox
        
        \medskip
        \end{multicols}
    \caption{Ideal of \(Z_\tau \cap \{\Delta_{12345} = 1\}\)}
    \label{a4pluck}
\end{figure}
\section*{Index of Notation}
%
%
{
\singlespacing
\begin{longtable}{>{\(}l<{\)}l}
    G & simple simply-connected complex algebraic group \\
    \g & the Lie algebra of \(G\) \\ 
    G^\vee & the Langlands dual group of \(G\) \\
    U,\, B,\, T & the usual associated suspects \(B = UT \subset G\) \\
    \fu,\, \fb,\, \t & their Lie algebras \\
    U^\vee,\, B^\vee,\, T^\vee & their Langlands dual groups \\
    P = X^\bullet(T) & weight lattice \\
    P^\vee = X_\bullet(T) & coweight lattice, \(P^\vee = X^\bullet(T^\vee) \) \\ 
    \t_\RR^\ast = P\otimes\RR & real weight lattice \\
    \t_\RR = P^\vee \otimes\RR & real coweight lattice \\
    w_0 \in W & the longest element in the Weyl group of \(G\) \\
    \ell = \len w_0 & length of any reduced word for \(w_0\) \\
    (I,H, H \overset{\varepsilon}{\to} \{\pm 1\}) & the doubled Dynkin quiver of \(G\) \\
    \{\alpha_i\}_I,\,\{\alpha_i^\vee\}_I & basis of simple roots for \(\t^\ast\), simple coroots for \(\t\) \\
    a_{ij} = \langle \alpha_i^\vee,\alpha_j\rangle = \alpha_j(\alpha_i^\vee) & \((i,j)\)th entry of Cartan matrix of \(\g\) \\ 
    s_i = s_{\alpha_i} & simple reflections in \(W\) acting by \(\lambda\mapsto \lambda - \langle\lambda,\alpha_i\rangle\alpha_i\) on \(P\) \\
    e_i = e_{\alpha_i} & Chevalley generators of \(\cU(\fu)\), root vectors weighted by \(\Pi\) \\
    \uvi = (i_1,\dots,i_{\ell}) & reduced expression for \(w_0\) \\
    w_k^{\uvi} = s_{i_1}\cdots s_{i_k} & subword \\ 
    \le_{\uvi} & associated convex order \\
    \Delta \supset \Delta_+ = \{\beta_k^{\uvi}\}_1^\ell & roots, positive roots; ordered by \(\le_{\uvi} \) with \(\beta_k^{\uvi} = w_{k-1}^{\uvi}\alpha_k \) \\ 
    \alpha_{i,j} = \alpha_i + \cdots + \alpha_j & shorthand for positive roots in type \(A\) \\ 
    Q_+\subset Q & ``positive'' root cone in root lattice \\
    \nu = \sum \nu_i & if \(\nu = \sum\nu_i \alpha_i\) is positive \\ 
    \Seq(\nu) & all \(\vi = (i_1,\dots,i_p)\) such that \(\sum \alpha_{i_k} = \nu\)\\
    \nu_\bullet = (\nu_w)_W & collection of vectors in \(P\) \\
    \{\omega_i\}_{I}, \{\omega_i^\vee\}_{I} & fundamental weights, fundamental coweights \\
    \mathring X_\tau \subset X_\tau & open dense subset of generalized orbital variety \\
    \mathring Z_\tau \subset Z_\tau & open dense subset of associated MV cycle, \(\mathring Z_\tau = \phi(X_\tau)\) \\
    \cZ(\infty) = \bigcup_{\nu\in Q_+} \cZ(\infty)_{\nu} & stable MV cycles \\ 
    \cZ(\lambda) = \bigcup_{\mu\in P_+} \cZ(\lambda)_\mu & MV cycles of type \(\lambda\) \\
    B(\infty) & Lusztig's dual canonical basis of \(\CC[U]\) \\ 
    c_Y \in B(\Lambda) & dual semicanonical basis of \(\CC[U]\) \\ 
    b_Z \in B(\Gr) & MV basis of \(\CC[U]\) \\ 
    v_\lambda\in L(\lambda) & highest weight vector in irreducible highest weight rep \\
    \Psi = \sum \Psi_\lambda & Berenstein--Kazhdan map \(\bigcup L(\lambda) \to \CC[U]\) \\
    \cO\subset \cK & rings of formal power series, Laurent series \\
    \Gr & affine Grassmannian\\
    \overline{\OO}_\lambda\cap\TT_\mu \overset{\Phi}{\to} \overline{\Gr^\lambda}\cap\Gr_\mu & Mirkovi\'c--Vybornov isomorphism \\ 
    \phi = \Phi\big|_{\overline{\OO}_\lambda\cap\TT_\mu\cap\n} & restricted Mirkovi\'c--Vybornov isomorphism \\ 
    \psi & moment map defining \(\Lambda\) \\ 
    \cT(\lambda)_\mu & semistandard Young tableaux of shape \(\lambda\) and weight \(\mu\)\\ 
    \overline{\OO}_\lambda\cap\TT_\mu\cap\n\overset{\tau}{\to}\cT(\lambda)_\mu & labeling of generalized orbital varieties by semi-standard Young tableaux \\
    \cA & the preprojective algebra of \((I,H,\varepsilon)\) \\
    \Lambda & Lusztig's nilpotent variety, elements are modules \(M\) for \(\cA\)\\
    F_\vi( M ) & variety of composition series of type \(\vi\) \\
    F_n(M) \supset F_{n,\mu}(M) & \(n+1\) step flags in; of \(\dimvec =\mu\) \\ 
    H^0 & sections \\
    \dimvec & dimension vector of a \(\cA\)-module in \(Q_+\)\\ 
    \scO,\, \scL & structure sheaf, line bundle \\ 
    n^{\uvi}(\phantom{0}) = (n^{\uvi}(\phantom{0})_\beta)_{\Delta_+} & \(\uvi\)-Lusztig datum \\
    \Seq(\nu) & ordered partitions of \(\nu\) \\
    \vi,\vk,\vj & elements of \(\Seq(\nu)\) \\
    \cPP & measures on \(\t_\RR^\ast\) \\
    \shuffle & shuffle product \\
    \delta_{\Delta^p} & Lebesgue measure on \(p\)-simplex \\
    D_\vi & pushfroward of Lebesgue measure on \(p\)-simplex \\
    \barD_\vi & the ``leading'' coefficient on its Fourier transform \\
\end{longtable}
}
\addcontentsline{toc}{chapter}{Bibliography}
\bibliographystyle{alpha}
\bibliography{master}
\end{document}